\documentclass[a4paper]{article}

\usepackage{amssymb}
\usepackage{latexsym}
\usepackage{amsfonts}
\usepackage{mathrsfs}
\usepackage{amscd}
\usepackage{amsmath}
\usepackage{amssymb}
\usepackage{wasysym}
\usepackage{bm}
\newtheorem{Th}{Theorem}[section]
\newtheorem{Prop}{Proposition}[section]
\newtheorem{Co}{Corollary}[section]
\newtheorem{Lma}{Lemma}[section]

\newtheorem{Rm}{Remark}

\newcommand{\be}{\begin{equation}}
\newcommand{\ee}{\end{equation}}
\newcommand{\bes}{\begin{equation*}}
\newcommand{\ees}{\end{equation*}}
\newcommand{\R}{\mathbb{R}}
\newcommand{\N}{\mathbb{N}}
\newcommand{\C}{\mathbb{C}}

\newcommand\res{\mathop{\hbox{\vrule height 7pt width .5pt depth 0pt
\vrule height .5pt width 6pt depth 0pt}}\nolimits}

\def\theequation{\thesection.\arabic{equation}}
\def\theTh{\Roman{section}.\arabic{Th}}
\def\theProp{\Roman{section}.\arabic{Prop}}
\def\theCo{\Roman{section}.\arabic{Co}}

\def\theRm{\Roman{section}.\arabic{Rm}}
\newcommand{\reset}{\setcounter{equation}{0}\setcounter{Th}{0}\setcounter{Prop}{0}\setcounter{Co}{0}\setcounter{Lma}{0}\setcounter{Rm}{0}}

\def\al{\alpha}
\def\la{\lambda}
\def\eps{\varepsilon}

\def\Om{\Omega}
\def\om{\omega}

\def\pro{\pi_{\vec{n}}}
\def\bn{\vec{n}}
\def\bna{\vec{n}_\al}

\def\bex{\bAe_1}
\def\bey{\bAe_2}

\def\bei{\bAe_i}
\def\bej{\bAe_j}
\def\bek{\bAe_k}

\def\px{\partial_{x_1}}
\def\py{\partial_{x_2}}
\def\pj{\partial_{x_j}}

\def\pk{\partial_{x_k}}

\def\bAe{\vec{e}}
\def\bH{\vec{H}}

\def\bC{\vec{C}}
\def\bA{\vec{A}}
\def\bAe{\vec{e}}

\def\bf{\vec{f}}
\def\bF{\vec{F}}

\def\bI{{I}}
\def\bL{\vec{L}}
\def\bG{\vec{G}}
\def\bR{\vec{R}}
\def\bV{\vec{V}}

\def\bX{\vec{X}}
\def\bPe{\vec{P}}

\def\bB{\vec{B}}

\def\bc{\vec{c}}

\def\bp{\vec{\Phi}}
\def\bP{\vec{\Phi}}

\def\bT{\vec{T}}
\def\bK{{K}}
\def\bQ{\vec{Q}}

\def\bet{\beta}
\def\bul{\bullet}
\def\di{{D}^2}
\def\res{\mathop{\hbox{\vrule height 7pt width .5pt 
depth 0pt\vrule height .5pt width 6pt depth 0pt}}\nolimits}

\begin{document}

%%%%%%%%%%%%%%%%%%%
%%%%%%%%%%%%%%%%%%%
%%%%%%%%%%%%%%%%%%%
%%%%%%%%%%%%%%%%%%%
%%%%%%%%%%%%%%%%%%%
%%%%%%%%%%%%%%%%%%%
%%%%%%%%%%%%%%%%%%%
%%%%%%%%%%%%%%%%%%%
%%%%%%%%%%%%%%%%%%%
%%%%%%%%%%%%%%%%%%%
%%%%%%%%%%%%%%%%%%%
%%%%%%%%%%%%%%%%%%%
%%%%%%%%%%%%%%%%%%%

\reset

\title{Asymptotic Analysis of
Branched Willmore Surfaces}
\author{Yann Bernard\footnote{Mathematisches Institut, Albert-Ludwigs-Universit\"at, 79004 Freiburg, Germany. }\:\:,\:Tristan Rivi\`ere\footnote{Department of Mathematics, ETH Zentrum,
8093 Z\"urich, Switzerland.}}
\date{ }
\maketitle
\noindent
$\textbf{Abstract:}$ {\it We consider a closed Willmore surface properly immersed in ${\R}^{m\ge3}$ with square-integrable second fundamental form, and with one point-singularity of finite arbitrary integer order. Using the ``conservative" reformulation of the Willmore equation introduced in \cite{Ri1}, we show that, in an appropriate conformal parametrization, the gradient of the Gauss map of the immersion has bounded mean oscillations if the singularity has order one, and is bounded if the order is at least two. We develop around the singular point local asymptotic expansions for the immersion, its first and second derivatives, and for the mean curvature vector. Finally, we exhibit an explicit condition ensuring the removability of the point-singularity.}

\medskip

\noindent{\it Math. Class.} 30C70, 58E15, 58E30, 49Q10, 53A30, 35R01, 35J35, 35J48, 35J50.

\section{Overview}

%\subsection{Introduction}

The Willmore energy of an immersed closed surface $\bp:\Sigma\rightarrow\R^{m\ge3}$ is given by
\be\label{wilen}
W(\bp\,)\;:=\;\int_{\Sigma}|\bH|^2\,d\text{vol}_g\:,
\ee
where $\bH$ denotes the weak mean curvature vector, and $d\text{vol}_g$ is the area form of the metric $g$ induced on $\bp(\Sigma)$ by the canonical Euclidean metric on $\R^m$. Critical points of the Lagrangian $W$ for perturbations of the form $\bp+t\,\vec{\xi}$, where $\vec{\xi}$ is an arbitrary compactly supported smooth map on $\Sigma$ into $\R^m$, are known as {\it Willmore surfaces}. Not only is the Willmore functional invariant under reparametrization, but more importantly, it is invariant under the group of M\"obius transformations of $\R^m\cup\{\infty\}$. This remarkable property prompts the use of the Willmore energy in various fields of science. A survey of the Willmore functional, of its properties, and of the relevant literature is available in \cite{Ri3}. \\

The study of singular points of Willmore immersions is primarily motivated by the fact that sequences of Willmore immersions with uniformly bounded energy converge everywhere except on a finite set of points where the energy concentrates (cf. \cite{BR2} and the references therein). Such point singularities of Willmore surfaces also occur as blow-ups of the Willmore flow (cf. \cite{KS1}). In their seminal paper \cite{KS1}, Ernst Kuwert and Reiner Sch\"atzle initiated the analytical study of point-singularities of Willmore immersions by first considering unit-density singularities in codimension 1. In a second paper \cite{KS2}, the authors studied singularities of higher order, still in codimension 1. Through a different approach, in \cite{Ri1}, the author recovered and extended the results from \cite{KS1} in arbitrary codimension. In the present paper, the original method developed in \cite{Ri1} is led to fruition in the study of point-singularities of arbitrary order in arbitrary codimension. Not only are all aforementioned results recovered, but new ones as well. Our goal is two-fold: understand the regularity of the Gauss map near a point-singularity of arbitrary (finite) integer order in arbitrary codimension, and develop precise asymptotics for the immersion and the mean curvature near that point. We also give an explicit condition ensuring that the point-singularity is removable. \\

Owing to the Gauss-Bonnet theorem, we note that the Willmore energy (\ref{wilen}) may be equivalently expressed as
\bes
W\big(\bP(\Sigma)\big)\,=\,\int_{\Sigma}\,\big|\vec{\mathbb{I}}\big|_g^2\,d\mu_g\;+\;\pi\chi(\Sigma)\:,
\ees
where $\vec{\mathbb{I}}$ is the second fundamental form, and $\chi(\Sigma)$ is the Euler characteristic of $\Sigma$, which is a topological invariant for a closed surface. From the variational point of view, Willmore surfaces are thus critical points of the energy
\bes
\int_{\Sigma}\,\big|\vec{\mathbb{I}}\big|_g^2\,d\text{vol}_g\:.
\ees
It then appears natural to restrict our attention on immersions whose second fundamental forms are locally square-integrable. \\

We assume that the point-singularity lies at the origin, and we localize the problem by considering a map $\bp:\di\rightarrow\R^{m\ge3}$, which is an immersion of $\di\setminus\{0\}$, and satisfying
\begin{itemize}
\item[(i)] $\:\:\:\bp\in C^0(\di)\cap C^\infty(\di\setminus\{0\})\:;$
\item[(ii)] $\:\:\:\mathcal{H}^2\big(\bp(\di)\big)\,<\,\infty\:;$
\item[(iii)] $\displaystyle{\:\:\:\int_{\di}|\vec{\mathbb{I}}|^2_g\,d\text{vol}_g\,<\,\infty}\:.$
\end{itemize}
By a procedure detailed in \cite{KS2}, it is possible to construct a parametrization $\zeta$ of the unit-disk such that $\bp\circ\zeta$ is {\it conformal}. To do so, one first extends $\bp$ to all of $\C\setminus\{0\}$ while keeping a bounded image and the second fundamental form square-integrable. One then shifts so as to have $\bp(0)=\vec{0}$, and inverts about the origin so as to obtain a complete immersion with square-integrable second fundamental form. Calling upon a result of Huber \cite{Hu} (see also \cite{MS} and \cite{To}), one deduces that the image of the immersion is conformally equivalent to $\C$. Inverting yet once more about the origin finally gives the desired conformal immersion\footnote{which degenerates at the origin in a particular way, see (\ref{immas}).}, which we shall abusively continue to denote $\bp$. It has the aforementioned properties (i)-(iii), and moreover,
\bes
\bp(0)\,=\,\vec{0}\qquad\text{and}\qquad\bp(\di)\subset B^m_R(0)\quad\text{for some $\:0<R<\infty$}\:.
\ees
Hence, $\bp\in L^\infty\cap W^{1,2}(\di\setminus\{0\})$. Away from the origin, we define the {\it Gauss map} $\bn$ via
\bes
\bn\;=\;\star\,\dfrac{\px\bp\wedge\py\bp}{|\px\bp\wedge\py\bp|}\:\:,
\ees
where $(x_1,x_2)$ are standard Cartesian coordinates on the unit-disk $\di$, and $\star$ is the Euclidean Hodge-star operator. The immersion $\bp$ is conformal, i.e.
\be\label{confcond}
|\px\bp|\;=\;\text{e}^{\la}\;=\;|\py\bp|\qquad\text{and}\qquad\px\bp\cdot\py\bp\;=\;0\:,
\ee
where $\la$ is the conformal parameter. 
An elementary computation shows that
\be\label{eraser}
d\text{vol}_g\;=\;\text{e}^{2\la}dx\qquad\text{and}\qquad|\nabla\bn|^2\,dx\;=\;\text{e}^{2\la}|\vec{\mathbb{I}}|_g^2\,dx\;=\;|\vec{\mathbb{I}}|^2_g\,d\text{vol}_g\:.
\ee
Hence, by hypothesis, we see that $\bn\in W^{1,2}(\di\setminus\{0\})$. In dimension two, the 2-capacity of isolated points is null, so we actually have $\bp\in W^{1,2}(\di)$ and $\bn\in W^{1,2}(\di)$ (note however that $\bp$ remains a non-degenerate immersion only away from the singularity). Rescaling if necessary, we shall henceforth always assume that 
\be\label{acheumeuneu}
\int_{\di}|\nabla\bn|^2\,dx\;<\;\eps_0\:,
\ee
where the adjustable parameter $\eps_0$ is chosen to fit our various needs (in particular, we will need it to be ``small enough" in Proposition \ref{morreydecay}). \\

For the sake of the following paragraph, we consider a conformal immersion $\bp:\di\rightarrow\R^m$, which is smooth across the unit-disk. We introduce the local coordinates $(x_1,x_2)$ for the flat metric on the unit-disk $\,\di=\big\{x=(x_1,x_2)\in{\R}^2\ ;\ x_1^2+x_2^2<1\big\}$. The operators $\nabla=(\px,\py)$, $\nabla^{\perp}=(-\py,\px)$, $\text{div}=\nabla\cdot\,$, and $\Delta=\nabla\cdot\nabla$  will be understood in these coordinates.
The conformal parameter $\la$ is defined as in (\ref{confcond}). We set
\be\label{lesvec}
\bej\::=\:\text{e}^{-\la}\pj\bp\qquad\text{for}\quad j\,\in\,\{1,2\}\:.
\ee
As $\bp$ is conformal, $\{\bex(x),\bey(x)\}$ forms an orthonormal basis of the tangent space $T_{\bp(x)}\bp(\di)$. Owing to the topology of $\di$, there exists for almost every $x\in\di$ a positively oriented orthonormal basis $\{\bn_1,\ldots,\bn_{m-2}\}$ of the normal space $N_{\bp(x)}\bp(\di)$, such that 
$\{\bex,\bey,\bn_1,\ldots,\bn_{m-2}\}$ forms a basis of $T_{\bp(x)}\R^m$. From the Pl\"ucker embedding, realizing the Grassmannian $Gr_{m-2}(\R^m)$ as a submanifold of the projective space of the $(m-2)^\text{th}$ exterior power $\,\mathbb{P}\big(\bigwedge^{m-2}\R^m\big)$, we can represent the Gauss map as the $(m-2)$-vector $\,\bn=\bigwedge_{\al=1}^{m-2}\bn_\al.$
Via the Hodge operator $\star\,$, we identify vectors and $(m-1)$-vectors in ${\R}^m$, namely:
\bes
%\label{II.0}
\star\,(\bn\wedge \bAe_1)\,=\,\bAe_2\quad,\qquad\star\,(\bn\wedge \bAe_2)\,=\,-\,\bAe_1\quad,\qquad\star\,(\bex\wedge\bey)\,=\,\bn\:.
\ees
In this notation, the second fundamental form $\vec{\mathbb{I}}$, which is a symmetric 2-form on $T_{\bp(x)}\bp(\di)$  into  $N_{\bp(x)}\bp(\di)$, is expressed as
\bes\label{vecb}
\vec{\mathbb{I}}\;=\;\sum_{\al,i,j}\ \text{e}^{-2\la}\,h^\al_{ij}\ \bn_\al\,dx_i\otimes dx_j\;\equiv\;\sum_{\al,i,j}\ h^\al_{ij}\ \bn_\al\,(\bAe_i)^\ast\otimes(\bAe_j)^\ast\:,
\ees
where
\bes
h^\al_{ij}\;=\;-\,\text{e}^{-\la}\,\bei\cdot\pj\bn_\al\:.
\ees\\
The mean curvature vector is
\bes
\bH\:=\:\sum_{\al=1}^{m-2}\,H^\al\,\bn_\al\:=\:\;\frac{1}{2}\,\sum_{\al=1}^{m-2}\,\big(h^\al_{11}+h^\al_{22}\big)\, \bn_\al\:.
\ees
\noindent
The Willmore equation \cite{We} is cast in the form 
\be\label{wil1}
\Delta_\perp\bH\,+\,\sum_{\al,\beta,i,j}\,h^\al_{ij}\,h^\bet_{ij}\,H^\bet\,\bn_\al\:-\:2\,\big|\vec{H}\big|^2\bH\:=\:0\:,
\ee
with
\bes
\Delta_\perp\bH\::=\:\text{e}^{-2\la}\,\pi_{\bn}\,\text{div}\big(\pi_{\bn}(\nabla\bH)\big)\:,
\ees
and $\pro$ is the projection onto the normal space spanned by $\{\bn_\al\}_{\al=1}^{m-2}$. \\

\noindent
The Willmore equation (\ref{wil1}) is a fourth-order nonlinear equation (in the coefficients of the induced metric, which depends on $\bp$). With respect to the coefficients $H^\al$ of the mean curvature vector, it is actually a strongly coupled nonlinear system whose study is particularly challenging. In codimension 1, there is one equation for the scalar curvature ; in higher codimension however, the situation becomes significantly more complicated, and one must seek different techniques to approach the problem. Fortunately, in a conformal parametrization, it is possible\footnote{this procedure requires to choose the normal frame $\{\bn_\al\}$ astutely. See \cite{Ri1} for details.} to recast the system (\ref{wil1}) in an equivalent, yet analytically more suitable form \cite{Ri1}. Namely, there holds
\be\label{wildiv}
\text{div}\big(\nabla\bH\,-\,3\,\pro(\nabla\bH)\,+\,\star\,(\nabla^\perp\bn\wedge\bH)\big)\;=\;0\:.
\ee
This reformulation in divergence form of the Willmore equation is the starting point of our analysis. In our singular situation, (\ref{wildiv}) holds only away from the origin, on $\di\setminus\{0\}$. In particular, we can define the constant $\,\bc_0\in\R^m$, called {\it residue}, by
\be\label{residudu}
\bc_0\;:=\;\int_{\partial\di}\vec{\nu}\cdot\Big(\nabla\bH\,-\,3\,\pro{(\nabla\bH)}\,+\,\star\,(\nabla^\perp\bn\wedge\bH)\Big)\:,
\ee
where $\vec{\nu}$ denotes the unit outward normal vector to $\partial\di$. We will see that the residue appears in the local asymptotic expansion of the mean curvature vector around the singularity (cf. Propositions \ref{Th3} and \ref{Th4}). \\

We next state the main result of the present paper. It concerns the regularity of the Gauss map around the point-singularity. 

\begin{Th}\label{Th1}
Let $\bp\in W^{1,2}\cap C^\infty(\di\setminus\{0\})\cap C^0(\di)$ be a conformal Willmore immersion of the punctured disk into $\R^m$, and whose Gauss map $\bn$ lies in $W^{1,2}(\di)$. Then $\nabla^2\bn\in L^{2,\infty}(\di)$, and thus in particular $\nabla\bn$ is an element of $BMO$. Furthermore, $\bn$ satisfies the pointwise estimate
\bes
|\nabla\bn(x)|\;\lesssim\;|x|^{-\epsilon}\qquad\forall\:\:\epsilon>0\:.
\ees
If the order of degeneracy of the immersion $\bp$ at the origin is at least two\footnote{Roughly speaking, $\nabla\bp(0)=\vec{0}$. The notion of ``order of degeneracy" is made precise below.}, then in fact $\nabla\bn$ belongs to $L^\infty(B_1(0))$. 
\end{Th}

\medskip

A conformal immersion of $\di\setminus\{0\}$ into $\R^m$ such that $\bp$ and its Gauss map $\bn$ both extend to maps in $W^{1,2}(\di)$ has a distinct behavior near the point-singularity located at the origin. One can show (cf. \cite{MS}, and Lemma A.5 in \cite{Ri2}) that there exists a positive natural number $\theta_0$ with
\be\label{immas}
|\bp(x)|\;\simeq\;|x|^{\theta_0}\qquad\text{and}\qquad|\nabla\bp(x)|\;\simeq\;|x|^{\theta_0-1}\qquad\text{near the origin}\:.
\ee
In addition, there holds
\bes
\la(x)\;:=\;\dfrac{1}{2}\log\Big(\frac{1}{2}\,\big|\nabla\bp(x)\big|^2\Big)\;=\;(\theta_0-1)\log|x|\,+\,u(x)\:,
\ees
where $u\in W^{2,1}(\di)$ ; and one has 
\be\label{lambdas}
\left\{\begin{array}{lcl}\nabla\la\,\in\,L^2(\di)&,&\text{when $\,\theta_0=1$}\\[1ex]
\nabla\la(x)\,\lesssim\,|x|^{-1}\,\in\,L^{2,\infty}(\di)&,&\text{when $\,\theta_0\ge2\:.$}\end{array}\right.
\ee
The integer $\theta_0$ is the density of the current $\bp_*[\di]$ at the image point $0\in\R^m$.\\

\noindent
When such a conformal immersion is Willmore on $\di\setminus\{0\}$, it is possible to refine the asymptotics (\ref{immas}). The following result describes the behavior of the immersion $\bp$ locally around the singularity at the origin. 

\begin{Prop}\label{Th2}
Let $\bp$ be as in Theorem \ref{Th1} with conformal parameter $\la$, and let $\theta_0$ be as in (\ref{immas}). 
There exists a constant vector $\bA=\bA^{1}+i\bA^{2}\in\R^2\otimes\R^m$ such that
\bes
\bA^{1}\cdot\bA^{2}\,=\,0\:\:,\qquad|\bA^{1}|\,=\,|\bA^{2}|\,=\,\theta_0^{-1}\lim_{x\rightarrow0}\,\dfrac{{e}^{\la(x)}}{|x|^{\theta_0-1}}\:\:,\qquad\pi_{\bn(0)}\bA\,=\,\vec{0}\:\:,
\ees
and
\begin{itemize}
\item[(i)] when $\theta_0=1$\:,
\be\label{rept}
\bp(x)\;=\;\Re\big(\bA\;{x}\big)\,+\,\vec{\zeta}(x)\:,%\text{O}(|x|^{\theta_0+1-\epsilon})
\ee
with $\,\vec{\zeta}\in\bigcap_{p<\infty}W^{2,p}(\di)\,$ and
\bes
\vec{\zeta}(x)\,=\,\text{O}(|x|^{2-\epsilon})\:\:,\qquad\nabla\vec{\zeta}(x)\,=\,\text{O}(|x|^{1-\epsilon})\:\:,\qquad\forall\:\:\epsilon>0\:.
\ees
\item[(ii)] when $\theta_0\ge2$\:,
\be\label{reptt}
\bp(x)\;=\;\Re\big(\bA\;{x}^{\theta_0}+\vec{B}\,x^{\theta_0+1}+\bC\,|x|^2{x}^{\,\theta_0-1}\big)\,+\,|x|^{\theta_0-1}\vec{\xi}(x)\:,
\ee
where $\bB$ and $\bC$ are constant vectors in $\C^m$. And for all $\epsilon>0$\,:
\bes
\vec{\xi}(x)\,=\,\text{O}(|x|^{3-\epsilon})\:\:,\qquad\nabla\vec{\xi}(x)\,=\,\text{O}(|x|^{2-\epsilon})\:\:,\qquad\nabla^2\vec{\xi}(x)\,=\,\text{O}(|x|^{1-\epsilon})\:.
\ees
\end{itemize}
\end{Prop}

\medskip
\noindent
The plane $span\{\bA^1,\bA^2\}$ is tangent to the surface at the origin. If $\theta_0=1$, this plane is actually $T_0\Sigma$. One can indeed show that the tangent unit vectors $\bej$ spanning $T_0\Sigma$ (defined in (\ref{lesvec})) satisfy $\bej(0)=\bA^j/|\bA^j|$. In contrast, when $\theta_0\ge2$, the tangent plane $T_0\Sigma$ does not exist in the classical sense, and the vectors $\bej(x)$ ``spin" as $x$ approaches the origin (cf. (\ref{locexphi}). More precisely, $T_0\Sigma$ is the plane $span\{\bA^1,\bA^2\}$ covered $\theta_0$ times. \\

\begin{Rm}\label{Rem1}
%From the previous proposition, one deduces
%\begin{itemize}
%\item[(i)] 
when $\theta_0=1$, the immersion $\bp$ belongs to $C^{1,\al}(\di)$ for all $\al\in[0,1)$. In general however, $\bp$ need \textbf{not} be $C^{1,1}(\di)$, as the following example shows. A conformal parametrization of the catenoid is
\bes
(r,\varphi)\;\longmapsto\;\Big(\big(r+{r}^{-1}\big)\cos(\varphi)\,,\,\big(r+{r}^{-1}\big)\sin(\varphi)\,,\,\log r\Big)\:.
\ees
Inverting the catenoid about the origin gives a Willmore surface\footnote{for it is the image of a minimal (thus Willmore) surface under a M\"obius transformation.} whose behavior near the origin consists of two identical graphs (mirror-symmetric about the $(x_1,x_2)$-plane) of order $\theta_0=1$ at the origin. One computes the inverted parametrization (for one graph) to be
\bes
\bp(r,\varphi)\;=\;\big(r\cos(\varphi)\,,\,r\sin(\varphi)\,,\,r^{2}\log r\big)\,+\,\text{O}\big(r^{3}\log^2r\big)\:.
\ees
Identifying (through $x=r{e}^{i\varphi}$) with (\ref{rept}) shows
\bes
\bA\,=\,(1\,,-\,i\,,0)\qquad\text{and}\qquad\vec{\zeta}(x)\,=\,\text{O}\big(|x|^2\log|x|\big)\:.
\ees
Thus, we cannot expect in general $\epsilon=0$ in (\ref{rept}). Moreover, $\bp\notin C^{1,1}(\di)$, and a computation reveals that
\bes
|\nabla\bn(x)|\;\simeq\;\log|x|\;\in\, BMO\setminus L^\infty(\di)\:.
\ees
%so that Theorem \ref{Th1} is sharp. 
%\item[(ii)] when $\theta_0\ge2$, the result is also sharp. 
%Consider, for example, the minimal surface conformally parametrized by its Weierstrass representation:
%\bes
%(r,\varphi)\;\longmapsto\;2\,\bigg(\log r+\dfrac{1}{2\,r^2}\cos(2\varphi)\,,\,-\,\varphi+\dfrac{1}{2\,r^2}\sin(2\varphi)\,,\,-\,\dfrac{1}{r}\cos(\varphi)\bigg)\:.
%\ees
%Inverting it about the origin gives a conformal parametrization of a Willmore surface whose immersion degenerates with order $\theta_0=2$ at the origin. Namely,
%\bes
%\bp(r,\varphi)\;=\;r^2\big(\cos(2\varphi)\,,\,\sin(2\varphi)\,,\,2\,r\cos(\varphi)\big)\,+\,\text{O}\big(1\,,1\,,r\big)\,r^4\log r\:.
%\ees
%Identifying this representation with (\ref{reptt}) gives
%\bes
%\bp(x)\;=\;\dfrac{1}{2}\,\Re\Big(\bA\,x^2+\bC\,|x|^2\,\overline{x}\Big)+\,\text{O}\big(|x|^3\log|x|\big)\:,
%\ees
%with
%\bes
%\bA=(1\,,-\,i\,,0)\quad,\quad\bB=\vec{0}\quad,\quad\bC=(0\,,0\,,2)\quad,\quad\vec{\xi}(x)=\text{O}(|x|^3\log|x|\big)\:,
%\ees
%thereby showing that we cannot in general expect $\epsilon=0$ in (\ref{reptt}).
%\end{itemize}
\end{Rm}

\medskip
It is also possible to obtain information on the local asymptotic behavior of the mean curvature vector near the origin. This is the object of the next proposition. 

%The mean curvature vector is directly related to the immersion via the identity $\,2\,\text{e}^{2\la}\bH=\Delta\bp$. As a consequence of Theorem \ref{Th2}, we thus obtain

\begin{Prop}\label{Th3}
Let $\bp$ be as in Theorem \ref{Th1}, $\la$ be its conformal parameter, and $\theta_0$ be as in (\ref{immas}). 
Locally around the singularity, the mean curvature vector satisfies
\begin{itemize}
\item[(i)] when $\theta_0=1$\:,
\bes
\bH(x)\,+\,\dfrac{\bc_0}{4\pi}\log|x|\;\in\bigcap_{p<\infty}W^{1,p}(\di)\:,
\ees
where $\bc_0$ is the residue defined in (\ref{residudu}). 
\item[(ii)] when $\theta_0\ge2$\:, 
\bes
{e}^{\la(x)}\bH(x)\;=\;f(x)\;\Re\Bigg[\bC\,\bigg(\dfrac{{x}}{|x|}\bigg)^{\!\!\theta_0-1\,}\Bigg]+\text{O}(|x|^{1-\epsilon})\qquad\forall\:\:\epsilon>0\:,
\ees
where $\bC\in\C^m$ is the same constant vector as in Proposition \ref{Th2}-(ii), and
\bes
f(x)\;:=\;2\,\theta_0\,|x|^{\theta_0-1}e^{-\la(x)}\,\in\,C^0\big(\di,(0,\infty)\big)\:.
\ees
In particular, since $\bH$ is a normal vector, we note that $\,\pi_{\bn(0)}\bC=\bC$. 
\end{itemize}
\end{Prop}

When $\theta_0\ge2$, the weighted mean curvature vector $\text{e}^{\la}\bH$ is thus bounded across the singularity (unlike in the case $\theta_0=1$, where it behaves logarithmically). But its limit may not exist: $\text{e}^{\la(x)}\bH(x)$ is a ``spinning vector" as $x$ approaches the origin\footnote{note however that the function $f(x)$ does have (positive) a limit at $x=0$, as shown in \cite{MS}.}. However, when this limit exists (and is thus necessarily zero), an interesting phenomenon occurs: both the mean curvature vector and the Gauss map undergo a ``leap of regularity". More precisely, 

\begin{Prop}\label{Th4}
Let $\bp$, $\bn$, $\bH$, and $\theta_0\ge2$ be as in Proposition \ref{Th3}. \\[.5ex]
If $\,\lim_{x\rightarrow0}{e}^{\la(x)}\bH(x)$ exists (i.e. if the vector $\bC$ from Proposition \ref{Th3}-(ii) vanishes), then there holds
\begin{itemize}
\item[(i)] $\nabla^{\theta_0+1}\bn\,\in\,L^{2,\infty}(\di)$, and hence $\,\nabla^{\theta_0}\bn\in BMO$. Furthermore,
\bes
\nabla^j\bn(x)\;=\;\text{O}\big(|x|^{\theta_0-j-\epsilon}\big)\qquad\forall\:\:\epsilon>0\:,\:\:j\in\{0,\ldots,\theta_0\}\:.
\ees
\item[(ii)] locally around the singularity\:,
\be\label{badphi}
\bp(x)\;=\;\sum_{j=0}^{\theta_0-1}\,\Re\big(\al_j\bA\;x^{\theta_0+j}\big)+\,\vec{\zeta}(x)\:,
\ee
where $\bA$ is as in Proposition \ref{Th2}, $\al_0=1$, and $\al_j\in\C^m$ are constants. The function $\vec{\zeta}$ satisfies
\bes
\nabla^j\vec{\zeta}(x)\;=\;\text{O}(|x|^{2\theta_0-j-\epsilon})\qquad\forall\:\:\epsilon>0\:\:,\:\:j\in\{0,\ldots,\theta_0\}\:;
\ees
\item[(iii)] the mean curvature vector satisfies
\bes
\bH(x)\,+\,\dfrac{\bc_0}{4\pi}\log|x|\;\in\bigcap_{p<\infty}W^{\theta_0,p}(\di)\:,
\ees
where $\bc_0$ is the residue defined in (\ref{residudu}). \\
\end{itemize}
\end{Prop}

\noindent
This apparent ``leap of regularity" is in some cases mildly surprising. It can indeed happen that the Willmore surface under consideration has been ``poorly" parametrized by $\bp$ (namely, $\bp$ parametrizes the same surface covered $\theta_0$ times), and therefore the mean curvature vector is just as regular as in the case when the point-singularity has order $\theta_0=1$. The following example sheds some light onto this phenomenon: we exhibit an ``unclever" conformal parametrization of the inverted catenoid, which degenerates at the origin with order $\theta_0\ge2$. As expected for the inverted catenoid (cf. Remark \ref{Rem1}), the mean curvature behaves logarithmically near the singularity, regardless of the order of degeneracy of the immersion.

\begin{Rm}\label{Rem2}
The result from the last proposition is sharp, as the following example shows.
We may conformally parametrize the $\theta_0$-times covered inverted catenoid by composing the parametrization of the single-covered inverted catenoid given in Remark \ref{Rem1} with $x^{\theta_0}$. Namely, 
%\bes
%(r,\varphi)\;\longmapsto\;\Big(\big(r^{\theta_0}+{r^{-\theta_0}}\big)\cos(\theta_0\,\varphi)\,,\,\big(r^{\theta_0}+{r^{-\theta_0}}\big)\sin(\theta_0\,\varphi)\,,\,\theta_0\log r\Big)\:.
%\ees
%After an inversion about the origin, it yields the same Willmore surface as in Remark \ref{Rem1}-(i), but covered $\theta_0$ times. 
the parametrization (for one graph) is
\bes
\bp(r,\varphi)\;=\;\big(r^{\theta_0}\cos(\theta_0\,\varphi)\,,\,r^{\theta_0}\sin(\theta_0\,\varphi)\,,\,r^{2\theta_0}\log r^{\theta_0}\big)\,+\,\text{O}\big(r^{3\theta_0}\log^2r\big)\:.
\ees
Identifying the latter (through $x=r{e}^{i\varphi}$) with (\ref{reptt}) shows that
\bes
\bA\,=\,(1\,,-\,i\,,0)\qquad\text{and}\qquad\vec{B}\,=\,\vec{0}\,=\,\bC\:,
\ees
so this example fits indeed within the context of Proposition \ref{Th4}. One computes explicitly the residue $\bc_0$ in this case, namely
\bes
\bc_0\;=\;-\,16\,\pi\,\theta_0\,(0\,,0\,,1)\:.
\ees
Moreover, there holds
\bes
|\nabla\bn(x)|\;\simeq\;|x|^{\theta_0-1}\log|x|\:,
\ees
thereby confirming that
\bes
|\nabla\bn(x)|\;\lesssim\;|x|^{\theta_0-1-\epsilon}\qquad\text{for all $\epsilon>0$, but not for $\epsilon=0$\:.}
\ees
\end{Rm}

\medskip
\noindent
It is currently unknown to the authors whether there exist Willmore immersions which degenerate at the origin with the order $\theta_0$, for which $\text{e}^{\la}\bH$ has a limit at the origin, and which do {\it not} parametrize a $\theta_0$-times-covered Willmore surface whose immersion degenerates at the origin with order 1. It seems never to be the case for branched inverted minimal surface (in $\R^3$ at least). Admittedly however, inverted minimal surfaces are a very special kind of Willmore surfaces.

\bigskip
Finally, when the residue $\bc_0=\vec{0}$ (\textbf{and} in the case $\theta_0\ge2$ the constant vector $\bC$ from Proposition \ref{Th3}-(ii) vanishes: $\bC=\vec{0}$), the singularity at the origin is removable. Namely,

\begin{Th}\label{Th5}
Under the hypotheses of Proposition \ref{Th3}, if $\theta_0=1$ and $\bc_0=\vec{0}$, or if $\theta_0\ge2$ and $\bc_0=\vec{0}=\bC$, then the immersion $\bp$ is smooth across the unit-disk. 
\end{Th}

\noindent
This is in particular the case for branched minimal immersions. \\[3ex]

\noindent
\textbf{Acknowledgments:\:\:} The first author is grateful to the DFG Collaborative Research Center SFB/Transregio 71 (Project B3) for fully supporting his research. Parts of this work were completed during the first author's visits to the Forschungsinstitut f\"ur Mathematik at the ETH. Welcoming facilities both at the Albert-Ludwigs-Universit\"at in Freiburg and at the ETH in Z\"urich have significantly and positively impacted the development of this work. The first author is also grateful to Ernst Kuwert for suggesting and discussing this problem.

\medskip

\reset

\section{Proof of Theorems}

\subsection{Fundamental Results and Reformulation}

We place ourselves in the situation described in the introduction. Namely, we have a Willmore immersion $\bp$ on the punctured disk which degenerates at the origin in such a way that
\bes
|\bp(x)|\;\simeq\;|x|^{\theta_0}\qquad\text{and}\qquad|\nabla\bp(x)|\;=\;\sqrt{2}\,\text{e}^{\la(x)}\;\simeq\;|x|^{\theta_0-1}\:,
\ees
for some $\theta_0\in\N\setminus\{0\}$.

\medskip

Amongst the analytical tools available to the study of weak Willmore immersions with square-integrable second fundamental form, the most important one is certainly the ``$\eps$-regularity". The version appearing in Theorem 2.10 and Remark 2.11 from \cite{KS3} (see also Theorem I.5 in \cite{Ri1}) states that there exists $\eps_0>0$ such that, if
\be\label{epsregcond}
\int_{B_1(0)}|\nabla\bn|^2\,dx\;<\;\eps_0\:,
\ee
then there holds
\be\label{epsreg}
\Vert\nabla\bn\Vert_{L^\infty(B^g_\sigma)}\;\le\;\dfrac{C}{\sigma}\,\Vert\nabla\bn\Vert_{L^2(B^g_{2\sigma})}\qquad\:\:\forall\:\:\;B^g_{2\sigma}\subseteq\Om:=\di\setminus\{0\}\:,
\ee
where $B^g_\sigma$ is a geodesic disk of radius $\sigma$ for the induced metric $g=\bp^*g_{\R^m}$, and $C$ is a universal constant. \\[1.5ex]
The $\eps$-regularity enables us to obtain the following result (already observed in \cite{KS2}), decisive to the remainder of the argument.

\begin{Lma}\label{delta}
The function $\,\delta(r):=r\sup_{|x|=r}|\nabla\bn(x)|\,$ satisfies
\bes
\lim_{r\searrow0}\:\delta(r)\,=\,0\hspace{1cm}\text{and}\hspace{1cm}\int_{0}^{1}\delta^2(r)\,\dfrac{dr}{r}\,<\,\infty\:.
\ees
\end{Lma}
$\textbf{Proof.}$ 
From (\ref{eraser}) and (\ref{immas}), the metric $g$ satisfies 
\bes
g_{ij}(x)\;\simeq\;|x|^{2(\theta_0-1)}\delta_{ij}\qquad\text{on}\:\:\:\:B_{2r}(0)\setminus B_{r/2}(0)\qquad\forall\:\:r\in(0\,,1/2)\:.
\ees
A simple computation then shows that
\be\label{barbitruc}
B^g_{2cr^{\theta_0}}(x)\;\subset\;B_{2r}(0)\setminus B_{r/2}(0)\qquad\forall\:\:\:x\in \partial B_{r}(0)\:,
\ee
where $\,0<2\theta_0\,c<1-2^{-\theta_0}$.\\
Since the metric $g$ does not degenerate away from the origin, given $r<1/2$, we can always cover the flat circle $\partial B_{r}(0)$ with finitely many metric disks:
\bes
\partial B_r(0)\;\subset\;\bigcup_{j=1}^{N} B^g_{cr^{\theta_0}}(x_j)\qquad\text{with}\quad x_j\in\partial B_r(0)\:.
\ees
Hence, per the latter, (\ref{epsreg}), and (\ref{barbitruc}), we obtain
\begin{eqnarray}\label{barbibulle}
r\sup_{|x|=r}|\nabla\bn(x)|&\leq&r\sup_{|x|=r}\Vert\nabla\bn\Vert_{L^\infty(B^g_{cr^{\theta_0}}(x))}\;\;\lesssim\;\;\sup_{|x|=r}\Vert\nabla\bn\Vert_{L^2(B^g_{2cr^{\theta_0}}(x))}\nonumber\\[1ex]
&\leq&\Vert\nabla\bn\Vert_{L^2(B_{2r}(0)\setminus B_{r/2}(0))}\:.
\end{eqnarray}
As $\nabla\bn$ is square-integrable by hypothesis, letting $r$ tend to zero in the latter yields the first assertion. \\
The second assertion follows from (\ref{barbibulle}), namely, 
\bes
\int_{0}^{1/2}\delta^2(r)\,\dfrac{dr}{r}\;\lesssim\;\int_{0}^{1/2}\Vert\nabla\bn\Vert^2_{L^2(B_{2r}(0)\setminus B_{r/2}(0))}\,\dfrac{dr}{r}\;=\;\log(4)\;\Vert\nabla\bn\Vert^2_{L^2(B_1(0))}\:,
\ees
%\begin{eqnarray*}
%\int_{0}^{2/3}\delta^2(r)\,\dfrac{dr}{r}&=&\int_{0}^{2/3}\Vert\nabla\bn\Vert^2_{L^\infty(\partial B_r(0))}\,r\,dr\\[1ex]
%&\lesssim&\int_{0}^{2/3}\Vert\nabla\bn\Vert^2_{L^2(B_{3r/2}(0)\setminus B_{r/2}(0))}\,\dfrac{dr}{r}\:\:=\:\:\dfrac{4}{3}\;\Vert\nabla\bn\Vert^2_{L^2(B_1(0))}\:\,,%\\[1ex]
%&\simeq&\int_{0}^{2\pi}\int_{0}^{2/3}\int_{r/2}^{3r/2}|\nabla\bn(s,\varphi)|^2\,ds\,dr\,d\varphi\\[1ex]
%&=&\dfrac{4}{3}\int_{0}^{2\pi}\int_{0}^{1}|\nabla\bn(s,\varphi)|^2\,s\,ds\,d\varphi\;\;\lesssim\;\;\Vert\nabla\bn\Vert^2_{L^2(B_1(0))}\:\,,
%\end{eqnarray*}
which is by hypothesis finite.\\[-3ex]

$\hfill\blacksquare$ \\

Recalling (\ref{eraser}) linking the Gauss map to the mean curvature vector and the fact that $\text{e}^{\la(x)}\simeq|x|^{\theta_0-1}$, we obtain from Lemma \ref{delta} that
\be\label{es2}
r^{\theta_0}\!\sup_{|x|=r}|\bH(x)|\:\leq\:r^{\theta_0}\!\sup_{|x|=r}\text{e}^{-\la(x)}|\nabla\bn(x)|\:\lesssim\:\delta(r)\:.
\ee
The Willmore equation (\ref{wildiv}) may be alternatively written
\bes
\text{div}\Big(\nabla\bH\,-\,3\,\pro{(\nabla\bH)}\,-\,\star\,(\bn\wedge\nabla^\perp\bH)\Big)\:=\:0\:\qquad\text{on}\quad\Om:=B_1(0)\setminus\{0\}\:.
\ees
It is elliptic \cite{Ri1}. Using the information on the gradient of $\bn$ given by (\ref{epsreg}), and some standard analytical techniques for second-order elliptic equations in divergence form (cf. [GW]), one deduces from (\ref{es2}) that
\be\label{es3}
r^{\theta_0+1}\!\sup_{|x|=r}|\nabla\bH(x)|\:\lesssim\:\delta(r)\:.
\ee
These observations shall be helpful in the sequel. \\

%As previously noted, there holds on $\Omega:=B_1(0) -\{0\}$ the identity
%\be\label{willy0}
%\text{div}\Big(\nabla\bH\,-\,3\,\pro{(\nabla\bH)}\,+\,\star\,(\nabla^\perp\bn\wedge\bH)\Big)\:=\:0\:.
%\ee
The equation (\ref{wildiv}) implies that for any ball $B_\rho(0)$ of radius $\rho$ centered on the origin and contained in $\Om$, there holds
\be\label{residue}
\int_{\partial B_\rho(0)} \vec{\nu}\cdot\Big(\nabla\bH\,-\,3\,\pro{(\nabla\bH)}\,+\,\star\,(\nabla^\perp\bn\wedge\bH)\Big)\:=\:\bc_0\:,\qquad\forall\:\:\rho\in(0,1)\:,
\ee
where $\bc_0$ is the residue defined in (\ref{residudu}). 
Here $\vec{\nu}$ denotes the unit outward normal vector to $\partial B_\rho(0)$.
An elementary computation shows that
\bes
\int_{\partial B_{\rho}(0)} \vec{\nu}\cdot\nabla\log|x|\:=\:2\pi\:,\qquad\forall\:\:\rho>0\:.
\ees
Thus, upon setting
\be\label{a1}
\bX\;:=\;\nabla\bH\,-\,3\,\pro{(\nabla\bH)}\,+\,\star\,(\nabla^\perp\bn\wedge\bH)\,-\,\dfrac{\bc_0}{2\pi}\,\nabla\log|x|\:,
\ee
we find
\bes
\text{div}\,\bX=0\quad\:\:\text{on}\:\:\Om\:,\hspace{.7cm}\text{and}\hspace{.5cm}\int_{\partial B_\rho(0)}\vec{\nu}\cdot\bX\,=\,{0}\qquad\forall\:\:\rho\in(0,1)\:.
\ees
As $\bX$ is smooth away from the origin, the Poincar\'e Lemma implies now the existence of an element $\bL\in C^{\infty}(\Om)$ such that
\be\label{a2}
\bX\;=\;\nabla^\perp\bL\qquad\text{on}\:\:\:\Om\:.
\ee
We deduce from Lemma \ref{delta} and (\ref{es2})-(\ref{a2}) that
\be\label{xx1}
\int_{B_1(0)}|x|^{2\theta_0}|\nabla\bL|^2\,dx\:\,\lesssim\:\,\int_{0}^{1}\delta^2(s)\,\dfrac{ds}{s}\:<\:\infty\:.
\ee
A classical Hardy-Sobolev inequality gives the estimate
\be\label{estimL}
\theta_0^2\int_{B_1(0)}|x|^{2(\theta_0-1)}|\bL|^2\,dx\:\le\:\int_{B_1(0)}|x|^{2\theta_0}|\nabla\bL|^2\,dx\;+\;\theta_0\!\int_{\partial B_1(0)}|\bL|^2\:, %r^{2\theta_0-1}\,
\ee
which is a finite quantity, owing to (\ref{xx1}) and to the smoothness of $\bL$ away from the origin. The immersion $\bP$ has near the origin the asymptotic behavior $\,|\nabla\bp(x)|\simeq|x|^{\theta_0-1}$. Hence (\ref{estimL}) yields that 
\be\label{xx2}
\bL\cdot\nabla\bP\:,\;\bL\wedge\nabla\bP\;\in\;L^2(B_1(0))\:.
\ee

We next set $\,\bF(x):=\dfrac{\bc_0}{2\pi}\log|x|\,$, and define the functions $g$ and $\bG$ via
\be\label{sys1}
\left\{\begin{array}{rclcrclcl}
\Delta g&=&\nabla\bF\cdot\nabla\bp&,\quad&\Delta\bG&=&\nabla\bF\wedge\nabla\bp&\quad&\text{in}\:\:\: B_1(0)\\[1ex]
g&=&0&,\quad&\bG&=&\vec{0}&\quad&\text{on}\:\:\:\partial B_1(0)\:.
\end{array}\right.
\ee
Since $\,|\nabla\bP(x)|\simeq|x|^{\theta_0-1}\,$ near the origin and $\bF$ is the fundamental solution of the Laplacian, by applying Calderon-Zygmund estimates to (\ref{sys1}), we find\footnote{The weak-$L^2$ Marcinkiewicz space $L^{2,\infty}(B_1(0))$ is defined as those functions $f$ which satisfy $\:\sup_{\al>0}\al^2\Big|\big\{x\in B_1(0)\,;\,|f(x)|\ge\al\big\}\Big|<\infty$. In dimension two, the prototype element of $L^{2,\infty}$ is $|x|^{-1}\,$. The space $L^{2,\infty}$ is also a Lorentz space, and in particular is a space of interpolation between Lebesgue spaces, which justifies the first inclusion in (\ref{regg}). See \cite{He} or \cite{Al} for details.} 
\be\label{regg}
\nabla^2g\;,\;\nabla^2\bG\:\;\in\:\left\{\begin{array}{lcl}L^{2,\infty}(B_1(0))&,&\theta_0=1\\[1.5ex]
BMO(B_1(0))&,&\theta_0\ge2\:.
\end{array}\right.
\ee

In the paper \cite{BR1} (cf. Lemma A.2), the authors derive the identities\footnote{Observe that $\,\nabla^\perp\bL+\nabla\bF\,$ is exactly the divergence-free quantity appearing in (\ref{wildiv}).}:
\be\label{constraints}
\left\{\begin{array}{lcl}
\nabla\bP\cdot(\nabla^\perp\bL+\nabla\bF)&=&0\\[1ex]
\nabla\bP\wedge(\nabla^\perp\bL+\nabla\bF)&=&-\,2\,\nabla\bP\wedge\nabla\bH\:.\end{array}\right.
\ee
Accounted into (\ref{sys1}), the latter yield that there holds in $\Om$\,:
\be\label{div1}
\left\{\begin{array}{rcl}
\text{div}\big(\bL\cdot\nabla^\perp\bP\,-\,\nabla g\big)&=&0\\[1.5ex]
\text{div}\big(\bL\wedge\nabla^\perp\bP\,-\,2\,\bH\wedge\nabla\bP\,-\,\nabla\bG\big)&=&\vec{0}\:,\end{array}\right.
\ee
where we have used the fact that
\bes
\Delta\bP\wedge\bH\;=\;2\,\text{e}^{2\la}\bH\wedge\bH\;=\;\vec{0}\:.
\ees
Note that the terms under the divergence symbols in (\ref{div1}) both belong to $L^2(B_1(0))$, owing to (\ref{xx2}) and (\ref{regg}). The distributional equations (\ref{div1}), which are {\it a priori} to be understood on $\Om$, may thus be extended to all of $B_1(0)$. Indeed, a classical result of Laurent Schwartz states that the only distributions supported on $\{0\}$ are linear combinations of derivatives of the Dirac delta mass. Yet, none of these (including delta itself) belongs to $W^{-1,2}$. We shall thus understand (\ref{div1}) on $B_1(0)$. 
It is not difficult to verify (cf. Corollary IX.5 in [DL]) that a divergence-free vector field in $L^2(B_1(0))$ is the curl of an element in $W^{1,2}(B_1(0))$. 
%Owing to (\ref{xx2}) and (\ref{regg}), 
We apply this observation to (\ref{div1}) so as to infer the existence of two functions\footnote{$S$ is a scalar while $\bR$ is $\bigwedge^2(\R^m)$-valued.} $S$ and of $\bR$ in the space $W^{1,2}(B_1(0))\cap C^\infty(\Om)$, with
\be\label{truxy}
\left\{\begin{array}{rclll}
\nabla^\perp S&=&\bL\cdot\nabla^\perp\bP\,-\,\nabla g&&\\[1.5ex]
\nabla^\perp\bR&=&\bL\wedge\nabla^\perp\bP\,-\,2\,\bH\wedge\nabla\bP\,-\,\nabla\bG\:.
\end{array}\right.
\ee
Moreover, $S$ and $\bR$ may be chosen to be constant on the boundary of the unit disk ; without loss of generality, we shall assume that $\,S\big|_{\partial B_1(0)}=0\,$ and $\,\bR\big|_{\partial B_1(0)}=\vec{0}.$ \\
According to the identities (\ref{sysSR0}) from the Appendix, the functions $S$ and $\bR$ satisfy on $B_1(0)$ the following system of equations, called {\it conservative conformal Willmore system}\footnote{refer to the Appendix for the notation and the operators used.}:
\be\label{sysSR}
\left\{\begin{array}{rclll}
\Delta S&=&-\,\nabla(\star\,\bn)\cdot\nabla^\perp\bR\,-\,\text{div}\big((\star\,\bn)\cdot\nabla\bG\big)&&\\[1.5ex]
\Delta\bR&=&-\,\nabla(\star\,\bn)\bul\nabla^\perp\bR\;+\;\nabla(\star\,\bn)\cdot\nabla^\perp S\\[.75ex]
&&\hspace{2.45cm}\,-\:\,\text{div}\big((\star\,\bn)\bul\nabla\bG\,-\,\star\,\bn\,\nabla g\big)\:.
\end{array}\right.
\ee

Not only is this system independent of the codimension, but it further displays two fundamental advantages. Analytically, (\ref{sysSR}) is uniformly elliptic. This is in sharp contrast with the Willmore equation (\ref{wil1}) whose leading order operator $\Delta_{\perp}$ degenerates at the origin, owing to the presence of the conformal factor $\text{e}^{\la(x)}\simeq|x|^{\theta_0-1}$. 
Structurally, the system (\ref{sysSR}) is in divergence form. We shall in the sequel capitalize on this remarkable feature to develop arguments of ``integration by compensation". {\it A priori}\, however, since $\bn$, $S$, and $\bR$ are elements of $W^{1,2}$, the leading terms on the right-hand side of the conservative conformal Willmore system (\ref{sysSR}) are critical. This difficulty is nevertheless bypassed using the fact that the $W^{1,2}(B_1(0))$-norm of the Gauss map $\bn$ is chosen small enough (cf. (\ref{acheumeuneu})). 
%This system is {\it a priori} only valid on the punctured disk $\Om=B_1(0) -\{0\}$. However, from the distributional point of view, it actually holds on the whole unit-disk $B_1(0)$. This is because the functions $g$ and $\bG$ belong to $W^{2,(2,\infty)}(B_1(0))$, the Gau\ss\, map $\bn$ is bounded and has its gradient in $L^{2}(B_1(0))$, and because $S$ and $\bR$ lie in $W^{1,2}(\Om)$. Hence all quantities involved in (\ref{sysSR}) belong to either $L^1(\Om)$ or to $W^{-1,2}(\Om)$. A classical result of Laurent Schwartz guarantees that there lies in these spaces no non-trivial distribution supported on $\{0\}$. The system (\ref{sysSR}) holds thus on all of $B_1(0)$.\\

\subsection{The general case when $\theta_0\ge1$}

We have gathered enough information about the functions involved to apply to the system (\ref{sysSR}) (a slightly extended version of) Proposition \ref{morreydecay} and thereby obtain that
\be\label{gradSinLp}
\nabla S\:,\:\nabla\bR\:\in\:L^{p}(B_{1}(0))\qquad\text{for some}\:\:p>2\:.
\ee
It is shown in the Appendix (cf. (\ref{delphi0})) that
\be\label{delphi}
-\,2\,\Delta\bP\;=\;\big(\nabla S-\nabla^\perp g\big)\cdot\nabla^\perp\bP\,-\,\big(\nabla\bR-\nabla^\perp\bG\big)\bul\nabla^\perp\bP\:.
\ee
Hence, as $\,|\nabla\bP(x)|\simeq\text{e}^{\la(x)}\simeq|x|^{\theta_0-1}$ around the origin, using (\ref{regg}) and (\ref{gradSinLp}), we may call upon Proposition \ref{CZpondere} with the weight $|\mu|=\text{e}^{\la}$ and $a=\theta_0-1$ to conclude that
\bes
\nabla\bp(x)\;=\;\bPe(\overline{x})+\,\text{e}^{\la(x)}\bT(x)\:,
\ees
where $\bPe$ is a $\C^m$-valued polynomial of degree at most $(\theta_0-1)$, and $\vec{T}(x)=\text{O}\big(|x|^{1-\frac{2}{p}-\epsilon}\big)$ for every $\epsilon>0$. Because $\,\text{e}^{-\la}\nabla\bp\,$ is a bounded function, we deduce more precisely that $\bPe(\overline{x})=\theta_0\bA^*\,\overline{x}^{\,\theta_0-1}$, for some constant vector $\bA\in\C^m$\, (we denote its complex conjugate by $\bA^*$), so that
\be\label{locexphi}
\nabla\bp(x)\;=\;\left(\begin{array}{c}\Re\\[.5ex]-\,\Im\end{array}\right)\big(\theta_0\,\vec{A}\;{x}^{\theta_0-1}\big)\,+\,\text{e}^{\la(x)}\bT(x)\:.
\ee
Equivalently, upon writing $\bA=\bA^{1}+i\bA^{2}$, where $\bA^{1}$ and $\bA^{2}$ are two vectors in $\R^m$, the latter may be recast as\footnote{$\varphi$ denotes the argument of $x$.}
\bes
\left\{\begin{array}{rcl}\px\bP(x)&\!\!=\!\!&\theta_0\,|x|^{\theta_0-1}\Big[\bA^{1}\cos\big((\theta_0-1)\varphi\big)-\bA^{2}\sin\big((\theta_0-1)\varphi\big)\Big]+\,\text{e}^{\la}\Re(\bT(x))\\[1.5ex]
-\py\bP(x)&\!\!=\!\!&\theta_0\,|x|^{\theta_0-1}\Big[\bA^{2}\cos\big((\theta_0-1)\varphi\big)+\bA^{1}\sin\big((\theta_0-1)\varphi\big)\Big] -\,\text{e}^{\la}\Im(\bT(x))     \:.
\end{array}\right.
\ees
The conformality condition of $\bp$ shows easily that
\be\label{condA1}
|\bA^{1}|\;=\;|\bA^{2}|\qquad\text{and}\qquad\bA^{1}\cdot\bA^{2}\;=\;0\:.
\ee
Yet more precisely, as $|\nabla\bp|^2=2\,\text{e}^{2\la}$, we see that
\be\label{condA2}
|\bA^{1}|\;=\;|\bA^{2}|\;=\;\dfrac{1}{\theta_0}\,\lim_{x\rightarrow0}\,\dfrac{\text{e}^{\la(x)}}{|x|^{\theta_0-1}}\,\in\:\:]0\,,\infty[\:.
\ee
Because $\bp(0)=\vec{0}$, we obtain from (\ref{locexphi}) the local expansion
\bes
\bp(x)\;=\;\Re\big(\bA\,x^{\theta_0}\big)\,+\,\text{O}\big(|x|^{\theta_0-\frac{2}{p}-\epsilon}\big)\:.
\ees
%Yet more precisely, it follows from (\ref{locexphi}) that
%\be\label{vecA}
%\dfrac{\bA_j}{|\bA_j|}\;=\;\bej(0)\:.
%\ee 
Since $\pro\nabla\bp\equiv\vec{0}$, we deduce from (\ref{locexphi}) that
\be\label{expA}
\pro\bA\;=\;-\,\theta_0^{-1}\,{x}^{1-\theta_0}\,\text{e}^{\la}\,\pro\bT^*(x)\;=\;\text{O}\big(|x|^{1-\frac{2}{p}-\epsilon}\big)\qquad\:\:\:\forall\:\:\epsilon>0\:.
\ee
Let now $\delta:=1-\frac{2}{p}\in(0,1)$, and let $0<\eta<p$ be arbitrary. We choose some $\epsilon$ satisfying
\bes
0\;<\;\epsilon\;<\;\frac{2\,\eta}{p(p-\eta)}\;\equiv\;\delta-1+\frac{2}{p-\eta}\:.
\ees
We have observed that $\,\pro\bA=\text{O}(|x|^{\delta-\epsilon})$, hence $\,
\pro\bA=\text{o}\big(|x|^{1-\frac{2}{p-\eta}}\big)$\,,
and in particular, we find
\be\label{pina}
\dfrac{1}{|x|}\,\pi_{\bn(x)}\bA\;\in\;L^{p-\eta}(B_1(0))\:\:\:\:\qquad\forall\:\:\eta>0\:.
\ee
This fact shall come helpful in the sequel.

\medskip
When $\theta_0=1$, one directly deduces from the standard Calderon-Zygmund theorem applied to (\ref{delphi}) that $\nabla^2\bp\in L^p$. In that case, $\text{e}^{\la}$ is bounded from above and below, and thus the identity (cf. (\ref{timide}) in the Appendix)
\be\label{modulegradn}
\big|\nabla\bn\big|\;=\;\text{e}^{-\la}\big|\pro\nabla^2\bp\big|
\ee
yields that $\nabla\bn\in L^p$. When now $\theta_0\ge2$, we must proceed slightly differently to obtain analogous results. From (\ref{lambdas}), we know that $|x|\nabla\la(x)$ is bounded across the unit-disk. We may thus apply Proposition \ref{CZpondere}-(ii) to (\ref{delphi}) with the weight $|\mu|=\text{e}^{\la}$ and $a=\theta_0$. The required hypothesis (\ref{hypw2}) is fullfilled, and we so obtain
%A simple computation (cf. (\ref{atchoum}) in the Appendix) shows that
%\bes
%\text{e}^{-\la}\,\pi_{T}\nabla^2\bP\;=\;\left(\begin{array}{cc}\px\la&\py\la\\[1ex]\py\la&-\,\px\la\end{array}\right)\vec{e}_1\;+\;\left(\begin{array}{cc}-\,\py\la&\px\la\\[1ex]\px\la&\py\la\end{array}\right)\vec{e}_2\:.
%\ees
%Since $\,|\nabla\la(x)|\simeq|x|^{-1}$ near the origin (for $\theta_0\ge2$), we see that
%\be\label{saumon}
%|x|\,\text{e}^{-\la}\,\pi_{T}\nabla^2\bP\;\in\;L^\infty(B_1(0))\:.
%\ee
%Identity (\ref{timide}) from the Appendix implies
%\be\label{modulegradn}
%\big|\nabla\bn\big|\;=\;\text{e}^{-\la}\big|\pro\nabla^2\bp\big|\:.
%\ee
%Hence, by hypothesis, 
%\be\label{ohla}
%\text{e}^{-\la}\,\pro\nabla^2\bp\;\in\;L^2(B_1(0))\:.
%\ee
%Moreover, using the $\eps$-regularity (\ref{epsreg}), it ensues that
%$\,|x|\text{e}^{-\la}\pro\nabla^2\bp$\, is a bounded function on $B_1(0)$. Coupled to (\ref{saumon}), this yields that
%\bes\label{obese}
%|x|\,\text{e}^{-\la}\nabla^2\bp\;\in\;L^\infty(B_1(0))\:,
%\ees
%which is the hypothesis (\ref{hypw2}) required to apply Proposition \ref{CZpondere}-(ii) to (\ref{delphi}). Thus,
\be\label{locexdelphi}
\nabla^2\bp(x)\;=\;\theta_0\,(1-\theta_0)\left(\begin{array}{cc}-\,\Re&\Im\\[.5ex]\Im&\Re\end{array}\right)\big(\vec{A}\;{x}^{\theta_0-2}\big)\,+\,\text{e}^{\la(x)}\vec{Q}(x)\:,
\ee
where $\vec{A}$ is as in (\ref{locexphi}), while $\vec{Q}$ belongs to $\mathbb{R}^4\otimes L^{p-\epsilon}(B_1(0))$ for every $\epsilon>0$. The exponent $p>2$ is the same as in (\ref{gradSinLp}).\\
Since $\,\text{e}^{\la(x)}\simeq|x|^{\theta_0-1}$, we obtain from (\ref{locexdelphi}) that
\bes
\text{e}^{-\la}\big|\pro\nabla^2\bP\big|\;\lesssim\;|x|^{-1}|\pro\vec{A}|\,+\,|\pro\,\vec{Q}|\:.
\ees
According to (\ref{pina}), the first summand on the right-hand side of the latter belongs to $L^{p-\eta}\,$ for all $\eta>0$. Moreover, we have seen that $\pro\bQ$ lies in $L^{p-\epsilon}$ for all $\epsilon>0$. Whence, it follows that $\,\text{e}^{-\la}\pro\nabla^2\bP\,$ is itself an element of $L^{p-\epsilon}$ for all $\epsilon>0$. Brought into (\ref{modulegradn}), this information implies that
%From this and (\ref{ohla}), we thus see that
%\bes
%\pro\big(\vec{A}\big)\,\text{e}^{-\la}\nabla\,\overline{x}^{\,\theta_0-1}\,=\;\text{e}^{-\la}\pro\nabla^2\bP\,+\,\pro\,\vec{Q}\;\in\;L^2(B_1(0))\:.
%\ees
%But the left-hand side of the latter behaves near the origin like $|x|^{-1}$ which is not square-integrable. Hence $\pro\big(\vec{A}\big)=\vec{0}$\,, thereby yielding from (\ref{locexdelphi}) that
%\be\label{zoumette}
%\pro\nabla^2\bp\;=\;-\,\text{e}^{\la}\pro\,\vec{Q}\:.
%\ee
%As $\,\vec{Q}\in L^{p-\epsilon}$\,, it follows from the latter and (\ref{modulegradn}) that
\be\label{regn}
\nabla\bn\;\in\;L^{p-\epsilon}(B_{1}(0))\:,\qquad\quad\forall\:\:\epsilon>0\:.
\ee
In light of this new fact, we may now return to (\ref{sysSR}). In particular, recalling (\ref{regg}), we find 
\bes
\Delta S\;\equiv\;-\,\nabla(\star\,\bn)\cdot\big(\nabla^\perp\bR\,+\nabla\bG\big)\,-\,(\star\,\bn)\cdot\Delta\bG\:\:\in\:\:L^{q}(B_{1}(0))\:,
\ees
with
\bes
\dfrac{1}{q}\:=\:\dfrac{1}{p}\,+\,\dfrac{1}{p-\epsilon}\:.
\ees
We attract the reader's attention on an important phenomenon occurring when $\theta_0=1$. In this case, if the aforementioned value of $q$ exceeds 2 (i.e. if $p>4$), then $\Delta S\notin L^q$, but rather only $\Delta S\in L^{2,\infty}$. This integrability ``barrier" stems from that of $\Delta\bG$, as given in (\ref{regg}). The same considerations apply naturally with $\bR$ and $g$ in place of $S$ and $\bG$, respectively. \\[1ex]
Our findings so far may be summarized as follows:
\be\label{regsr}
\nabla S\;\,,\,\nabla\bR\:\,\in\:\,\left\{\begin{array}{lcl}W^{1,(2,\infty)}&,&\text{if}\:\:\:\theta_0=1\:\:\:\text{and}\:\:\:p>4\\[1ex]
W^{1,q}&,&\text{otherwise}.%\\[1ex]
%W^{1,q}&,&\text{if}\:\:\:\theta_0\ge 2\:.
\end{array}\right.
\ee
With the help of the Sobolev embedding theorem\footnote{we also use a result of Luc Tartar \cite{Ta} stating that $\,W^{1,(2,\infty)}\subset{BMO}$.}, we infer
\be\label{improvesr}
\nabla S\;\,,\,\nabla\bR\:\,\in\:\,\left\{\begin{array}{lcl}BMO&,&\text{if}\:\:\:\theta_0=1\:\:\:\text{and}\:\:\:p>4\\[1ex]
%L^s&,&\text{if}\:\:\:\theta_0=1\:\:\:\text{and}\:\:\:p\le4\\[1ex]
L^\infty&,&\text{if}\:\:\:\theta_0\ge2\:\:\:\text{and}\:\:\:p>4\\[1ex]
L^s&,&\text{if}\:\:\:\theta_0\ge1\:\:\:\text{and}\:\:\:p\le4\:,
\end{array}\right.
\ee
with
\bes
\dfrac{1}{s}\:=\:\dfrac{1}{q}\,-\,\dfrac{1}{2}\:\,=\:\,\dfrac{1}{p}\,+\,\dfrac{1}{p-\epsilon}\,-\,\dfrac{1}{2}\:\,<\,\:\dfrac{1}{p}\:.
%\dfrac{(2-p)(p-\epsilon)+2p}{2p(p-\epsilon)}\:<\:\dfrac{1}{p}\:.
\ees
Comparing (\ref{improvesr}) to (\ref{gradSinLp}), we see that the integrability has been improved. The process may thus be repeated until reaching that
\bes
\nabla S\:\,,\,\nabla\bR\:\,\in\:L^{b}(B_{1}(0))\qquad\:\:\forall\:\:\:b<\infty
\ees
holds in all configurations. With the help of this newly found fact, we reapply Proposition \ref{CZpondere} so as to improve (\ref{regsr}) and (\ref{regn}) to
\be\label{srinsobolev}
\nabla S\;\,,\,\nabla\bR\:\,\in\:\,\left\{\begin{array}{lcl}W^{1,(2,\infty)}(B_{1}(0))&,&\text{if}\:\:\:\theta_0=1\\[1ex]
W^{1,b}(B_{1}(0))&,&\text{if}\:\:\:\theta_0\ge2\:,\qquad\forall\:\:b<\infty
\end{array}\right.
\ee
and
\be\label{regn2}
\nabla\bn\,\in\,L^{b}(B_{1}(0))\,\quad\quad\forall\:\:b<\infty\:.
\ee
The $\eps$-regularity in the form (\ref{barbibulle}) then yields pointwise estimates for the Gauss map. Namely, in a neighborhood of the origin, 
\bes
|\nabla\bn(x)|\;\lesssim\;|x|^{-\epsilon}\qquad\forall\:\:\epsilon>0\:.
\ees

\subsection{The case $\theta_0=1$}\label{tet1}

We shall now investigate further the case $\theta_0=1$, when $|\nabla\bp|\simeq\text{e}^{\la}$ is bounded from both above and below around the origin. Setting
\be\label{f1f20}
\bF_1\;:=\;\nabla^\perp\bR\,+\nabla\bG\qquad\text{and}\qquad F_2\;:=\;\nabla^\perp S\,+\nabla g
\ee
in (\ref{delphi}) gives
\be\label{delphi20}
2\,\Delta\bP\;=\;\,F_2\cdot\nabla\bP\,-\,\bF_1\bul\nabla\bP\:.
\ee
According to (\ref{regg}) and (\ref{regsr}), the right-hand side of the latter has bounded mean oscillations. Hence $\nabla^2\bp\in\bigcap_{p<\infty}L^p$. Using the fact that $2\,\text{e}^{2\la}\bH=\Delta\bp$, we differentiate (\ref{delphi20}) to obtain
\bes
4\,\nabla\big(\text{e}^{2\la}\bH\big)\;=\;\nabla F_2\cdot\nabla\bP\,-\,\nabla\bF_1\bul\nabla\bP\,+\,F_2\cdot\nabla^2\bP\,-\,\bF_1\bul\nabla^2\bP\;\in\;L^{2,\infty}\:.
\ees
This shows that $\bH\in BMO$. Moreover, since $\nabla\la\in L^{2}$, it follows that $\nabla\bH\in L^{2,\infty}\subset\bigcap_{1\le p<2}L^p$. We shall now obtain an asymptotic expansion for $\bH(x)$ near the origin. To achieve this, we use a ``generic" procedure, which will be called upon again in section \ref{lim0}.

\begin{Prop}\label{logex}
Let the immersion $\bp$ satisfy an expansion of the type (\ref{locexphi}), for all $p<\infty$. Suppose that $\,\bn\in\bigcap_{p<\infty} W^{1,p}(B_1(0))\,$ and $\,\bH\in\bigcap_{p<2}W^{1,p}(B_1(0))$. Then locally around the origin,
\bes
\bH(x)+\,\frac{\bc_0}{4\pi}\,\log|x|\;\in\bigcap_{p<\infty}W^{1,p}(B_1(0))\:,
\ees
where $\bc_0$ is the residue defined in (\ref{residue}). 
\end{Prop}
$\textbf{Proof.\:}$ In order to demonstrate this result, one must return to the formalism developed in \cite{Ri1}, where it is shown that
\bes
\mathcal{L}(\bH)\;:=\;\text{div}\Big(\nabla\bH\,-\,3\,\pro\nabla\bH\,+\,\star\,(\nabla^\perp\bn\wedge\bH)\Big)\:=\:0\qquad\text{on}\:\:B_1(0)\setminus\{0\}\:.
\ees
Owing to the hypotheses on $\bn$ and $\bH$, this equation has a distributional sense. Since $\mathcal{L}(\bH)$ is supported on the origin and it belongs to $W^{-1,p}$ for $p<2$, it must be proportional to the Dirac mass $\delta_0$.
% (indeed, as we previously explained, distributions supported at one point are linear combinations of derivatives of the delta mass at that point. But only $\delta_0$ belongs to $W^{-1,p}$ for $p<2$. None other of its higher order derivatives does). 
\noindent
From (\ref{residue}), we deduce that
\bes
\mathcal{L}(\bH)\;=\;-\,\bc_0\,\delta_0\:.
\ees 
Let $\bA\in\C^m$ be the constant vector appearing in the expansion (\ref{locexphi}). Since $\pi_{\bn(0)}\bA=\vec{0}$ (cf. (\ref{expA})), an elementary computation gives
\begin{eqnarray}\label{sheen}
\bA\cdot\bc_0\,\delta_0&=&\pi_{T}\bA\cdot\bc_0\,\delta_0\;=\;-\,\pi_{T}\bA\cdot\mathcal{L}(\bH)\nonumber\\[1ex]
&=&-\,\text{div}\big(\!-\bH\cdot\nabla\pi_T\bA\,+\,\pi_{T}\bA\cdot\star\,(\nabla^\perp\bn\wedge\bH)\big)\nonumber\\[0ex]
&&\quad+\:\nabla\pi_{T}\bA\cdot\Big(\nabla\bH\,-\,3\,\pro\nabla\bH\,+\,\star\,(\nabla^\perp\bn\wedge\bH)\Big)\:,\qquad
\end{eqnarray}
where we have used the fact that $\pi_{T}\bH\equiv\vec{0}$.\\[1ex]
Because $\bA$ is constant and $\nabla\bn\in\bigcap_{p<\infty} L^p$, it follows from (\ref{divpro}) that
$\nabla\pro\bA$ and thus $\nabla\pi_T\bA$ lie in $\bigcap_{p<\infty}L^p$.
Moreover, $\nabla\bH\in\bigcap_{1\le p<2}L^{p}\,$ by hypothesis. Introducing this information into (\ref{sheen}), we note that its right-hand side belongs to $W^{-1,p}$ for all $p<\infty$. Yet, its left-hand side is proportional to the Dirac mass, which does {\it not} belong to any $W^{-1,p}$ for $p\ge2$. We accordingly conclude that $\,\bA\cdot\bc_0\,=\,0.$
Returning to the expansion (\ref{locexphi}) reveals now that 
\bes
\bc_0\cdot\left(\begin{array}{c}\bex(x)\\[1ex]\bey(x)\end{array}\right)\,\simeq\;\bc_0\cdot\bT(x)\;=\;\text{O}(|x|^{1-\epsilon})\qquad\forall\:\:\epsilon>0\:,
\ees
whence
\be\label{dilig}
|x|^{-1}\pi_T(\bc_0)\;\in\;\bigcap_{p<\infty}L^p(B_1(0))\:.
\ee
A direct computation gives
\begin{eqnarray*}
\mathcal{L}\big(\bc_0\log|x|\big)&=&4\pi\,\bc_0\,\delta_0\,+\,\text{div}\big(3\,\pi_T(\bc_0)\,\nabla\log|x|\,+\,\star\,(\nabla^\perp\bn\wedge\bc_0)\log|x|\big)\\[1ex]
&=&-\,4\pi\,\mathcal{L}(\bH)\,+\,\text{div}\big(3\,\pi_T(\bc_0)\,\nabla\log|x|\,+\,\star\,(\nabla^\perp\bn\wedge\bc_0)\log|x|\big)\:.
\end{eqnarray*}
Using the fact that $\nabla\bn\in\bigcap_{p<\infty}L^p$ and (\ref{dilig}) shows that
\bes
\mathcal{L}\bigg(\bH+\frac{\bc_0}{4\pi}\,\log|x|\bigg) \;\in\;\bigcap_{p<\infty}W^{-1,p}\:.
\ees
It is established in \cite{Ri1} that the operator $\mathcal{L}$ is elliptic and in particular that it satisfies $\mathcal{L}^{-1}W^{-1,p}\subset W^{1,p}$. The desired claimed therefore ensues:
\bes
\bH(x)\,+\,\frac{\bc_0}{4\pi}\,\log|x|\;\in\;\bigcap_{p<\infty}W^{1,p}\:,%\subset\bigcap_{\al\in[0,1)}\!\!C^{0,\al}\:,
\ees\\[-8ex]

$\hfill\blacksquare$\\

We end our study of the case $\theta_0=1$ by a slight improvement on the regularity of the Gauss map $\bn$. In the Appendix (cf. (\ref{prof4})), it is shown that the $\bigwedge^{m-2}(\mathbb{S}^{m-1})$-valued Gauss map $\bn$ satisfies a perturbed harmonic map equation, namely
\begin{eqnarray}\label{prof444}
\Delta\bn\,+\,|\nabla\bn|^2\,\bn&=&2\star\text{e}^{\la}\Big[\bAe_1\wedge\pro\,\py\bH\,-\,\bAe_2\wedge\pro\,\px\bH\Big]\nonumber\\[.5ex]
&&\hspace{.5cm}-\;2\star\text{e}^{2\la}\,\vec{h}_{12}\wedge\big(\vec{h}_{11}-\vec{h}_{22}\big)\:.
\end{eqnarray}
Recall that
\bes
|\nabla\bn|\;=\;\text{e}^{-\la}\big|\pro\nabla^2\bp\big|\;=\;\text{e}^{\la}\left|\begin{array}{cc}\vec{h}_{11}&\vec{h}_{12}\\[1ex]\vec{h}_{21}&\vec{h}_{22}\end{array}\right|\:,
\ees
so that $\text{e}^{\la}\vec{h}_{ij}$ inherits the regularity of $\nabla\bn\in\bigcap_{p<\infty}L^p$. Bringing this information and the expansion given in Proposition \ref{logex} into (\ref{prof444}) shows that
\bes
|\Delta\bn|\;\lesssim\;|x|^{-1}+\,\text{terms in\:\:$\bigcap_{p<\infty}L^p$}\;\in\;L^{2,\infty}\:.
\ees
Hence $\nabla^2\bn\in L^{2,\infty}$, and thus $\nabla\bn\in BMO$. 
%Using the $\eps$-regularity (\ref{epsreg}) then implies that in a neighborhood of the origin, 
%\bes
%|\nabla\bn(x)|\;\lesssim\;|x|^{-\epsilon}\qquad\forall\:\:\epsilon>0\:.
%\ees

\subsection{The case $\theta_0\ge2$}

We now return to (\ref{delphi}) in the case when $\theta_0\ge2$. Setting
\be\label{f1f2}
\bF_1\;:=\;\nabla^\perp\bR\,+\nabla\bG\qquad\text{and}\qquad F_2\;:=\;\nabla^\perp S\,+\nabla g\:,
\ee
it reads
\be\label{delphi2}
2\,\Delta\bP\;=\;\,F_2\cdot\nabla\bP\,-\,\bF_1\bul\nabla\bP\:.
\ee
Owing to (\ref{regg}) and (\ref{regsr}), the functions $\bF_1$ and $F_2$ are H\"older continuous of any order $\al\in(0,1)$. It thus makes sense to define the constants
\bes
\bf_1\;:=\;\bF_1(0)\qquad\text{and}\qquad f_2\;:=\;F_2(0)\:.
\ees
They are elements of $\,\R^2\otimes\bigwedge^2(\R^m)\,$ and of $\,\R^2$, respectively. We will in the sequel view $\bf_1$ as an element of $\bigwedge^2(\C^m)$ and $f_2$ as an element of $\C$. \\[.75ex]
For future purposes, let us define $\vec{\Gamma}$ via
\be\label{eqgamma}
\Delta\vec{\Gamma}\;=\;4\,\theta_0\,\Re\big(\bC\,{x}^{\,\theta_0-1}\big)\:\qquad\text{with}\qquad8\,\bC\;:=\;f_2\cdot\bA\,-\,\bf_1\bul\bA\\:,
%2\,\Delta\vec{\Gamma}\;=\;\bf_1\bul\big(\bA\;\overline{x}^{\,\theta_0-1}\big)\,+\,f_2\cdot\big(\bA\;\overline{x}^{\,\theta_0-1}\big)\:,
\ee
where $\bA$ is the constant vector in (\ref{locexphi}). This equation is solved explicitly (up to an unimportant harmonic function):
%To do so, we proceed as before and decompose the constant vector $\bA$
%into its real and imaginary parts, namely $\bA=\bA^{1}+i\bA^{2}$, so that
%\bes
%\bA\;\overline{x}^{\,\theta_0-1}\;=\;|x|^{\theta_0-1}\!\left(\begin{array}{c}\bA^{1}\cos((\theta_0-1)\varphi)\,+\,\bA^{2}\sin((\theta_0-1)\varphi)\\[1ex]\bA^{2}\cos((\theta_0-1)\varphi)\,-\,\bA^{1}\sin((\theta_0-1)\varphi) \end{array}\right)\:.
%\ees
%Since $\,\bf_1:=(\nabla^\perp\bR+\nabla\bG)(0)\,$ and $\,f_2:=(\nabla^\perp S+\nabla g)(0)$, it then follows that
\be\label{didi}
\vec{\Gamma}(x)\;=\;\Re\Big(\bC\,|x|^2{x}^{\,\theta_0-1}\Big)\:.
\ee
%Hence
%\be\label{diva2}
%2\,\Delta\vec{\Gamma}\;=\;|x|^{\theta_0-1}\big(\bu_c\cos((\theta_0-1)\varphi)-\bu_s\sin((\theta_0-1)\varphi)\big)\:,
%\ee
%with the constant vectors\footnote{with the understanding that $\bA=(\bA^{1}\,,\bA^{2})$ and $\bA^\perp=(-\bA^{2}\,,\bA^{1})$.}
%\be\label{ucus}
%\left\{\begin{array}{lcl}\bu_c&=&(\nabla^\perp\bR+\nabla\bG)(0)\bul\bA\,-\,(\nabla^\perp S+\nabla g)(0)\cdot\bA\\[1ex]
%\bu_s&=&(\nabla^\perp\bR+\nabla\bG)(0)\bul\bA^\perp\,-\,(\nabla^\perp S+\nabla g)(0)\cdot\bA^\perp\:.
%\end{array}\right.
%\ee
%\be\label{ucus}
%\bu_c\;:=\;\bf_1\bul\bA\,+\,f_2\cdot\bA\qquad\text{and}\qquad
%\bu_s\;:=\;\bf_1\bul\bA^\perp\,+\,f_2\cdot\bA^\perp\:.
%\ee
%Solving (\ref{diva2}) yields\footnote{as before, $\nabla$ is understood as $\,\px+i\,\py$.}
%\be\label{gaby1}
%\nabla\vec{\Gamma}(x)\;=\;\dfrac{1}{8\theta_0}\Big((\bu_c-i\bu_s)\,x^{\theta_0}\,+\,\theta_0(\bu_c+i\bu_s)|x|^2\,\overline{x}^{\,\theta_0-2}\Big)\:.
%\bC^{*}x^{\theta_0}\,+\,\theta_0\,\bC\,|x|^2\,\overline{x}^{\,\theta_0-2}\:,
%\ee
Note next that (\ref{delphi2}) and (\ref{eqgamma}) give
\be\label{ipam}
2\,\Delta(\bp-\vec{\Gamma})\;=\;(F_2-f_2)\cdot\nabla\bp\,-\,(\bF_1-\bf_1)\bul\nabla\bp\;+\;\text{e}^{\la}\big[f_2\cdot\bT\,-\,\bf_1\bul\bT\big]\:,
\ee
%\begin{eqnarray*}
%-\,2\,\Delta(\bp-\vec{\Gamma})&=&(\bF_1-\bf_1)\bul\nabla\bp\,-\,(F_2-f_2)\cdot\nabla\bp\\[.7ex]
%&&\quad+\;\text{e}^{\la}\big[\bf_1\bul\bT\,+\,f_2\cdot\bT\big]
%\end{eqnarray*}
where we have used the representation (\ref{locexphi}). We have seen (compare (\ref{locexphi}) to (\ref{locexdelphi})) that $\pj(\text{e}^{\la}\bT)=\text{e}^{\la}\bQ_j$, where $\bQ_j$ belongs to $L^p$ for all $p<\infty$. Differentiating (\ref{ipam}) throughout with respect to $x_j$ gives
\begin{eqnarray}\label{ipam2}
2\,\Delta\pj(\bp-\vec{\Gamma})&=&\pj F_2\cdot\nabla\bp\,-\,\pj\bF_1\bul\nabla\bp\;+\;\text{e}^{\la}\big[f_2\cdot\bQ_j\,-\,\bf_1\bul\bQ_j\big]\nonumber\\[.75ex]
&&\quad+\:\:(F_2-f_2)\cdot\nabla\pj\bp\,-\,(\bF_1-\bf_1)\bul\nabla\pj\bp\:.
\end{eqnarray}
Since $\pj\bF_1$, $\pj F_2$, and $\bQ_j$ belong to $L^p$ for every finite $p$, while $|\nabla\bp|\simeq\text{e}^{\la}\simeq|x|^{\theta_0-1}$, we may apply Proposition \ref{CZpondere}-(i) to the first three summands on the right-hand side of (\ref{ipam2}). Moreover, $\,|\bF_1(x)-\bf_1|+|F_2(x)-f_2|\lesssim|x|^{\al}$\, for all $\al\in(0,1)$\, while $\,|\nabla\pj\bp(x)|\simeq|x|^{\theta_0-2}$\,, so that the last two summands on the right-hand side of (\ref{ipam2}) fit within the frame of Corollary \ref{CZcoro}. Accordingly, 
\be\label{diva1}
\nabla\pj\big(\bp-\vec{\Gamma}\big)(x)\;=\;\vec{P}_j(\overline{x})\,+\,\text{e}^{\la(x)}\vec{U}_j(x)\:,
\ee
where $\vec{P}_j$ is a polynomial of degree at most $(\theta_0-1)$, and $\vec{U}_j(x)=\text{O}(|x|^{1-\epsilon})$, for every $\epsilon>0$.\\
One sees in (\ref{didi}) that $\nabla\pj\vec{\Gamma}(x)=\text{O}(|x|^{\theta_0-1})$. Hence, from (\ref{diva1}) and the fact that $|\nabla\pj\bp|(x)\simeq|x|^{\theta_0-2}$, it follows that the polynomial $\vec{P}_j$ contains exactly one monomial of degree of $(\theta_0-2)$ and one monomial of degree of $(\theta_0-1)$. More precisely, identifying the representation (\ref{locexdelphi}) with (\ref{diva1}) yields
%another of  $\,\vec{P}(\overline{x})=\vec{A}\;\overline{x}^{\,\theta_0-1}+\vec{B}\;\overline{x}^{\,\theta_0}\,$ for two $\C^m$-valued vectors $\bA$ and $\bA$ (and of course $\bA$ is as in (\ref{locexphi})). Accordingly, we have the representation
%\bes
%\nabla\bp(x)\;=\;\vec{A}\;\overline{x}^{\,\theta_0-1}\,+\,\vec{B}\;\overline{x}^{\,\theta_0}\,+\,\vec{C}^*x^{\theta_0}\,+\,\theta_0\,\vec{C}\,|x|^2\,\overline{x}^{\,\theta_0-2}\,+\,\text{e}^{\la(x)}\vec{U}(x)\:,
%\ees
%where $\vec{C}=-\frac{1}{8\theta_0}(\vec{u}_c-i\vec{u}_s)$ and $\vec{C}^*$ is its complex conjugate. From the latter we obtain
\begin{eqnarray}\label{ddphi0}
\nabla^2\bp(x)&=&\left(\begin{array}{cc}-\,\Re&\Im\\[.5ex]\Im&\Re\end{array}\right)\Big(\theta_0(1-\theta_0)\bA\;{x}^{\theta_0-2}\,-\,\theta_0(1+\theta_0)\vec{B}\;{x}^{\theta_0-1}\Big)\nonumber\\[-1ex]
%&&+\:(\theta_0-1)\left(\begin{array}{cc}0&1\\1&0\end{array}\right)\Im\Big(\bA\;\overline{x}^{\,\theta_0-2}+\vec{B}\;\overline{x}^{\,\theta_0-1}+\,\bC\,|x|^2\,\overline{x}^{\,\theta_0-3}\Big)\nonumber\\[1ex]
&&\hspace{2cm}+\:\nabla^2\vec{\Gamma}(x)\,+\,\text{e}^{\la(x)}\vec{U}(x)\:,
\end{eqnarray}
where $\bB\in\C^m$ is a constant vector and $\vec{U}(x)=\text{O}(|x|^{1-\epsilon})$ for all $\epsilon>0$. The constant vector $\bA$ is as in (\ref{locexphi}). \\
We deduce from (\ref{ddphi0}) and (\ref{didi}) the expansion (recall that $|\nabla\bp(0)|={0}=|\bp(0)|$)\,:
%\bes
%\nabla\bp(x)\;=\;\left(\begin{array}{l}\Re\\[1.5ex]\Im\end{array}\right)\big(\bA\;\overline{x}^{\,\theta_0-1}+\,\theta_0^{-1}{\bC}\,{x}^{\theta_0}\,+\,\bC\,|x|^2{x}^{\,\theta_0-2}\big)\;+\;\text{O}\big(|x|^{\theta_0+1-\epsilon}\big)\:,
%\ees
%and
\be\label{xbp}
\bp(x)\;=\;\Re\Big(\bA\;x^{\theta_0}+\bB\,x^{\theta_0+1}+\bC\,|x|^2{x}^{\,\theta_0-1}\Big)+\,|x|^{\theta_0-1\,}\vec{\xi}(x)\:,
\ee
where
\bes
\vec{\xi}(x)=\text{O}(|x|^{3-\epsilon})\:\:,\:\:\:\:\:\:\nabla\vec{\xi}(x)=\text{O}(|x|^{2-\epsilon})\:\:,\:\:\:\:\:\:\nabla^2\vec{\xi}(x)=\text{O}(|x|^{1-\epsilon})\qquad\forall\:\:\epsilon>0\:.
\ees
Moreover, as $\,2\,\text{e}^{2\la}\bH=\Delta\bp\,$, the representation (\ref{ddphi0}) along with (\ref{eqgamma}) gives the local asymptotic expansion
\be\label{h2}
\text{e}^{\la(x)}\bH(x)\;=\;f(x)\,\Re\Bigg[\,\bC\bigg(\dfrac{{x}}{|x|}\bigg)^{\!\theta_0-1\,}\Bigg]+\,\text{O}(|x|^{1-\epsilon})\:,
\ee
where $\bC$ is as above, and $\,f(x):=2\,\theta_0|x|^{\theta_0-1}\text{e}^{-\la(x)}$, which is known to have a positive limit at the origin. 
%\begin{eqnarray}\label{h2}
%2\,\text{e}^{2\la(x)}\bH(x)&\equiv&\text{Tr}\,\nabla^2\bp(x)\\[1ex]
%&=&|x|^{\theta_0-1}\big(\bu_c\cos((\theta_0-1)\varphi)-\bu_s\sin((\theta_0-1)\varphi)\big)\,+\,\text{O}(|x|^{\theta_0-\epsilon})\nonumber\:.
%\end{eqnarray}
This shows that $\text{e}^{\la(x)}\bH(x)$ is a bounded function. However, it ``spins" as $x$ approaches the origin: its limit need not exist ; and, if it does exist, it must be zero (i.e. $\bC=\vec{0}$). This possibility is studied in details below. \\

We close this section by proving that $\nabla^2\bn\in L^{2,\infty}$ and that $\nabla\bn\in L^\infty$. We have seen that $\text{e}^{\la}\bH$ is bounded. Applying standard elliptic techniques to (\ref{wildiv}) then yields that $\,|x|\,\text{e}^{\la}\nabla\bH$ is bounded as well, and hence that $\text{e}^{\la}\nabla\bH\in L^{2,\infty}$. Going back to the perturbed harmonic map equation (\ref{prof444}) satisfied by the Gauss map $\bn$, and using the fact that $\text{e}^{\la}\vec{h}_{ij}$ inherits the regularity of $\nabla\bn\in\bigcap_{p<\infty}L^p$, we deduce that $\Delta\bn$ lies in $L^{2,\infty}$, and therefore indeed that $\nabla^2\bn\in L^{2,\infty}$. In particular, this implies that $\nabla\bn\in BMO$. It is actually possible to show that $\nabla\bn\in L^\infty(B_1(0))$. To see this, we first note that (\ref{xbp}) yields
\bes
\nabla\bp(x)\;=\;\left(\begin{array}{c}\Re\\-\,\Im\end{array}\right)\big(\theta_0\bA\,x^{\theta_0-1}\big)\,+\,\nabla\big(|x|^{\theta_0-1}\vec{\xi}(x)\big)\,+\,\text{O}(|x|^{\theta_0})\:.
\ees
Since $\,\pro\nabla\bp\equiv0$, the latter and the estimates on $\vec{\xi}$ give
\bes
|\pi_{\bn(x)}\bA|\;=\;\text{O}(|x|)\:.
\ees
A quick inspection of the identity (\ref{ddphi0}) then reveals that
\bes
\big|\pro\nabla^2\bp(x)\big|\:\lesssim\:\pro(\bA)\,|x|^{\theta_0-2}\;=\;\text{O}(|x|^{\theta_0-1})\:.
\ees
Combining this to (\ref{modulegradn}) gives thus that $\nabla\bn$ is bounded across the singularity.

\subsubsection{When $\text{e}^{\la}\bH$ has a limit at the origin}\label{lim0}

We shall now consider the case when $\lim_{|x|\searrow0}\text{e}^{\la(x)}\bH(x)$ exists. Then, as seen in (\ref{h2}), we automatically have $\bC=\vec{0}$, and accordingly
\be\label{goodH}
\text{e}^{\la(x)}\bH(x)\;=\;\text{O}(|x|^{1-\epsilon})\qquad\forall\:\:\epsilon>0\:.
\ee
We draw the reader's attention on the fact that when $\theta_0=2$, the latter implies
$\bH(x)=\text{O}(|x|^{-\epsilon})$ for all $\epsilon>0$.\\[-1ex]

We now show that $\bC=\vec{0}$ implies that the constants $\bf_1$ and $f_2$ are also trivial. To see this, recall (\ref{f1f2}) and (\ref{truxy}), namely
\bes
\bf_1\;=\;\big(\nabla^\perp\bR+\nabla\bG\big)(0)\;=\;\big(\bL\wedge\nabla^\perp\bP\,-\,2\,\bH\wedge\nabla\bp)(0)
\ees
and
\be\label{cst2}
\bf_2\;=\;\big(\nabla^\perp S+\nabla g\big)(0)\;=\;\big(\bL\cdot\nabla^\perp\bP)(0)\:.
\ee
From $|\nabla\bP|\simeq\text{e}^{\la}$ and (\ref{goodH}), we know that $\,(\bH\wedge\nabla\bP)(0)=\vec{0}$. To obtain $\bf_1=0=f_2$, it thus suffices to show that $\,\lim_{|x|\searrow0}\text{e}^{\la(x)}\bL(x)=\vec{0}$. This is what we shall do. 
Using a standard argument from elliptic analysis (identical to that enabling to deduce (\ref{es3}) from (\ref{es2})), it follows from (\ref{goodH}) that $\,\text{e}^{\la(x)}\nabla\bH(x)=\text{O}(|x|^{-\epsilon})$ for all $\epsilon>0$. Bringing this information into (\ref{a1}) and (\ref{a2}), along with the fact that $\nabla\bn\in L^p$ for all $p<\infty$, shows now that $\,\text{e}^{\la}\nabla\bL\in L^p$ for all finite $p$. The Hardy-Sobolev inequality (\ref{estimL}) with $\theta_0-1$ in place of $\theta_0$ implies in particular that $\,|x|^{-1}\text{e}^{\la}\bL\in L^2$. Owing to (\ref{cst2}), the limit $\,f_2=\lim_{|x|\searrow0}(\bL\cdot\nabla^\perp\bp)(x)\,$ exists. Yet, we have seen the function $\,|x|^{-1}(\bL\cdot\nabla^\perp\bp)(x)\,$ is square-integrable near the origin. This is only possible if $f_2=0$. We proceed {\it mutatis mutandis} to show that $\,\lim_{|x|\searrow0}(\bL\wedge\nabla^\perp\bp)(x)=\vec{0}$, thereby yielding $\bf_1=\vec{0}$. \\

%For future purposes, it is useful to compute $\,\pro\,\text{div}\,\pro\nabla\pj\bp$. Putting $\bC=\vec{0}$ in (\ref{ddphi}) and carrying an elementary computation gives
%\begin{eqnarray}\label{divx}
%\pro\,\text{div}\,\pro\nabla\px\bp&=&(\theta_0-1)|x|^{\theta_0-2}\cos((\theta_0-2)\varphi)\,\pro\big(\px\pro\bA^{1}+\py\pro\bA^{2})\nonumber\\[1ex]
%&+&\!\!\!(\theta_0-1)|x|^{\theta_0-2}\sin((\theta_0-2)\varphi)\,\pro\big(\px\pro\bA^{2}-\py\pro\bA^{1})\:.\nonumber\\
%\end{eqnarray}
%Because $\bA$ is a constant vector, there holds\footnote{implicit summations on repeated indices is understood, with $1\le\al\le m-2$ and $1\le l\le 2$. Refer to the Appendix XXX for the notation.} 
%\begin{eqnarray*}
%\pro\,\pj\pro\bA_k&=&-\,\pro\,\pj\pi_{T}\bA_k\;\;=\;\;-\,\big(\bn_\al\cdot\pj\pi_T\bA_k\big)\,\bn_\al\nonumber\\[.75ex]
%&=&\big(\bA_k\cdot\pi_T(\pj\bn_\al)\big)\,\bn_\al\;\;=\;\;-\,\text{e}^{\la}(\bA_k\cdot\bAe_l)\,\vec{h}_{jl}\:.
%\end{eqnarray*}
%Whence the estimate
%\bes
%\big|\pro\,\pj\pro\bA_k\big|\;\leq\;|\bA_k|\,\text{e}^{\la}\big(|\vec{h}_{j1}|+|\vec{h}_{j2}|\big)\;\lesssim\;\text{e}^{-\la}\big|\pro\nabla^2\bp\big|\;=\;|\nabla\bn|\:.
%\ees
%Brought into (\ref{divx}), the latter and (\ref{regn2}) and  yields
%\be\label{trip1}
%\big|\text{e}^{-\la}\,\pro\,\text{div}\,\pro\nabla\px\bp\big|\;\lesssim\;|\nabla\bn|\;\in\;L^p\:,\qquad\forall\:\:p<\infty\:.
%\ee
%{\it Mutatis mutandis}, we obtain
%\be\label{trip2}
%\big|\text{e}^{-\la}\,\pro\,\text{div}\,\pro\nabla\py\bp\big|\;\lesssim\;|\nabla\bn|\;\in\;L^p\:,\qquad\forall\:\:p<\infty\:.
%\ee
In the Appendix (cf. (\ref{prof4})), it is shown that the $\bigwedge^{m-2}(\mathbb{S}^{m-1})$-valued Gauss map $\bn$ satisfies a perturbed harmonic map equation:
\begin{eqnarray}\label{prof44}
\Delta\bn\,+\,|\nabla\bn|^2\,\bn&=&2\star\text{e}^{\la}\Big[\bAe_1\wedge\pro\,\py\bH\,-\,\bAe_2\wedge\pro\,\px\bH\Big]\nonumber\\[.5ex]
&&\hspace{.5cm}-\;2\star\text{e}^{2\la}\,\vec{h}_{12}\wedge\big(\vec{h}_{11}-\vec{h}_{22}\big)\:.
\end{eqnarray}
Moreover, as explained at the end of Section \ref{tet1}, $\text{e}^{\la}\vec{h}_{ij}$ inherits the regularity of $\nabla\bn$. Namely,
\be\label{wipe0}
|\nabla\bn|\;=\;\text{e}^{-\la}\big|\pro\nabla^2\bp\big|\;=\;\text{e}^{\la}\left|\begin{array}{cc}\vec{h}_{11}&\vec{h}_{12}\\[1ex]\vec{h}_{21}&\vec{h}_{22}\end{array}\right|\:.
\ee
We then deduce from (\ref{prof44}) the estimate
\bes
|\Delta\bn|\;\lesssim\;|\nabla\bn|^2\,+\,\text{e}^{\la}\big|\pro\nabla\bH\big|\:.
\ees
In proving (\ref{cst2}), we have seen that $\text{e}^{\la(x)}\nabla\bH(x)=\text{O}(|x|^{-\epsilon})$ for all $\epsilon>0$. Furthermore, from (\ref{regn2}), we know that $\nabla\bn$ has as much integrability as we please. The right-hand side of the equation (\ref{prof44}) thus belongs to $L^p$ for all finite $p$, thereby showing that (\ref{prof44}) is subcritical, and thus yielding
\be\label{regd2n}
\nabla^2\bn\:\in\:L^p(B_1(0))\qquad\forall\:\:p<\infty\:.
\ee

When $\theta_0=2$, the argument comes to a halt at this point. However, if $\theta_0\ge3$, we note from (\ref{sys1}) that the regularity of $\nabla g$ and $\nabla\bG$ improves to $W^{2,p}$ for all $p<\infty$. Introducing this new information along with (\ref{srinsobolev}) and (\ref{regd2n}) into (\ref{sysSR}) shows that $\nabla S$ and $\nabla\bR$ are elements of $W^{2,p}$ for all finite $p$. Hence the functions $\bF_1$ and $F_2$ defined in (\ref{f1f2}) now lie in $W^{2,p}$ for all $p<\infty$. Moreover, we have seen that they both vanish at the origin, so that
\bes
|\bF_1(x)|\,+\,|F_2(x)|\;\lesssim\;|x|^{1+\al}\qquad\forall\:\:\al\in(0,1)\:.
\ees
Returning to (\ref{delphi2}) and applying Corollary \ref{CZcoro2} with $\mu=\text{e}^{\la}$, $a=\theta_0-1$, $n=1=J$, and $r=\al$ gives
\be\label{rlb}
\nabla^2\bp(x)\;=\;\nabla\vec{P}(\overline{x})\,+\,\text{e}^{\la(x)}\bV(x)\:,
\ee
where $\vec{P}$ is a polynomial of degree at most $(\theta_0+1)$ and $\bV(x)=\text{O}(|x|^{2-\epsilon})$ for all $\epsilon>0$. Note that since $|\nabla\bp(x)|\simeq|x|^{\theta_0-1}$ by hypothesis, $\vec{P}$ has no terms of degree smaller than $(\theta_0-1)$. Being a (nonlinear) polynomial of the variable $\overline{x}$, the polynomial $P$ has traceless gradient. Whence we deduce from (\ref{rlb}) that
\bes
2\,\text{e}^{\la(x)}\bH(x)\;\equiv\;\text{e}^{-\la(x)}\,\text{Tr}\,\nabla^2\bp(x)\;=\;\text{Tr}\,\bV(x)\;=\;\text{O}(|x|^{2-\epsilon})\:.
\ees
In particular, when $\theta_0=3$, we arrive at $\,\bH(x)=\text{O}(|x|^{-\epsilon})$\, for all $\epsilon>0$.\\[1.5ex]
As we did in the paragraph following (\ref{cst2}), we deduce from the asymptotics of $\bH$ those of $\nabla\bH$, namely $\,\text{e}^{\la(x)}\nabla\bH(x)=\text{O}(|x|^{1-\epsilon})$. To further improve on the regularity of the mean curvature, we may differentiate (\ref{wildiv}) throughout with respect to $x_j$. We obtain an equation for $\pj\bH$ in divergence form valid on $B_1(0)\setminus\{0\}$. The coefficients involve $\bn$, its first and its second derivatives, all of which belong to $L^p$ for every $p<\infty$. As previously done, we can now deduce local asymptotics for $\nabla\pj\bH$ from those of $\pj\bH$. More precisely, $\text{e}^{\la(x)}\nabla^2\bH(x)=\text{O}(|x|^{-\epsilon})$ for all $\epsilon>0$. Since $\nabla\la(x)\lesssim|x|^{-1}$, we find
\be\label{wipe1}
\nabla\big(\text{e}^{\la}\pj\bH\big)\;\equiv\;\text{e}^{\la}\big(\nabla\pj\bH\,+\,\nabla\la\,\pj\bH\big)\:\lesssim\:|x|^{-\epsilon}\qquad\forall\:\:\epsilon>0\:.
\ee
In general for a vector $\bV$, there holds $\,\pro\bV=\bn\res(\bn\res\bV)$. Hence, 
\be\label{divpro}
\big|\nabla\pro\bV\big|\;\lesssim\;|\nabla\bV|\,+\,|\bV|\,|\nabla\bn|\:.
\ee
Since $\nabla^2\bn\in\bigcap_{p<\infty}L^p$, (\ref{wipe1}) gives in particular
\be\label{wipe2}
\big|\nabla\big(\text{e}^{\la}\pro\,\pj\bH\big)\big|\;\lesssim\;|x|^{-\epsilon}\qquad\forall\:\:\epsilon>0\:.
\ee
In addition,
\bes
\nabla\bek\;\equiv\;\text{e}^{-\la}\big(\pk\nabla\bp\,-\,\nabla\la\,\pk\bp)\:\lesssim\:|x|^{-1}\:.
\ees
Combining the latter to (\ref{wipe2}) shows that
\bes
\Big|\nabla\big(\text{e}^{\la}\,\bek\wedge\pro\,\pj\bH\big)\Big|\;\lesssim\;|x|^{-\epsilon}\qquad\forall\:\:\epsilon>0\:.
\ees
We now introduce this information along with (\ref{wipe0}) and (\ref{regd2n}) into (\ref{prof44}) to obtain
\bes
\big|\Delta\nabla\bn\big|\;\;\lesssim\;\;\big|\nabla\bn\big|^3+\big|\nabla\bn\big|\,\big|\nabla^2\bn\big|\,+\,|x|^{-\epsilon}\:\in\:\bigcap_{p<\infty}L^p\:,
\ees
so that
\bes
\nabla^3\bn\,\in\,L^p(B_1(0))\qquad\forall\:\:p<\infty\:.
\ees
Note also that
\bes
|\nabla^3\bp|\;=\;\big|\nabla^2|\nabla\bp|\big|\;\simeq\;|\nabla^2\text{e}^{\la}|\;=\;\text{e}^{\la}\,\big|\nabla^2\la+(\nabla\la)^2\big|\:.
\ees
The expansion (\ref{rlb}) thus gives
\bes
|x|^2|\nabla^2\la|\:\,\lesssim\:\,\big||x|\nabla\la\big|^2+\,|x|^2|\text{e}^{-\la}\nabla^2\vec{P}|\,+\,|x|^2\text{e}^{-\la}\big|\nabla(\text{e}^{\la}\bV)\big|\:.
\ees
We know that $|x|\la$ is a bounded function. Moreover, since $\text{e}^{-\la}\simeq|x|^{1-\theta_0}$ and $\vec{P}$ is a polynomial containing no terms of degree less than $(\theta_0-1$), we get
\bes
|x|^2|\nabla^2\la|\:\,\leq\:\,C\,+\,|x|^2\text{e}^{-\la}\big|\nabla(\text{e}^{\la}\bV)\big|\:\qquad\text{for some constant $C$}. 
\ees
Corollary \ref{CZcoro2} states that $\,|x|^{\epsilon-1}\text{e}^{-\la}\nabla(\text{e}^{\la}\bV)$ belongs to $L^p$ for all $p<\infty$ and all $\epsilon>0$. However, by tracking the way this estimate is obtained, it is not difficult to verify that $\,|x|^{2}\text{e}^{-\la}\nabla(\text{e}^{\la}\bV)$ tends to zero as $x$ moves towards the origin. Hence,
\bes
|x|^2\,\nabla^2\la(x)\;\in\;L^\infty(B_1(0))\:.
\ees

\bigskip
This procedure continues on. As $\theta_0$ increases, so does the regularity of $g$ and $\bG$, thereby improving that of $S$ and $\bR$. Repeating the above argument through Corollary \ref{CZcoro2}\footnote{namely, every time $\theta_0$ increases by an increment of one, so does the parameter $n$ in Corollary \ref{CZcoro2}, and we increase accordingly the parameter $J$ by one. The procedure allows up to $n=\theta_0-2$ and $J=\theta_0-1$.} yields that
\bes
\nabla^{j+1}\bp(x)\;=\;\nabla^{j}\vec{P}(\overline{x})\,+\,\text{e}^{\la(x)}\bV_j(x)\:,\qquad\forall\:\:j\in\{0,\ldots,\theta_0-1\}\:,
\ees
where $\vec{P}$ is a polynomial of degree at most $(2\theta_0-2-j)$  and $\bV_j(x)=\text{O}(|x|^{\theta_0-j-\epsilon})$ for all $\epsilon>0$. In particular, using the hypothesis that $|\nabla\bp(x)|\simeq|x|^{\theta_0-1}$, we deduce that
\bes
\overline{\vec{P}(\overline{x})}\;=\;\theta_0\bA\,{x}^{\,\theta_0-1}+\,\sum_{j=1}^{\theta_0-1}\,(\theta_0+j)\bA_j\,{x}^{\,\theta_0-1+j}\:,
\ees
where $\bA$ is as in (\ref{locexphi}) and $\bA_j\in\C^m$ are constant vectors. Altogether, since $|\bp(0)|={0}=|\nabla\bp(0)|$ by hypothesis, we obtain the representation
\be\label{baddphi}
\bp(x)\;=\;\Re\bigg(\bA\,x^{\theta_0}+\,\sum_{j=1}^{\theta_0-1}\bA_j\,x^{\theta_0+j}\bigg)\,+\,\vec{\zeta}(x)\:,
\ee
where the function $\vec{\zeta}$ satisfies
\be\label{zetata}
|\nabla^j\vec{\zeta}(x)|\;=\;\text{e}^{\la(x)}|\bV_{j-1}(x)|\;=\;\text{O}(|x|^{2\theta_0-j-\epsilon})\qquad\forall\:\:\epsilon>0\:\:,\:\:\:j\in\{0,\ldots,\theta_0\}\:.
\ee
The last item from Corollary \ref{CZcoro2} also gives
\be\label{intermeme}
\big|\nabla^{\theta_0-1}\Delta\vec{\zeta}(x)\big|\:\lesssim\:|x|^{\theta_0-1-\epsilon}\qquad\forall\:\:\epsilon>0\:.
\ee
%Note that the conformality condition on $\bp$ yields
%\bes
%\big\langle\bA_j\,,\bA_k\big\rangle_{\C^{2m}}=\;0\qquad\forall\:\:\:\:0\le j,k\le\theta_0-1\:\:\qquad\text{with}\:\:\:\:\bA:=\bA\:.
%\ees
From the above we obtain inductively that
\be\label{lambdada}
|x|^j\nabla^j\la(x)\;\in\;L^{\infty}(B_1(0))\qquad\forall\:\:j\in\{0,\ldots,\theta_0-1\}\:.
\ee
Moreover,
\bes
\bH(x)\;=\;\dfrac{1}{2}\,\text{e}^{-2\la(x)}\Delta\vec{\zeta}(x)\:.
\ees
Combining (\ref{zetata})-(\ref{lambdada}) then yields
\be\label{betterH}
\nabla^j\bH(x)\;=\;\text{O}(|x|^{-j-\epsilon})\qquad\forall\:\:\epsilon>0\:\:\:\:\text{and}\:\:\:\:j\in\{0,\ldots,\theta_0-1\}\:.
\ee
In particular, we have $\,\nabla\bH(x)=\text{O}(|x|^{-1-\epsilon})\in\bigcap_{p<2}L^p$. With this fact at our disposal and (\ref{regn}), we call upon Proposition \ref{logex} and obtain the local expansion
\be\label{strausskahn}
\bH(x)\,+\,\frac{\bc_0}{4\pi}\,\log|x|\;\in\bigcap_{p<\infty}W^{1,p}(B_1(0))\:,
\ee
where $\bc_0$ is the residue defined in (\ref{residue}). \\[1ex]
In addition, the above procedure implies
\be\label{highn}
\nabla^{\theta_0}\bn\,\in\,L^p(B_1(0))\qquad\forall\:\:p<\infty\:.
\ee
To obtain pointwise information about the Gauss map, we use a ``higher order" version of the $\eps$-regularity which appears in \cite{Ri1} (cf. Theorem I.5) along with the same technique as in the proof of Lemma \ref{delta}. Under the hypothesis (\ref{epsregcond}), there holds
\be\label{epsreg2}
r^j\!\sup_{|x|=r}|\nabla^j\bn(x)|\:\lesssim\:\Vert\nabla\bn\Vert_{L^2(B_{2r}(0)\setminus B_{r/2}(0))}\qquad\:\:\forall\:\:\;r\in(0\,,1/2)\:\:,\:\:j\in\mathbb{N}^*\:.
%\;\lesssim\;\dfrac{1}{\sigma^j}\,\Vert\nabla\bn\Vert_{L^2(B_{2\sigma})}
\ee
It then follows easily from (\ref{highn}) that
\be\label{ocean}
|\nabla^j\bn(x)|\;\lesssim\;|x|^{\theta_0-j-\epsilon}\qquad\forall\:\:\epsilon>0\:,\:\:j\in\{0,\ldots,\theta_0\}\:.
\ee

\medskip
%\textbf{On the Constant Vectors $\bA_j$\,.}\\
We have seen that
\be\label{baddphi2}
\bp(x)\;=\;\Re\bigg(\bA\;x^{\theta_0}+\,\sum_{j=1}^{\theta_0-1}\bA_j\,x^{\theta_0+j}\bigg)\,+\,\vec{\zeta}(x)\:.
\ee
From $\,\pro\nabla\bp\equiv\vec{0}$, we obtain after a few computations
\begin{eqnarray*}
(\theta_0+1)\,\pro\bA_1&=&x^{1-\theta_0}\pro\big(\partial_{x_1x_1}\bp-i\,\partial_{x_1x_2}\bp\big)\,+\,x^{1-\theta_0}\pro\big(\partial_{x_1x_1}\vec{\zeta}-i\,\partial_{x_1x_2}\vec{\zeta}\,\big)\\[.5ex]
&&+\:\;x^{-\theta_0}\pro\big(\px\vec{\zeta}+i\,\py\vec{\zeta}\big)\,-\,\sum_{j=2}^{\theta_0-1}j(\theta_0+j)\,\pro\bA_j\,x^{j-1}\:.
\end{eqnarray*}
Hence,
\bes
|\pro\bA_1|\;\;\lesssim\;\;|x|^{1-\theta_0}|\pro\nabla^2\bp|\,+\,|x|^{1-\theta_0}|\pro\nabla^2\vec{\zeta}|\,+\,|x|^{-\theta_0}|\pro\nabla\vec{\zeta}|\,+\,\text{O}(|x|)\:.
\ees
Using (\ref{zetata}) and the fact that $\,|\nabla\bn|=\text{e}^{-\la}|\pro\nabla^2\bp|\simeq|x|^{1-\theta_0}|\pro\nabla^2\bp|\,$ yields
\bes
|\pi_{\bn(x)}\bA_1|\;\lesssim\;|\nabla\bn(x)|\,+\,\text{O}(|x|^{\theta_0-1-\epsilon})\:\:\:\qquad\forall\:\:\epsilon>0\:.
\ees
Then, from (\ref{ocean}), we deduce
\bes
\pi_{\bn(0)}\bA_1\;=\;\vec{0}\:.
\ees
This process is repeated by taking successive derivatives and using (\ref{zetata}). Namely, for each $k\in\{0,\ldots,\theta_0-1\}$\,:
\begin{eqnarray*}
|\pro\bA_k|&\lesssim&\sum_{j=1}^{k+1}\,|x|^{j-k-\theta_0}\Big[|\pro\nabla^{j}\bp|+|\pro\nabla^{j}\vec{\zeta}|\Big]\,+\,\text{O}(|x|)\\[.5ex]
&\lesssim&\sum_{j=1}^{k+1}\,|x|^{j-k-\theta_0}|\pro\nabla^{j}\bp|\,+\,\text{O}\big(|x|^{\theta_0-k-\epsilon}\big)\qquad\quad\forall\:\:\epsilon>0\:.
\end{eqnarray*}
Because $|\nabla\bn|=\text{e}^{-\la}|\pro\nabla^2\bp|$, we obtain through a simple calculation that
\bes
|\pro\nabla^j\bp|\:\;\lesssim\:\;\text{e}^{\la}\,\sum_{q=1}^{j-1}\,|\nabla^{j-1-q}\la|\,|\nabla^{q}\bn|\;=\;\text{O}\big(|x|^{2\theta_0-j-\epsilon}\big)\qquad\quad\forall\:\:\epsilon>0\:\:,
\ees
where we have used (\ref{lambdada}) and (\ref{ocean}). Combining altogether the latter two estimates yields
\bes
|\pi_{\bn(x)}\bA_k|\:=\:\text{O}\big(|x|^{\theta_0-k-\epsilon}\big)\qquad\quad\forall\:\:\epsilon>0\:\:,\quad k\in\{0,\ldots,\theta_0-1\}\:.
\ees
Hence in particular,
\be\label{pronull}
\pi_{\bn(0)}\bA_k\;=\;\vec{0}\qquad\quad\forall\:\:k\in\{0,\ldots,\theta_0-1\}\:.
\ee

\medskip
We have seen in general (cf. (\ref{condA1})-(\ref{condA2})) that $\{\bA^{1},\bA^{2}\}$ forms an orthogonal basis of a plane through the origin (this plane may {\it a priori} not be viewed as {\it the} tangent plane $T_0\Sigma$ which need not exist). Moreover, from the expansion (\ref{locexphi}), we obtain that
\bes
\star\,\bn\;:=\;\dfrac{\px\bp\wedge\py\bp}{|\px\bp\wedge\py\bp|}\;\,\simeq\;\,-\,\dfrac{\bA^{1}\wedge\bA^{2}}{|\bA^{1}\wedge\bA^{2}|}\:,
\ees
so that the Gauss map is well-defined at the origin. Hence,
%This is yet one additional clue indicating that in the case when $\lim_{x\rightarrow0}\text{e}^{\la(x)}\bH(x)$ exists, our parametrization $\bp$ is a bad one: it seems to cover multiple times the same surface. From the point of view of manifolds, there is thus no tangent plane at the origin. But from the point of view of varifolds, the surface has at the origin a tangent $\theta_0$-plane. Fortunately, as seen in the above computation and in Proposition \ref{Th4}-(iii), the Gau\ss\;map, and hence the mean curvature vector, are ``blind" to these bad parametrizations. We will in this section try to quantify in more precise terms how close the parametrization (\ref{badphi}) is to being a bad one. \\
any constant vector $\bV\in\R^m$ has a representation
\bes
\bV\;=\;\dfrac{1}{a}(\bV\cdot\bA^1)\,\bA^1\,+\,\dfrac{1}{a}(\bV\cdot\bA^2)\,\bA^2\,+\,\pi_{\bn(0)}\bV\:,
\ees
where $\,a:=|\bA^1|=|\bA^2|$. For two vectors $\vec{U}=\vec{U}^1+i\,\vec{U}^2$ and $\bV=\bV^1+i\,\bV^2$ in $\R^2\otimes\R^m\simeq\C^m$, we define the product
\bes
\big\langle\vec{U}\,,\bV\big\rangle_{\C^{2m}}\;:=\;\big(\vec{U}^1\cdot\vec{V}^1-\vec{U}^2\cdot\vec{V}^2\big)\,+\,i\,\big(\vec{U}^1\cdot\vec{V}^2+\vec{U}^2\cdot\vec{V}^1\big)\;\in\;\C\:.
\ees
The conformality condition applied to (\ref{baddphi2}) then easily yields that for every $\,s\in\{0,\ldots,\theta_0-1\}\,$ there holds
\be\label{moogli}
2\,\sum_{{0\le j<s}}\big\langle\bA_{j}\,,\bA_{s-j}\big\rangle_{\C^{2m}}-\,\big\langle\bA_{s/2}\,,\bA_{s/2}\big\rangle_{\C^{2m}}=\;0\:,
\ee
where $\bA_0:=\bA$, and last term is to be ignored when $s$ is odd.\\
For $s=0$, we recover what we have previously observed, namely $\langle\bA,\bA\rangle_{\C^{2m}}=0$, thereby yielding (\ref{condA1}). Using $s=1$ in (\ref{moogli}) gives $\langle\bA,\bA_1\rangle_{\C^{2m}}=0$, so that from (\ref{pronull}) we deduce $\bA_1=\al_1\bA\,$ for some $\al_1\in\C$. Putting now $s=2$ in (\ref{moogli}) gives
\bes
2\,\big\langle\bA\,,\bA_2\big\rangle_{\C^{2m}}\;=\;\big\langle\bA_1\,,\bA_1\big\rangle_{\C^{2m}}\;=\;\al_1^2\,\big\langle\bA\,,\bA\big\rangle_{\C^{2m}}\;=\;0\:,
\ees
so that $\bA_2=\al_2\bA\,$ for some $\al_2\in\C$. Proceeding inductively reveals
\bes
\bA_j\;=\;\al_j\bA\qquad\text{for some $\al_j\in\C$}\qquad\quad\forall\:\:j\in\{0,\ldots,\theta_0-1\}\:. 
\ees
The representation (\ref{baddphi2}) thus becomes
\be\label{supernewphi}
\bp(x)\;=\;\sum_{j=0}^{\theta_0-1}\Re\big(\al_j\bA\;x^{\theta_0+j}\big)\,+\,\vec{\zeta}(x)\:,\qquad\quad\text{with $\:\:\al_0=1$}\:. 
\ee
We will use this formulation to obtain a result describing the behavior of the mean curvature near the singularity and improving (\ref{strausskahn}). To do so, we return to the proof of Proposition \ref{logex}. Setting
\bes
\mathcal{L}(\bH)\;:=\;\text{div}\Big(\nabla\bH\,-\,3\,\pro\nabla\bH\,+\,\star\,(\nabla^\perp\bn\wedge\bH)\Big)\:,
\ees
we saw that
\bes
\mathcal{L}(\bH)\;=\;-\,\bc_0\,\delta_0\:,
\ees 
where $\bc_0$ is the residue defined in (\ref{residudu}). We also proved that
\bes
\bc_0\cdot\bA\;=\;0\:.
\ees
Introducing this fact into (\ref{supernewphi}) then yields
\begin{eqnarray}
\pi_{T(x)}\bc_0&:=&\text{e}^{-2\la(x)}\big(\bc_0\cdot\px\bp(x)\big)\,\px\bp(x)\,+\,\text{e}^{-2\la(x)}\big(\bc_0\cdot\py\bp(x)\big)\,\py\bp(x)\nonumber\\[.5ex]
&=&\text{e}^{-2\la(x)}\big(\bc_0\cdot\px\vec{\zeta}(x)\big)\,\px\bp(x)\,+\,\text{e}^{-2\la(x)}\big(\bc_0\cdot\py\vec{\zeta}(x)\big)\,\py\bp(x)\:.\nonumber
\end{eqnarray}
Hence, calling upon (\ref{zetata}) and (\ref{lambdada}) gives the estimate
\bes
|\nabla^k\pi_T\bc_0|\:\lesssim\:|x|^{\theta_0-k-\epsilon}\qquad\quad\forall\:\:\epsilon>0\:\:,\:\:k\in\{0,\ldots,\theta_0-1\}\:.
\ees
Once combined to (\ref{ocean}), the latter implies
\bes
\bigg|\,\nabla^{\theta_0-2}\,\text{div}\Big[3\,\pi_T\bc_0\,\nabla\log|x|\,+\,\star\,\big(\nabla^\perp\bn\wedge\bc_0\big)\log|x|\Big]\,\bigg|\:\lesssim\:|x|^{-\epsilon}\;\in\;\bigcap_{p<\infty}L^p(B_1(0))\:.
\ees
We saw in the course of the proof of Proposition \ref{logex} that
\bes
\mathcal{L}\Big(4\pi\bH+{\bc_0}\log|x|\Big)\;=\;\text{div}\Big[3\,\pi_T\bc_0\,\nabla\log|x|\,+\,\star\,\big(\nabla^\perp\bn\wedge\bc_0\big)\log|x|\Big]\:.
\ees
Hence
\bes
\mathcal{L}\Big(\nabla^{\theta_0-2}\big(4\pi\bH+{\bc_0}\log|x|\big)\Big)\;\in\;\bigcap_{p<\infty}L^p\:,
\ees
so that (since $\mathcal{L}$ is second-order elliptic and essentially behaves like the Laplacian \cite{Ri1}), 
\be\label{superH}
\bH(x)+\,\dfrac{\bc_0}{4\pi}\log|x|\;\in\;\bigcap_{p<\infty}W^{\theta_0,p}\:.
\ee

Finally, to obtain that $\nabla^{\theta_0+1}\bn\in L^{2,\infty}$, and thus that $\nabla^{\theta_0}\bn\in BMO$, we return to (\ref{prof444}). Since $\text{e}^{\la}\vec{h}_{ij}$ inherits the regularity of $\nabla\bn$, it is not hard to obtain
\bes
\big|\Delta\nabla^{\theta_0-1}\bn\big|\:\:\lesssim\:\:\sum_{k=0}^{\theta_0-1}\,|\nabla^{k+1}\bn|\,|\nabla^{\theta_0-k}\bn|+|\nabla^{k+1}\bp|\,|\nabla^{\theta_0-k}\bH|
\ees
Bringing (\ref{baddphi}), (\ref{zetata}), (\ref{ocean}), and (\ref{superH}) into the latter shows that for all $\epsilon>0$, 
\bes
\big|\Delta\nabla^{\theta_0-1}\bn\big|\:\:\lesssim\:\:|x|^{\theta_0-1-\epsilon}+\,|x|^{-1}+\,\text{terms in}\bigcap_{p<\infty}L^p\:\in\:L^{2,\infty}\:,
\ees
thereby yielding the sought out result.

\subsection{When the residue $\bc_0$ vanishes: point removability}

This last section is devoted to proving Theorem \ref{Th5}. We shall assume that the residue defined in (\ref{residudu}) satisfies $\bc_0=\vec{0}$, and furthermore when $\theta_0\ge2$ that the constant vector appearing in Proposition \ref{Th2}-(ii) is also null: $\bC=\vec{0}$. \\[1.5ex]
%In such cases, we have shown that
%\bes
%\bH\;\in\;\bigcap_{p<\infty}W^{1,p}(B_1(0))\:.
%\ees
When $\bc_0=\vec{0}$, the functions $g$ and $\bG$ vanish identically, so the conservative conformal Willmore system (\ref{sysSR})-(\ref{delphi}) reads
\be\label{sysSRc0}
\left\{\begin{array}{rcl}
-\,\Delta S&=&\nabla(\star\,\bn)\cdot\nabla^\perp\bR\\[1ex]
-\,\Delta\bR&=&\nabla(\star\,\bn)\bul\nabla^\perp\bR\;-\,\nabla(\star\,\bn)\cdot\nabla^\perp S\\[1.5ex]
-\,2\,\Delta\bP&=& \nabla S\cdot\nabla^\perp\bP\,-\,\nabla\bR\bul\nabla^\perp\bP\:.
\end{array}\right.
\ee
When $\theta_0=1$, the immersion $\bp$ is non-degenerate at the origin: its gradient is bounded from above and below. In this case, it was shown in \cite{BR1}\footnote{cf. last paragraph of Section III.2.1.} that the system (\ref{sysSRc0}) yields that $\bp$ is smooth across the unit disk. In the case when $\theta_0\ge2$, we have shown in the previous section that
\bes
\nabla^{\theta_0}\bn\;,\;\nabla^{\theta_0}S\;,\;\nabla^{\theta_0}\bR\:\in\bigcap_{p<\infty}\!L^p\:.
\ees
It immediately follows from the first two equations in (\ref{sysSRc0}) that $S$ and $\bR$ lie in $W^{\theta_0+1,p}$ for all $p<\infty$. We are thus in the position of applying the procedure given in the previous section (since $g$ and $\bG$ no longer obstruct) and deduce that $\bn\in W^{\theta_0+1,p}$ for all $p<\infty$. A bootstrapping argument is then implemented to increase the regularity of all functions involved up to $C^\infty(B_1(0))$. This is in particular the case for the immersion $\bp$, thereby completing the proof of Theorem \ref{Th5}.

\renewcommand{\theequation}{A.\arabic{equation}}
\renewcommand{\theTh}{A.\arabic{Th}}
\renewcommand{\theProp}{A.\arabic{Prop}}
\renewcommand{\theLma}{A.\arabic{Lma}}
\renewcommand{\theCo}{A.\arabic{Co}}
\renewcommand{\theRm}{A.\arabic{Rm}}
\renewcommand{\theequation}{A.\arabic{equation}}
\setcounter{equation}{0} 
\reset
\appendix
\section{Appendix}
\subsection{Notational Conventions}

We append an arrow to all the elements belonging to $\R^m$. To simplify the notation, by $\bp\in X(\di)$ is meant $\bp\in X(\di,\R^m)$ whenever $X$ is a function space. Similarly, we write $\nabla\bp\in X(\di)$ for $\nabla\bP\in \mathbb{R}^2\otimes X(\di,\R^m)$.\\[1.5ex]
Although this custom may seem at first odd, we allow the differential operators classically acting on scalars to act on elements of $\R^m$. Thus, for example, $\nabla\bp$ is the element of $\R^2\otimes\R^m$ that can be written $(\px\bp,\py\bp)$. If $S$ is a scalar and $\bR$ an element of $\R^m$, then we let
\begin{eqnarray*}
\bR\cdot\nabla\bP&:=&\big(\bR\cdot\px\bP\,,\,\bR\cdot\py\bP\big)\:\\[1ex]
\nabla^\perp S\cdot\nabla\bP&:=&\px S\,\py\bp\,-\,\py S\,\px\bp\:\\[1ex]
\nabla^\perp\bR\cdot\nabla\bP&:=&\px\bR\cdot\py\bp\,-\,\py\bR\cdot\px\bp\:\\[1ex]
\nabla^\perp\bR\wedge\nabla\bP&:=&\px\bR\wedge\py\bp\,-\,\py\bR\wedge\px\bp\:.
\end{eqnarray*}
Analogous quantities are defined according to the same logic. \\

Two operations between multivectors are useful. The {\it interior multiplication} $\res$ maps a pair comprising a $q$-vector $\gamma$ and a $p$-vector $\beta$ to a $(q-p)$-vector. It is defined via
\bes
\langle \gamma\res\beta\,,\alpha\rangle\;=\;\langle \gamma\,,\beta\wedge\alpha\rangle\:\qquad\text{for each $(q-p)$-vector $\al$.}
\ees
Let $\al$ be a $k$-vector. The {\it first-order contraction} operation $\bul$ is defined inductively through 
\bes
\al\bul\beta\;=\;\al\res\beta\:\:\qquad\text{when $\beta$ is a 1-vector}\:,
\ees
and
\bes
\al\bul(\beta\wedge\gamma)\;=\;(\al\bul\beta)\wedge\gamma\,+\,(-1)^{pq}\,(\al\bul\gamma)\wedge\beta\:,
\ees
when $\beta$ and $\gamma$ are respectively a $p$-vector and a $q$-vector.

\subsection{Miscellaneous Facts}

\subsubsection{On the Gauss Map}

Let $\bp$ be a conformal immersion of the unit-disk into $\R^m$. 
By definition, for $j\in\{1,2\}$,
\bes
\vec{e}_j\;:=\;\text{e}^{-\la}\,\pj\bP\qquad\text{with}\quad 2\,\text{e}^{2\la}\,=\,|\nabla\bp|^2\:.
\ees
One easily verifies (cf. details in \cite{BR1} Section III.2.2) that\footnote{$\pi_{T}$ denotes projection onto the tangent space spanned by $\{\bAe_1\,,\bAe_2\}$.}
\be\label{atchoum}
\pi_T\nabla\bAe_j\;=\;(\nabla^\perp\la)\,\bAe_{j'}\qquad\quad\text{where}\qquad(\bAe_{1'}\,,\bAe_{2'})\,:=\,(\bAe_{2}\,,-\bAe_{1})\:.
\ee
Moreover\footnote{$\pro$ denotes projection onto the normal space, namely $\pro=\text{id}-\pi_T$.}, 
\be\label{simplet}
\pro\nabla\bAe_j\;\equiv\:\text{e}^{-\la}\pro\nabla\partial_j\bP\:=:\;\text{e}^\la\!\left(\begin{array}{c}\vec{h}_{1j}\\[1ex]\vec{h}_{2j}\end{array}\right)\:.
\ee
With this notation, the mean curvature vector takes the form
\be\label{weingarten}
\bH\;=\;\dfrac{1}{2}\,\big(\vec{h}_{11}+\vec{h}_{22}\big)\:.%\qquad\text{and}\qquad\bH_{0}\;=\;\dfrac{1}{2}\,\big(\vec{h}_{11}-\vec{h}_{22}\big)\,+\,i\,\vec{h}_{12}\:.
\ee
The $(m-2)$-vector $\bn$ satisfies $\,\bn:=\star\,(\bAe_1\wedge\bAe_2)$. Accordingly, using (\ref{atchoum}), there holds
\be\label{timide}
\nabla\bn\;=\;\star\,\Big[\big(\pro\nabla\bAe_1\big)\wedge\bAe_{2}\,+\,\bAe_{1}\wedge\big(\pro\nabla\bAe_2\big)\Big]\:,
\ee
so that
\begin{eqnarray*}
\Delta\bn&=&\star\,\Big[\text{div}\big(\pro\nabla\bAe_1\big)\wedge\bAe_{2}\,+\,\bAe_{1}\wedge\text{div}\big(\pro\nabla\bAe_2\big)\Big]\;+\;2\star\big[\pro\nabla\bAe_1\wedge\pi_n\nabla\bAe_{2}\big]\nonumber\\
&&\quad+\:\star\Big[\pro\nabla\bAe_1\wedge\pi_T\nabla\bAe_{2}\,+\,\pi_T\nabla\bAe_{1}\wedge\pro\nabla\bAe_2\Big]\:.
\end{eqnarray*}
The identities (\ref{atchoum}) yield
\bes
\pi_T\nabla\bAe_k\wedge\pro\nabla\bAe_{l}\;=\;(\nabla^\perp\la)\cdot\big(\bAe_{k'}\wedge\pro\nabla\bAe_l\big)\:,
\ees
and thus
\begin{eqnarray}\label{prof}
\Delta\bn&=&\star\,\Big[\text{div}\big(\pro\nabla\bAe_1\big)\wedge\bAe_{2}\,+\,\bAe_{1}\wedge\text{div}\big(\pro\nabla\bAe_2\big)\Big]\;+\;2\star\big[\pro\nabla\bAe_1\wedge\pi_n\nabla\bAe_{2}\big]\nonumber\\
&&\quad+\:\star(\nabla^\perp\la)\cdot\Big[\bAe_1\wedge\pro\nabla\bAe_{1}\,+\,\bAe_{2}\wedge\pro\nabla\bAe_2\Big]\:.
\end{eqnarray}
Next, using the definition of $\bAe_k$ and again (\ref{atchoum}), we obtain\footnote{implicit summations over repeated indices are understood.}
\begin{eqnarray*}
\text{div}\,\pro\nabla\bAe_k&\equiv&\pro\,\text{div}\,\pro\nabla\bAe_k\;+\;\pi_T\,\text{div}\,\pro\nabla\bAe_k\\[1ex]
&=&\pro\,\text{div}\,\pro\,\nabla\big(\text{e}^{-\la}\partial_k\bP\big)\;+\;\big(\bAe_l\cdot\text{div}\,\pro\nabla\bAe_k\big)\,\bAe_l\\[1ex]
&=&\text{e}^{-\la}\,\pro\,\text{div}\,\pro\nabla\partial_k\bP\;-\;\text{e}^{-\la}\,\pro\big(\nabla\la\cdot\nabla\partial_k\bp\big)\\
&&\quad\;-\;\big(\pro\,\nabla\bAe_l\cdot\pro\nabla\bAe_k\big)\,\bAe_l\\[1ex]
&=&\text{e}^{-\la}\,\pro\,\text{div}\,\pro\nabla\partial_k\bP\;-\;\big(\pro\nabla\bAe_l\cdot\pro\nabla\bAe_k\big)\,\bAe_l\,-\,\nabla\la\cdot\pro\nabla\bAe_k\:.
\end{eqnarray*}
Introducing the latter into (\ref{prof}) gives after a few elementary manipulations,
\begin{eqnarray*}
\Delta\bn&=&\star\,\text{e}^{-\la}\Big[\pro\,\text{div}\big(\pro\nabla\px\bp\big)\wedge\bAe_{2}\,+\,\bAe_{1}\wedge\pro\,\text{div}\big(\pro\nabla\py\bp\big)\Big]\\
&&-\;\;\Big[|\pro\nabla\bAe_1\big|^2+|\pro\nabla\bAe_2\big|^2\Big]\star(\bAe_1\wedge\bAe_2)\;+\;2\star\big[\pro\nabla\bAe_1\wedge\pi_n\nabla\bAe_{2}\big]\nonumber\\
&&+\:\star(\nabla^\perp\la)\cdot\Big[\bAe_1\wedge\pro\big(\nabla\bAe_{1}-\nabla^\perp\bAe_2\big)\,-\,\pro\big(\nabla^\perp\bAe_1+\nabla\bAe_2\big)\wedge\bAe_{2}\Big]\nonumber\:.
\end{eqnarray*}
Owing to (\ref{simplet}) and (\ref{timide}), we find
\begin{eqnarray*}\label{prof2}
\Delta\bn\,+\,|\nabla\bn|^2\,\bn&=&\star\,\text{e}^{-\la}\Big[\pro\,\text{div}\big(\pro\nabla\px\bp\big)\wedge\bAe_{2}\,+\,\bAe_{1}\wedge\pro\,\text{div}\big(\pro\nabla\py\bp\big)\Big]\nonumber\\
&&\hspace{-2.5cm}+\;2\star\text{e}^{-2\la}\big[\pro\nabla\px\bP\wedge\pro\nabla\py\bp\big]\:+\:2\star\text{e}^{\la}\bH\wedge\big[\py\la\;\bAe_1\,-\,\px\la\;\bAe_2\big]\,.
\end{eqnarray*}
Equivalently, 
\begin{eqnarray}\label{prof3}
\Delta\bn\,+\,|\nabla\bn|^2\,\bn&=&\star\,\bAe_1\wedge\pro\Big[\text{e}^{-\la}\,\text{div}\big(\pro\nabla\py\bp\big)\,-\,2\,\text{e}^{\la}\bH\,\py\la\Big]\nonumber\\
&&-\;\star\bAe_2\wedge\pro\Big[\text{e}^{-\la}\,\text{div}\big(\pro\nabla\px\bp\big)\,-\,2\,\text{e}^{\la}\bH\,\px\la\Big]\nonumber\\[.5ex]
&&\hspace{.5cm}+\;2\star\text{e}^{-2\la}\big[\pro\nabla\px\bP\wedge\pro\nabla\py\bp\big]\:.
\end{eqnarray}
%Using the fact that $\,\Delta\bp=2\,\text{e}^{2\la}\bH$, we note that
%\be\label{bailey}
%\text{e}^{-\la}\,\text{div}\big(\pro\nabla\pj\bp\big)\,-\,2\,\text{e}^{\la}\bH\,\pj\la\;=\;\text{e}^{\la}\pj\bH\,-\,\text{e}^{-\la}\text{div}\,\pi_T\nabla\pj\bp\:.
%\ee
Moreover, (\ref{atchoum}) gives $\,\pi_{T}\nabla\pj\bp=\nabla(\text{e}^{\la})\bAe_j+\nabla^\perp(\text{e}^{\la})\bAe_{j'}$. Hence, calling upon (\ref{simplet}) implies
\bes
\pro\,\text{div}\,\pi_T\nabla\pj\bp\;=\;\nabla(\text{e}^{\la})\cdot\pro\nabla\bAe_j\,+\,\nabla^\perp(\text{e}^{\la})\cdot\pro\nabla\bAe_{j'}\;=\;\bH\,\pj\text{e}^{2\la}\:,
\ees
and thus, as $\Delta\bp=2\text{e}^{2\la}\bH$, 
\begin{eqnarray*}
\pro\,\text{div}\,\pro\nabla\pj\bp&\equiv&\pro\,\pj\Delta\bp\,-\,\pro\,\text{div}\,\pi_T\nabla\pj\bp\\[1ex]
&=&2\,\pro\,\pj\big(\text{e}^{2\la}\bH\big)\,-\,\bH\,\pj\text{e}^{2\la}\:.
\end{eqnarray*}
The interested reader will note that this equation is equivalent to the Codazzi-Mainardi identities. Substituted into (\ref{prof3}), the latter gives
\begin{eqnarray}\label{prof4}
\Delta\bn\,+\,|\nabla\bn|^2\,\bn&=&2\star\text{e}^{\la}\Big[\bAe_1\wedge\pro\,\py\bH\,-\,\bAe_2\wedge\pro\,\px\bH\Big]\nonumber\\[.5ex]
&&\hspace{1.5cm}+\;2\star\text{e}^{-2\la}\,\big[\pro\nabla\px\bP\wedge\pro\nabla\py\bp\big]\nonumber\\[1ex]
&=&2\star\text{e}^{\la}\Big[\bAe_1\wedge\pro\,\py\bH\,-\,\bAe_2\wedge\pro\,\px\bH\Big]\nonumber\\[.5ex]
&&\hspace{1.5cm}-\;2\star\text{e}^{2\la}\,\vec{h}_{12}\wedge(\vec{h}_{11}-\vec{h}_{22})\:.
\end{eqnarray}

\subsubsection{Conservative Conformal Willmore System}

We establish in this section a few general identities. As before, we let $\bp$ be a (smooth) conformal immersion of the unit-disk into $\R^m$, and set $\bej:=\text{e}^{-\la}\pj\bp$, where $\la$ is the conformal parameter. Since $\bp$ is conformal, $\{\bex,\bey\}$ forms an orthonormal basis of the tangent space. As $\bn=\star(\bex\wedge\bey)$, if $\bV$ is a 1-vector, we find
\bes
(\star\,\bn)\cdot(\bV\wedge\pj\bp)\;=\;\text{e}^{-\la}(\bex\wedge\bey)\cdot(\bV\wedge\bej)\;=\;-\,\text{e}^{-\la}\,\bAe_{j'}\cdot\bV\;=\;-\,\partial_{x_{j'}}\bp\cdot\bV\:,
\ees
where
\bes
(\bAe_{1'}\,,\bAe_{2'})\,:=\,(\bAe_{2}\,,-\,\bAe_{1})\:.
\ees
%and similarly,
%\bes
%(\star\,\bn)\cdot(\py\bp\wedge\bV)\;=\;\text{e}^{-\la}(\bex\wedge\bey)\cdot(\bey\wedge\bV)\;=\;-\,\text{e}^{-\la}\,\bex\cdot\bV\;=\;-\,\px\bp\cdot\bV\:.
%\ees
Whence,
\be\label{idd1}
\left\{\begin{array}{lcr}
(\star\,\bn)\cdot(\bV\wedge\nabla\bp)&=&\bV\cdot\nabla^\perp\bp\\[1ex]
(\star\,\bn)\cdot(\bV\wedge\nabla^\perp\bp)&=&-\,\bV\cdot\nabla\bp\:.
\end{array}\right.
\ee

\medskip
We choose next an orthonormal basis $\{\bn_\al\}^{m-2}_{\al=1}$ of the normal space such that $\,\{\bex,\bey,\bn_1,\ldots,\bn_{m-2}\}$ is a positive oriented orthonormal basis of $\R^m$. \\
%The definition of the Hodge star operator in Euclidean space shows that $\bn$ is the $(m-2)$-vector $\,\bn=\bigwedge_{1\le\al\le m-2}\bn_\al$. \\
Recalling the definition of the interior multiplication operator $\res$ given in Section A.1, it is not hard to obtain
\bes
(\star\,\bn)\res\bej\;=\;(\bex\wedge\bey)\res\bej\;=\;\delta_{j2}\,\bex\,-\,\delta_{j1}\,\bey\:,
\ees
and
\bes
(\star\,\bn)\res\bn_\al\;=\;0\:.
\ees
Hence,
\begin{eqnarray*}
(\star\,\bn)\bul(\bej\wedge\bn_\al)&\equiv&\big((\star\,\bn)\res\bej\big)\wedge\bn_\al\,+\,\big((\star\,\bn)\res\bn_\al\big)\wedge\bej\\[.5ex]
&=&\delta_{j2}\,\bex\wedge\bn_\al\,-\,\delta_{j1}\,\bey\wedge\bn_\al\:.
\end{eqnarray*}
Moreover, there holds trivially
\bes
(\star\,\bn)\bul(\bej\wedge\bek)\;=\;\pm\,(\star\,\bn)\bul(\star\,\bn)\;=\;0\:.
\ees
From this one easily deduces for every 1-vector $\bV$, one has
\be\label{idd20}
\left\{\begin{array}{lcl}
(\star\,\bn)\bul\big(\bV\wedge\nabla\bp\big)&=&\pro\bV\wedge\nabla^\perp\bp\\[1.75ex]
(\star\,\bn)\bul\big(\bV\wedge\nabla^\perp\bp\big)&=&-\,\pro\bV\wedge\nabla\bp\:.
\end{array}\right.
\ee

There holds furthermore
\bes
%A(\bA,j)\;:=\;
(\bV\wedge\bej)\bul\bei\;=\;(\bei\res\bV)\wedge\bej\,+\,\bV\wedge(\bei\res\bej)\;=\;(\bei\cdot\bV)\,\bej\,+\,\delta_{ij}\,\bV\:.
\ees
From this, and $\,\bei:=\text{e}^{-\la}\partial_{x^i}\bp\,$, it follows that whenever $
\bV=V^{i}\bei+V^{\al}\bna$ then
\be\label{idd3}
\left\{\begin{array}{lll}
\big(\bV\wedge\nabla^\perp\bp\big)\bul\nabla^\perp\bp&=&\text{e}^{2\la}\,\big(3\,V^i\,\bei\,+\,2\,V^\al\,\bna\big)\\[1.5ex]
\big(\bV\wedge\nabla\bp\big)\bul\nabla^\perp\bp&=&\text{e}^{2\la}\,\big(V^{2}\,\bex\,-\,V^{1}\,\bey\big)\;\;\equiv\;\;\big(\bV\cdot\nabla\bp\big)\cdot\nabla^\perp\bp\:.
\end{array}\right.
\ee

\bigskip
We are now sufficiently geared to prove

\begin{Lma}\label{identities}
Let $\bp$ be a smooth conformal immersion of the unit-disk into $\R^m$ with corresponding mean curvature vector $\bH$, and let $\bL$ be a 1-vector. We define $A\in \R^2\otimes\bigwedge^0(\R^m)$ and $\bB\in\R^2\otimes\bigwedge^2(\R^m)$ via
\bes
\left\{\begin{array}{rclll}
A&=&\bL\cdot\nabla\bP&&\\[1.5ex]
\bB&=&\bL\wedge\nabla\bP\,+\,2\,\bH\wedge\nabla^\perp\bP\:.
\end{array}\right.
\ees
Then the following identities hold:
\be\label{hyperid}
\left\{\begin{array}{rclll}
A&=&-\,(\star\,\bn)\cdot\bB^\perp&&\\[1.5ex]
\bB&=&-\,(\star\,\bn)\bul\bB^\perp\,+\,(\star\,\bn)\,A^\perp\:,
\end{array}\right.
\ee
where $\,\star\,\bn:=(\px\bp\wedge\py\bp)/|\px\bp\wedge\py\bp|\,$.\\[1ex]
Moreover, we have
\be\label{hyperdel}
-\,2\,\Delta\bp\;=\;A\cdot\nabla^\perp\bp\,-\,\bB\bul\nabla^\perp\bp\:.
\ee
\end{Lma}
$\textbf{Proof.}$
The identities (\ref{idd1}) give immediately (recall that $\bH$ is a normal vector, so that $\,\bH\cdot\nabla^\perp\bp=0$) the required
\bes
(\star\,\bn)\cdot\bB^\perp\;=\;-\,\bL\cdot\nabla\bp\,+\,2\,\bH\cdot\nabla^\perp\bp\;=\;-\,\bL\cdot\nabla\bp\;=\;-\,A\:.
\ees
Analogously, the identities (\ref{idd20}) give (again, $\bH$ is normal, so $\pro\bH=\bH$),
\begin{eqnarray*}
(\star\,\bn)\bul\bB^\perp&=&-\,\pro\bL\wedge\nabla\bp\,-\,2\,\bH\wedge\nabla^\perp\bp\;\;=\;\;-\,\bB\,+\,\pi_T\bL\wedge\nabla\bP\nonumber\\[1ex]
&=&-\,\bB\,+\,\text{e}^{\la}\big((\bL\cdot\bex)\,\bex\,+\,(\bL\cdot\bey)\,\bey\big)\wedge\left(\begin{array}{r}\bex\\\bey\end{array}\right)\nonumber\\[1ex]
&=&-\,\bB\,+\,\text{e}^{\la}\left(\begin{array}{r}-\,\bL\cdot\bey\\\bL\cdot\bex\end{array}\right)\bex\wedge\bey\nonumber\\[1ex]
&=&-\,\bB\,+\,(\bL\cdot\nabla^\perp\bp)\,(\star\,\bn)\;\;=\;\;-\,\bB\,+\,(\star\,\bn)\,A^\perp\:,
\end{eqnarray*}
which is the second equality in (\ref{hyperid}). \\
In order to prove (\ref{hyperdel}), we will use (\ref{idd3}). Namely, since $\,\bH=H^\al\bna\,$, we find
\bes
\bB\,\bul\nabla^\perp\bp\;=\;\big(\bL\cdot\nabla\bp\big)\cdot\nabla^\perp\bp\,+\,4\,\text{e}^{2\la}\bH\;=\;A\cdot\nabla^\perp\bp\,+\,4\,\text{e}^{2\la}\bH\:.
\ees
Hence,
\bes
\bB\,\bul\nabla^\perp\bp\,-\,A\cdot\nabla^\perp\bp\;=\;4\,\text{e}^{2\la}\bH\:.
\ees
Finally, there remains to recall that $\,\Delta\bp=2\,\text{e}^{2\la}\bH\,$
to reach the desired identity.\\[-1.5ex]

$\hfill\blacksquare$ \\

\noindent
We choose now
\bes
A\;=\;\nabla S\,-\nabla^\perp g\qquad\text{and}\qquad \bB\;=\;\nabla\bR\,-\nabla^\perp\bG\:,
\ees
where $S$ and $g$ are scalars and $\bR$ and $\bG$ are 2-vectors. Then Lemma \ref{identities} yields
\bes
\left\{\begin{array}{rclll}
\nabla S&=&-\,(\star\,\bn)\cdot\big(\nabla^\perp\bR\,+\nabla\bG\big)\,+\,\nabla^\perp g&&\\[1.5ex]
\nabla\bR&=&-\,(\star\,\bn)\bul\big(\nabla^\perp\bR\,+\nabla\bG\big)\,+\,(\star\,\bn)\,\big(\nabla^\perp S\,+\nabla g\big)\,+\,\nabla^\perp\bG\:,
\end{array}\right.
\ees
thereby giving
\be\label{sysSR0}
\left\{\begin{array}{rclll}
-\,\Delta S&=&\nabla(\star\,\bn)\cdot\nabla^\perp\bR\,+\,\text{div}\big((\star\,\bn)\cdot\nabla\bG\big)&&\\[1.5ex]
-\,\Delta\bR&=&\nabla(\star\,\bn)\bul\nabla^\perp\bR\;-\;\nabla(\star\,\bn)\cdot\nabla^\perp S\\[.75ex]
&&\hspace{2.15cm}\,+\:\,\text{div}\big((\star\,\bn)\bul\nabla\bG\,-\,\star\,\bn\,\nabla g\big)\:.
\end{array}\right.
\ee
Furthermore, there holds,
\be\label{delphi0}
-\,2\,\Delta\bp\;=\;(\nabla S-\nabla^\perp g)\cdot\nabla^\perp\bp\,-\,(\nabla\bR-\nabla^\perp\bG)\bul\nabla^\perp\bp\:.
\ee

\medskip

\subsection{Nonlinear and weighted elliptic results}

\begin{Prop}\label{morreydecay}
Let $u\in W^{1,2}(B_1(0))\cap C^1(B_1(0)\setminus\{0\})$ satisfy the equation
\be\label{equ}
-\,\Delta u\:=\:\nabla b\cdot\nabla^\perp u\,+\text{div}\,(b\,\nabla f)\qquad\quad\text{on}\:\:\:\:B_1(0)\:,
\ee
where $\,f\in W_0^{2,(2,\infty)}(B_1(0))$, and moreover
\be\label{hypn}
b\,\in\,W^{1,2}\cap L^{\infty}(B_1(0))\qquad\text{with}\qquad\Vert\nabla b\Vert_{L^2(B_1(0))}\;<\;\eps_0\:,
\ee
for some $\eps_0$ chosen to be ``small enough". Then
\bes
\nabla u\,\in\,L^p(B_{1/4}(0))\qquad\text{for some $\,p>2$}\:.
\ees
\end{Prop}
$\textbf{Proof.}$ Before delving into the proof of the statement, one important remark is in order. Let $D$ be any disk included (properly or not) in $B_1(0)$. From the very definition of the space $L^{2,\infty}$ (cf. \cite{Ta}), there holds
\be\label{estimf1}
\Vert\Delta f\Vert_{L^1(D)}\:\le\:|D|^{\frac{1}{2}}\Vert\Delta f\Vert_{L^{2,\infty}(D)}\:\lesssim\:|D|^{\frac{1}{2}}\Vert\nabla^2f\Vert_{L^{2,\infty}(D)}\:.
\ee
Moreover, an embedding result of Luc Tartar \cite{Ta} states that $\nabla f$ has bounded mean oscillations. Whence in particular
\be\label{estimf2}
\Vert\nabla f\Vert_{L^2(D)}\;\lesssim\;|D|^{\frac{1}{2}-\epsilon}\qquad\forall\:\:\epsilon>0\:.
\ee
These inequalities shall come helpful in the sequel. \\[1ex]
We now return to the proof of the proposition. Let us fix some point $x_0\in B_{1/2}(0)$ and some radius $\sigma\in(0\,,\frac{1}{2})$, and we let $k\in(0\,,1)$. Note that $B_{k\sigma}(x_0)$ is properly contained in $B_1(0)$. To reach the desired result, we decompose the solution to (\ref{equ}) as the sum $\,u=u_0+u_1\,$, where 
\bes
\left\{\begin{array}{rcl}
-\,\Delta u_0&=&\text{div}\,(b\,\nabla f)\\[1ex]
u_0&=&u\end{array}\right.\quad,\quad
\left.\begin{array}{rclcl}
-\,\Delta u_1&=&\nabla b\cdot\nabla^\perp u_{}&\quad&\text{in\:\:\:} B_\sigma(x_0)\\[1ex]
u_{1}&=&0&\quad&\text{on\:\:\:}\partial B_\sigma(x_0)\:.
\end{array}\right.
\ees
Accounting for the hypotheses (\ref{hypn}) and (\ref{estimf2}) into standard elliptic estimates (cf. Proposition 4 in \cite{Al}) yields
\begin{eqnarray}\label{b1}
\Vert\nabla u_0\Vert_{L^2(B_{k\sigma}(x_0))}&\lesssim&\Vert b\,\nabla f\Vert_{L^2(B_{k\sigma}(x_0))}\,+\,k\,\Vert\nabla u\Vert_{L^2(B_\sigma(x_0))}\nonumber\\[1ex]
&\lesssim&(k\sigma)^{1-\epsilon}\,+\,k\,\Vert\nabla u\Vert_{L^2(B_\sigma(x_0))}\:,
\end{eqnarray}
up to some unimportant multiplicative constants. 
On the other hand, applying Wente's inequality (see \cite{He} Theorem 3.4.1) gives
\begin{eqnarray}\label{b2}
\Vert\nabla u_1\Vert_{L^2(B_{k\sigma}(x_0))}&\le&\Vert\nabla u_1\Vert_{L^2(B_{\sigma}(x_0))}\nonumber\\[1ex]
&\lesssim&\Vert\nabla b\Vert_{L^2(B_{\sigma}(x_0))}\,\Vert\nabla u_{}\Vert_{L^2(B_{\sigma}(x_0))}\nonumber\\[1ex]
&\leq&\eps_0\,\Vert\nabla u_{}\Vert_{L^2(B_\sigma(x_0))}\:,
\end{eqnarray}
again up to some multiplicative constant without bearing on the sequel. Hence, combining (\ref{b1}) and (\ref{b2}), we obtain the estimate
\begin{eqnarray*}\label{b3}
\Vert\nabla u\Vert_{L^2(B_{k\sigma}(x_0))}&\leq&\Vert\nabla u_0\Vert_{L^2(B_{k\sigma}(x_0))}\;+\;\Vert\nabla u_1\Vert_{L^2(B_{k\sigma}(x_0))}\nonumber\\[1ex]
&\lesssim&(k+\eps_0)\,\Vert\nabla u\Vert_{L^2(B_\sigma(x_0))}\,+\,(k\sigma)^{1-\epsilon}\:.
\end{eqnarray*}
Because $\eps_0$ and $\epsilon$ are small adjustable parameters, we may always choose $k$ so as to arrange for $(k+\eps_0)$ to be less than 1. A standard ``controlled-growth" argument (see e.g. Lemma 3.5.11 in \cite{He}) enables us to conclude that there exists some $\beta\in(0\,,1)$ for which
%\bes
%\Vert\nabla u\Vert_{L^2(B_{\sigma}(x_0))}\:\leq\:C_0\,\sigma^\beta\:,
%C(\epsilon)\,k\,\Vert\nabla u\Vert_{L^2(B_{\sigma}(x_0))}\,+\,C_1\,(k\sigma)^{1-\epsilon}\:.
%\ees
%for some constant $C_0$. We have thus reached the growth estimate
\be\label{decu}
\Vert\nabla u\Vert_{L^2(B_{\sigma}(x))}\:\leq\:C_0\,\sigma^\beta\:,\qquad\forall\:\:\sigma\in\Big(0\,,\frac{1}{2}\Big)\:\:,\:\:\:x\in B_{1/2}(0)\:,
\ee
and for some constant $C_0$. \\
With the help of the Poincar\'e inequality, this estimate may be used to show that $u$ is locally H\"older continuous. We are however interested in another implication of (\ref{decu}). Consider the maximal function
\be\label{maxfun}
M_{2-\beta}\,g(x)\;:=\;\sup_{\sigma>0}\;\sigma^{-\beta}\!\int_{B_\sigma(x)}|g(y)|\,dy\:.
\ee
We recast the equation (\ref{equ}) in the form
\bes
-\,\Delta u\;=\;b\,\Delta f\;+\;\nabla b\cdot\big(\nabla^\perp u+\nabla f\big)\:.
\ees
Calling upon (\ref{hypn})-(\ref{estimf2}) and upon the estimate (\ref{decu}), we derive that for $x\in B_{1/2}(0)$, there holds
\begin{eqnarray}
&&M_{2-\beta}\big(\chi_{B_{1/2}(0)}\Delta u\big)(x)\:\:\leq\:\:\Vert b\Vert_{L^\infty(B_{1}(0))}\,\sup_{0<\sigma<\frac{1}{2}}\,\sigma^{-\beta}\,\Vert\Delta f\Vert_{L^1(B_{\sigma}(x))}\nonumber\\
&&\hspace{1cm}+\:\:\Vert\nabla b\Vert_{L^2(B_{1}(0))}\,\sup_{0<\sigma<\frac{1}{2}}\sigma^{-\beta}\Big(\Vert\nabla u\Vert_{L^2(B_{\sigma}(x))}\,+\,\Vert\nabla f\Vert_{L^2(B_{\sigma}(x))}\Big)  \nonumber\\[1ex]
&&\hspace{1cm}\lesssim\:\:\sup_{0<\sigma<\frac{1}{2}}\sigma^{-\beta+1}\;+\;\eps_0\sup_{0<\sigma<\frac{1}{2}}\big(\sigma^{-\beta+\beta}+\sigma^{-\beta+1-\epsilon}\big)\:\:<\:\:\infty\:,\qquad\qquad
\end{eqnarray}
for all $\,0<\epsilon\le1-\beta$. 
Moreover, it is clear that $\Delta u$ is integrable on $B_{1/2}(0)$. We may thus use Proposition 3.2 from \cite{Ad}\footnote{namely, $\,\big\Vert|x|^{-1}\!*g\big\Vert^{r}_{L^{r,\infty}}\lesssim\,\Vert M_{2-\beta}\,g\Vert_{L^\infty}^{1-\frac{1}{r}\,}\Vert g\Vert_{L^1}^{\frac{1}{r}}\,$ for $r=\frac{2-\beta}{1-\beta}\,$ and $\beta\in(0,1)$.} %We use in particular $\,g=\chi_{B_{1/2}}\Delta u$.} 
to deduce that
\bes
\dfrac{1}{|x|}*\chi_{B_{1/2}(0)}\Delta u\;\in\;L^{r,\infty}(B_{1/2}(0))\qquad\text{with}\quad r\,:=\,\dfrac{2-\beta}{1-\beta}\,>\,2\:.
\ees
A classical estimate about Riesz kernels states there holds in general
\bes
|\nabla u|(y)\:\lesssim\:\dfrac{1}{|x|}*\chi_{B_{1/2}(0)}\Delta u\,+\,C\:,\qquad\forall\:\:y\in B_{1/4}(0))\:,
\ees
where $C$ is a constant depending on the $C^1$-norm of $u$ on $\partial B_{1/2}(0)$, hence finite by hypothesis. It follows in particular that 
\bes
\nabla u\,\in\,L^p(B_{1/4}(0))\quad\:\:\:\text{for all}\:\:\:\quad p<r\:,
\ees
as announced.\\[-3ex]

$\hfill\blacksquare$ \\[-1.5ex]

\begin{Prop}\label{CZpondere}
Let $u\in C^2(B_1(0)\setminus\{0\})$ solve
\be\label{eqw}
\Delta u(x)\;=\;\mu(x)f(x)\qquad\text{in\:\:\:} B_1(0)\:,
\ee
where $f\in L^p(B_1(0))$ for some $p>2$. The weight $\mu$ satisfies
\be\label{hypw}
|\mu(x)|\;\simeq\;|x|^{a}\qquad\quad\text{for some}\:\: a\in\mathbb{N}\:.
\ee
%Suppose that
%\be\label{hypu}
%|\mu(x)|^{-1}\,\nabla u(x)\;\in\;L^\infty(B_1(0))\:.
%\ee
%We also assume that $f$ is identically vanishing outside of $B_1(0)$.\\
%If
%\be\label{hypw}
%|x|^{-a}\,\nabla^2 w(x)\;\in\; L^2(B_1(0))\:,
%\ee
Then 
\begin{itemize}
\item[(i)] there holds\footnote{$\overline{x}$ is the complex conjugate of $x$. Namely, we parametrize $B_1(0)$ by $x=x_1+i\,x_2$, and then $\overline{x}=x_1-i\,x_2$. With this notation, $\nabla u$ on the left-hand side of (\ref{stmt}) is understood as $\partial_{x_1}u+i\,\partial_{x_2}u$.}
\be\label{stmt}
\nabla u(x)\;=\;P(\overline{x})\,+\,|\mu(x)|\,T(x)\:,
\ee
where $P(\overline{x})$ is a complex-valued polynomial of degree at most $a$, and near the origin $\,T(x)=\text{O}\big(|x|^{1-\frac{2}{p}-\epsilon}\big)$ for every $\epsilon>0$. \\[-1.5ex]
\item[(ii)] furthermore, if $\,\mu\in C^1(B_1(0)\setminus\{0\})$, if $a\ne0$, and if
\be\label{hypw2}
|x|^{1-a}\,\nabla\mu(x)\,\in\,L^{\infty}(B_1(0))\:,
\ee
there holds
\be\label{stmt2}
\nabla^2 u(x)\;=\;\nabla P(\overline{x})\,+\,|\mu(x)|\,Q(x)\:,
\ee
where $P$ is as in (i), and 
\bes
Q\;\in\;L^{p-\epsilon}(B_1(0),\C^2)\qquad\quad\forall\:\:\epsilon>0\:.
\ees
As a $(2\times2)$ real-valued matrix, $Q$ satisfies in addition
\bes
\text{Tr}\;Q\;\in\;L^p(B_1(0))\:.
\ees

Naturally, if $a=0$, the standard Calderon-Zygmund theorem yields that $u\in W^{2,p}(B_1(0))$. The hypothesis (\ref{hypw2}) becomes unnecessary, and (\ref{stmt2}) holds with $P$ being constant and $\epsilon=0$. 
%\item[(iii)] under the conditions (\ref{hypw}) and (\ref{hypw2}), in the special case when $p=\infty$, the $(2\times2)$ real-valued matrix $Q(x)$ from (\ref{stmt2}) satisfies $\text{Tr}\;Q\in L^\infty(B_1(0))$. 
\end{itemize}
%\bes
%T(x)\;=\;\mathbf{O}\big(|x|^{1-\frac{2}{p}-\epsilon}\big)\:\qquad\quad\forall\:\:\:\:0<\epsilon\,\le\,1-\dfrac{2}{p}\:.
%\ees
%Furthermore, if
%\be\label{hypu2}
%|x|^{-a}\,\pro\nabla^2\bu(x)\;\in\;L^2(B_1(0))\:,
%\ee
%then in fact
%\be\label{hypu2}
%|x|^{-a}\,\pro\nabla^2\bu(x)\;\in\;L^{p-\epsilon}(B_1(0))\:,
%\ee
%for all positive $\epsilon$ small enough. 
\end{Prop}
$\textbf{Proof.}$ Using Green's formula for the Laplacian, an exact expression for the solution $u$ may be found and used to obtain
\begin{eqnarray}\label{green}
\nabla u(x)&=&\dfrac{1}{2\pi}\,\int_{\partial B_1(0)}\bigg[\dfrac{x-y}{|x-y|^2}\,\partial_{\vec{\nu}}\,u(y)\,-\,u(y)\,\partial_{\vec{\nu}}\,\dfrac{x-y}{|x-y|^2}\bigg]\,d\sigma(y)\nonumber\\[.5ex]
&&\hspace{1cm}-\:\:\dfrac{1}{2\pi}\int_{B_1(0)}\dfrac{x-y}{|x-y|^2}\,\mu(y)f(y)\,dy\nonumber\\[1ex]
&=:&J_0(x)\,+\,J_1(x)\:,\hspace{3cm}\forall\:\:x\in B_1(0)\:,
\end{eqnarray}
where $\vec{\nu}$ is the outer normal unit-vector to the boundary of $B_1(0)$. 
Without loss of generality, and to avoid notational clutter, because $u$ is twice differentiable away from the origin, we shall henceforth assume that $|x|<1/2$.\\
We will estimate separately $J_0$ and $J_1$, and open the discussion by noting that when $|y|>|x|$, we have the expansion
\bes
\dfrac{x-y}{|x-y|^2}\;=\;-\,\sum_{m\ge0}P^m(x,y)\qquad\text{with}\qquad P^m(x,y)\,:=\,\overline{x}^{\,m}\,\overline{y}^{\,-(m+1)}\:.
\ees
Hence, we deduce the identity
\begin{eqnarray}\label{J0}
J_0(x)&=&-\,\dfrac{1}{2\pi}\,\sum_{m\ge0}\,\int_{\partial B_1(0)}\big[P^m(x,y)\,\partial_{\vec{\nu}}\,u(y)\,-\,u(y)\partial_{\vec{\nu}}\,P^m(x,y)\big]\,dS(y)
\nonumber\\[1ex]
&=&-\,\dfrac{1}{2\pi}\,\sum_{m\ge 0}\,\overline{x}^{\,m}\int_{0}^{2\pi}\big[(m+1)\,u(\text{e}^{i\varphi})\,-\,(\partial_{\vec{\nu}}\,u)(\text{e}^{i\varphi})\big]\,\text{e}^{i(m+1)\varphi}d\varphi\;\nonumber\\[1ex]
&=&\sum_{m\ge0}C_m\,\overline{x}^{\,m}\:,
\end{eqnarray}
where $C_m$ are (complex-valued) constants depending only on the $C^1$-norm of $u$ along $\partial B_1(0)$. 
As $u$ is continuously differentiable on the boundary of the unit disk by hypothesis, and $|x|<1$, it is clear that $|J_0(x)|$ is bounded above by some constant $C$ for all $x\in B_1(0)$. Since $\,|C_m|$ grows sublinearly in $m$, we can surely find two constants  $\gamma$ and $\delta$ such that
\bes
|C_m|\;<\;\gamma\,\delta^m\qquad\quad\forall\:\:\:m\,\ge\,0\:.
\ees
Hence, when $\,|x|\le R<\delta^{-1}$, there holds
\bes
\bigg|\sum_{m\ge a+1}C_m\,\overline{x}^{\,m}\bigg|\;\le\;\gamma\,\delta^{a+1}\,|x|^{a+1}\,\sum_{m\ge0}(\delta R)^m\;\lesssim\;|x|^{a+1}\:.
\ees
And because $J_0$ is bounded, when $\,R<|x|<1$, we find some large enough constant $K=K(C,a,\gamma,\delta)$ such that
\begin{eqnarray*}
\bigg|\sum_{m\ge a+1}C_m\,\overline{x}^{\,m}\bigg|&\le&|J_0(x)|\,+\!\sum_{0\le m\le a}C_m\,|x|^m\;\;\leq\;\;C\,+\,(a+1)\,\gamma\,\delta^a\nonumber\\[1ex]
&\leq&K\,\delta^{a+1}\;\leq\; K\big(R^{-1}\delta\big)^{\!a+1}\,|x|^{a+1}\;\lesssim\;|x|^{a+1}\:.
\end{eqnarray*}
As by hypothesis $\,|\mu(x)|\simeq|x|^a\,$, we may now return to (\ref{J0}) and write
\be\label{estJ0}
J_0(x)\:=\:P_0(\overline{x})\,+\,|\mu(x)|\,T_0(x)\:,
\ee
where $P_0$ is a polynomial of degree at most $a$, and the remainder $T_0$
%\bes
%T_0(x)\:=\:C\,\sum_{m\ge 0}(m+a)\,|x|^{m}\,\text{e}^{-i(m+a+1)\theta}
%\ees
is controlled by some constant depending on the $C^1$-norm of $u$ on $\partial B_1(0)$. Moreover, $T_0(x)=\text{O}(|x|)$ near the origin. \\[1ex]
We next estimate the integral $J_1$. To do so, we proceed as above and write
\be\label{J1}
J_1(x)\;=\;I_1(x)\,+\sum_{m=a+1}^\infty\!\!I_2^m(x)\,-\sum_{m=0}^{a}I_1^m(x)\,+\sum_{m=0}^{a}I_1^m(x)+I_2^m(x)\:,
\ee
where we have put
\bes
I_1(x)\,:=\,\dfrac{1}{2\pi}\int_{B_1(0)\cap B_{2|x|}(0)}\dfrac{x-y}{|x-y|^2}\,\mu(y)f(y)\,dy\:,
\ees
\bes
I^m_1(x)\,:=\,\dfrac{1}{2\pi}\int_{B_1(0)\cap B_{2|x|}(0)}P^m(x,y)\,\mu(y)f(y)\,dy\:,
\ees
\bes
I^m_2(x)\,:=\,\dfrac{1}{2\pi}\int_{B_1(0)\setminus B_{2|x|}(0)}P^m(x,y)\,\mu(y)f(y)\,dy\:.
\ees\\
We first observe that the last sum in (\ref{J1}) may be written
\begin{eqnarray*}\label{poly}
P_1(x)&:=&\sum_{0\le m\le a}I_1^m(x)+I_2^m(x)\:\:=\:\:\sum_{0\le m\le a}\int_{B_1(0)}P^m(x,y)\,\mu(y)f(y)\,dy\nonumber\\[0ex]
&=&\sum_{0\le m\le a}A_m\,\overline{x}^{\,m}\:,
\end{eqnarray*}
where
\bes
A_m\;:=\;-\int_{B_1(0)}\overline{y}^{\,-(m+1)}\mu(y)f(y)\,dy\:.
\ees\\
From the fact that $f\in L^p(B_1(0))$ for $p>2$, and the hypothesis $|\mu(y)|\simeq|y|^a$, it follows easily that $|A_m|<\infty$ for $m\le a$, and thus that $P_1$ is a polynomial of degree at most $a$. \\
We have next to handle the other summands appearing in (\ref{J1}), beginning with $I_1$. We find
\begin{eqnarray}\label{estI1}
|I_1(x)|&\lesssim&|\mu(x)|\,\int_{B_{2|x|}(0)}\dfrac{|f(y)|}{|x-y|}\,dy\:\:\lesssim\:\:|\mu(x)|\,\int_{B_{3|x|}(x)}\dfrac{|f(y)|}{|x-y|}\,dy\nonumber\\[1ex]
&\lesssim&|\mu(x)|\,|x|\,M_0f(x)\;\;\lesssim\;\;|x|^{1-\frac{2}{p}}|\mu(x)|\:,
\end{eqnarray}
where we have used the fact that $B_{2|x|}(0)\subset B_{3|x|}(x)$, and a classical estimate bounding convolution with the Riesz kernel by the maximal function\footnote{cf. (\ref{maxfun}) for the definition of $M_0f$.} (cf. Proposition 2.8.2 in \cite{Zi}). We have also used the simple estimate $M_0f(x)\lesssim |x|^{-\frac{2}{p}}\Vert f\Vert_{L^p}$. \\[1ex]
Next, let $q\in[1,2)$ be the conjugate exponent of $p$. We immediately deduce for $0\le m\le a$ that
\begin{eqnarray}\label{estI1a}
|I^m_1(x)|&\lesssim&|x|^{m}\int_{B_{2|x|}(0)}|y|^{-1-m+a}|f(y)|\,dy\nonumber\\[1ex]
&\lesssim&|x|^{a}\,\big\Vert |y|^{-1}\big\Vert_{L^q(B_{2|x|}(0))}\,\Vert f\Vert_{L^p(B_{1}(0))}\:\:\lesssim\:\:|x|^{1-\frac{2}{p}}|\mu(x)|\:.
\end{eqnarray}
We next estimate $I^m_2$. As $m\ge a+1$, we note that for any $\epsilon>0$, there holds
\bes
a+1-m-\epsilon-\dfrac{2}{p}\;<\;0\:.
\ees
With again $q$ being the conjugate exponent of $p$, we find thus
\begin{eqnarray}\label{estI2m}
|I^m_2(x)|&\lesssim&|x|^{m}\int_{B_1(0)\setminus B_{2|x|}(0)}|y|^{a-1-m}|f(y)|\,dy\nonumber\\[1ex]
&=&|x|^{m}\int_{B_1(0)\setminus B_{2|x|}(0)}|y|^{a+1-m-\epsilon-\frac{2}{p}}|y|^{\epsilon-\frac{2}{q}}|f(y)|\,dy\nonumber\\[1ex]
&\le&2^{a+1-m-\epsilon-\frac{2}{p}}\,|x|^{a+1-\frac{2}{p}-\epsilon}\,\Big\Vert |y|^{\epsilon-\frac{2}{q}}\Big\Vert_{L^q(B_1(0))}\,\Vert f\Vert_{L^p(B_{1}(0))}\nonumber\\[1ex]
&\lesssim&2^{a+1-m-\epsilon-\frac{2}{p}}\,|x|^{1-\frac{2}{p}-\epsilon}\,|\mu(x)|\:.
\end{eqnarray}
Combining altogether in (\ref{J1}) our findings (\ref{estI1})-(\ref{estI2m}), we obtain that
\be\label{estJ1}
J_1(x)\;=\;P_1(\overline{x})\,+\,|\mu(x)|\,T_1(x)\:,
\ee
where $P_1$ is a polynomial of degree at most $a$, and the remainder $T_1$ satisfies the estimate
\be\label{estR1}
|T_1(x)|\:\,\lesssim\:\,|x|^{1-\frac{2}{p}-\epsilon}\:,\qquad\forall\:\:\epsilon>0\:.
\ee
Altogether, (\ref{estJ0}) and (\ref{estJ1}) put into (\ref{green}) show that there holds
\be\label{formw}
\nabla u(x)\;=\;P(\overline{x})\,+\,|\mu(x)|\,T(x)\:,
\ee
where $P:=P_0+P_1$ is a polynomial of degree at most $a$, and the remainder $T:=T_0+T_1$ satisfies the same estimate (\ref{estR1}) as $T_1$. 
%Hence
%\bes
%|\mu(x)|^{-1}|\nabla u(x)|\;=\;|\mu(x)|^{-1}P(\overline{x})\,+\,T(x)\:.
%\ees
%
%\noindent
%By hypothesis, the left-hand side of the latter is bounded. Moreover, by choosing $\,0<\epsilon\le1-\frac{2}{p}$, we insure that $T$ is also bounded. Accordingly, the term $\,|\mu(x)|^{-1}P(\overline{x})\in L^\infty$. As $P$ is a polynomial of degree at most $a$, and $|\mu|\simeq|x|^a$, we deduce that $P(\overline{x})=A\,\overline{x}^{\,a}$, for some constant $A$. Altogether, we thus find that
%\be\label{formw2}
%\nabla u(x)\;=\;A\,\overline{x}^{\,a}\,+\,|\mu(x)|\,T(x)\:,
%\ee
%where 
%\bes
%|T(x)|\:\lesssim\:|x|^{1-\frac{2}{p}-\epsilon}\:\qquad\forall\:\:\epsilon>0\:.
%\ees
The announced statement {(i)} ensues immediately.\\

We prove next statement {(ii)}. Comparing (\ref{stmt2}) to (\ref{formw}), we see that
\begin{eqnarray}\label{Q}
|\mu(x)|\,Q(x)&=&\nabla\big(|\mu(x)|\,T(x)\big)\\[1ex]
&=&\nabla\big(|\mu(x)|\,T_0(x)\big)\,+\,\nabla\bI_1(x)\,+\sum_{m\ge a+1}\nabla\bI_2^m(x)\,-\sum_{0\le m\le a}\nabla\bI_1^m(x)\nonumber\:.
\end{eqnarray}
By definition, 
\bes
|\mu(x)|\,T_0(x)\:=\:\sum_{m\ge a+1}C_m\,\overline{x}^{\,m}\:,
\ees
with the constants $C_m$ depending only on the $C^1$-norm of $u$ along $\partial B_1(0)$ and growing sublinearly in $m$. Using similar arguments to those leading to (\ref{estJ0}), it is  clear from (\ref{hypw}) that
\be\label{gradR0}
|\mu(x)|^{-1}\nabla\big(|\mu(x)|\,T_0(x)\big)\;\in\;L^\infty(B_1(0))\:.
\ee
Controlling the gradients of $\bI_1^m$ and $\bI_2^m$ is done {\it mutatis mutandis} the estimates (\ref{estI1a}) and (\ref{estI2m}). For the sake of brevity, we only present in details the case of $\bI_1^m$. Namely,
\begin{eqnarray}\label{alex0}
&&\hspace{-2cm}\nabla\bI^m_1(x)\;\;=\;\;\dfrac{1}{2\pi}\int_{B_1(0)\cap B_{2|x|}(0)}\nabla_x P^m(x,y)\,\mu(y)f(y)\,dy\nonumber\\[0.5ex]
&&\hspace{1cm}+\:\:\dfrac{1}{2\pi}\,\dfrac{x}{|x|}\,\otimes\int_{\partial B_{2|x|}(0)}P^m(x,y)\,\mu(y)f(y)\,dy\:.
\end{eqnarray}
After some elementary computations, and using the hypothesis $\,|\mu(y)|\simeq|y|^a$, we reach
\begin{eqnarray*}
|\nabla\bI_1^m(x)|&\lesssim&m\,|x|^{a-2}\int_{B_1(0)\cap B_{2|x|}(0)}|f(y)|\,dy\;+\;|x|^{a-1}\!\!\int_{\partial B_{2|x|}(0)}|f(y)|\,dy\nonumber\\[.5ex]
&\lesssim&m\,|x|^{a-\frac{2}{p}}\,\Vert f\Vert_{L^p(B_{1}(0))}\;+\;|x|^{a-1}\!\!\int_{\partial B_{2|x|}(0)}|f(y)|\,dy\:,
\end{eqnarray*}
so that immediately
\bes
\big\Vert |x|^{-a}\nabla\bI^m_1(x)\big\Vert_{L^{p-\epsilon}(B_1(0))}\:\,<\;\,\infty\:,\quad\qquad\forall\:\:\epsilon>0\:.
\ees
Proceeding analogously for $\nabla\bI^m_2$, we reach that for any $\epsilon>0$ there holds
\be\label{estSm}
\sum_{m\ge a+1}\big\Vert |x|^{-a}\nabla\bI^m_2(x)\big\Vert_{L^{p-\epsilon}(B_1(0))}\,+\sum_{0\le m\le a}\big\Vert |x|^{-a}\nabla\bI^m_1(x)\big\Vert_{L^{p-\epsilon}(B_1(0))}\:\:<\:\:\infty\:.
\ee
Hence, there remains only to estimate $\nabla\bI_1$. This is slightly more delicate. For notational convenience, we write
\begin{eqnarray}\label{cecondesolers}
\nabla\bI_1(x)&=&\dfrac{1}{2\pi}\,\nabla\int_{B_1(0)\cap B_{2|x|}(0)}\dfrac{x-y}{|x-y|^2}\,\mu(y)f(y)\,dy\nonumber\\[1ex]
&=:&\dfrac{1}{2\pi}\,\big(L(x)\,+\,\bK(x)\big)\:,
\end{eqnarray}
with%\footnote{for notational convenience, we write the vector $(\cos{\varphi}\,,\sin{\varphi})$ under its complex form $\text{e}^{i\varphi}$.}
\bes
\bK(x)\:=\:\chi_{B_{1/2}(0)}(x)\,\dfrac{x}{|x|}\,\otimes\int_{\partial B_{2|x|}(0)}\dfrac{x-y}{|x-y|^2}\,\mu(y)f(y)\,dy\:,
\ees
and the convolution
\bes
L(x)\:=\:\big(\Om*f(y)\,\mu(y)\,\chi_{B_1(0)\cap B_{2|x|}(0)}(y)\big)(x)\:,
\ees
where $\Om$ is the $(2\times 2)$-matrix made of the Calderon-Zygmund kernels:
\bes
\Om(z)\;:=\;\dfrac{|z|^2\,\mathbb{I}_2\,-2\,z\otimes z}{|z|^4}\:.
\ees
The boundary integral $\bK$ is easily estimated:
\bes
|x|^{-a}|\bK(x)|\;\;\lesssim\;\;\dfrac{1}{|x|}\int_{\partial B_{2|x|}(0)}|f(y)|\,dy\:,
\ees
thereby yielding
\be\label{estK}
\big\Vert |x|^{-a}\bK(x)\big\Vert_{L^p(B_1(0))}\;\;\lesssim\;\;\Vert f\Vert_{L^p(B_1(0))}\:.
\ee
To estimate $L$, we proceed as follows
\begin{eqnarray}\label{nonnon}
&&\hspace{-2cm}L(x)\,-\,\mu(x)\big(\Om*f\,\chi_{B_1(0)\cap B_{2|x|}(0)}\big)(x)\nonumber\\[1ex]
&&\hspace{1cm}=\;\;\int_{B_1(0)\cap B_{2|x|}(0)}\Om(x-y)\,f(y)\big(\mu(y)-\mu(x)\big)\,dy\:.
%&&\hspace{1cm}=\;\;\int_{B_{2|x|(0)}}\dfrac{|x-y|^2\,\mathbb{I}_2\,-\,2(x-y)\otimes(x-y)}{|x-y|^4}\,\big(w(y)-w(x)\big)f(y)\,dy\nonumber\\[1ex]
%&&\hspace{1cm}\lesssim\;\;|x|^{a-1}\int_{B_{2|x|}(0)}\dfrac{|f(y)|}{|x-y|}\,dy\;\;\lesssim\;\;|x|^{a}M_0(f)\:,
\end{eqnarray}
Let $S_x$ be the cone with apex the point $x/2$ and such that the disk $B_{|x|/4}(0)$ is inscribed in it. Note that for $y\in S_x$, there holds\, $2|x-y|>|x|$. Hence, we find
\begin{eqnarray}\label{nigo}
&&\hspace{-2cm}\int_{S_x\cap B_1(0)\cap B_{2|x|}(0)}\Om(x-y)f(y)\big(\mu(y)-\mu(x)\big)\,dy\nonumber\\[0ex]
&&\hspace{3cm}\;\lesssim\;|\mu(x)|\,|x|^{-2}\!\int_{B_{2|x|}(0)}|f(y)|\,dy\:.
\end{eqnarray}
%\begin{eqnarray*}\label{nigo}
%&&\hspace{-.5cm}\int_{S_x\cap B_1(0)\cap B_{2|x|}(0)}\Om(x-y)\,f(y)\big(\mu(y)-\mu(x)\big)\,dy\nonumber\\[1ex]
%&&\hspace{2cm}\lesssim\;\;|\mu(x)|\,|x|^{-2}\int_{B_{2|x|}(0)}|f(y)|\,dy\nonumber\\[1ex]
%&&\hspace{2cm}\lesssim\;\;|x|^{a-1}\int_{B_{2|x|}(0)}|y|^{-1}|f(y)|\,dy\:\:\,\lesssim\:\:|\mu(x)|\,|x|^{-\frac{2}{p}}\:,\qquad
%\end{eqnarray*}

\noindent
\smallskip
\noindent
By hypothesis, the function $\mu$ is continuously differentiable away from the origin. Thus, to each point $y$ in the complement of the cone $S_x$, there corresponds some $\al\equiv\al(x,y)\in[0,1]$ with
\bes
\mu(y)-\mu(x)\;=\;(x-y)\cdot\nabla\mu\big(\al\,x+(1-\al)y\big)\:.
\ees
Using (\ref{hypw2}), we deduce easily
\bes
|\mu(y)-\mu(x)|\;\lesssim\;|x|^{a-1}|x-y|\qquad\:\:\forall\:\:y\in S_x^c\cap B_1(0)\cap B_{2|x|}(0)\:.
\ees
Accordingly, there holds
\begin{eqnarray}\label{douille}
&&\hspace{-2cm}\int_{S^c_x\cap B_1(0)\cap B_{2|x|}(0)}\Om(x-y)\,f(y)\big(\mu(y)-\mu(x)\big)\,dy\nonumber\\[0ex]
&&\hspace{1cm}\lesssim\;\;|x|^{a-1}\!\!\int_{B_{2|x|}(0)}\dfrac{|f(y)|}{|x-y|}\,dy\:\:\,\lesssim\:\:\,|\mu(x)|\,M_0f(x)\:,
\end{eqnarray}
where we have used the same estimate as in (\ref{estI1}). Bringing (\ref{nigo}) and (\ref{douille}) into (\ref{nonnon}) and using the fact that $|\mu(x)|\simeq|x|^a$ yields
\begin{eqnarray*}
|\mu(x)|^{-1}|L(x)|&\lesssim&\big(\Om*f(y)\,\chi_{B_1(0)\cap B_{2|x|}(0)}(y)\big)(x)\nonumber\\[1ex]
&&\hspace{1.5cm}+\:\;\dfrac{1}{|x|^2}\int_{B_{2|x|}(0)}|f(y)|\,dy\;+\,M_0f(x)\:.
\end{eqnarray*}
Because $f\in L^p$, standard estimates on Calderon-Zygmund operators, on the maximal function, and a classical Hardy inequality then give us
\bes
\big\Vert |\mu|^{-1}L\big\Vert_{L^{p}(B_1(0))}\:\lesssim\:\Vert f\Vert_{L^p(B_1(0))}\:<\:\infty\:.
\ees
Owing to the latter and to (\ref{estK}), we obtain from (\ref{cecondesolers}) that
$\,|\mu|^{-1}\nabla\bI_1\in L^{p}(B_1(0))$.
With (\ref{gradR0}) and (\ref{estSm}), the identity (\ref{Q}) thus implies that $Q$ belongs to $L^{p-\epsilon}$ for all $\epsilon>0$. This completes the first part of statement (ii).\\

We shall now prove the second part of (ii), and show that the trace of $Q$ is in $L^p$.  To this end, let us note that
\be\label{gizmo}
\text{Tr}\,\nabla\overline{x}\;=\;\text{Tr}\left(\begin{array}{cc}1&0\\0&-1\end{array}\right)\;=\;0\:.
\ee
We have seen in (\ref{Q}) that
\be\label{aalex}
|\mu|\,Q\;=\;\nabla\big(|\mu|\,T_0\big)\,+\,\nabla\bI_1\,+\sum_{m\ge a+1}\nabla\bI_2^m\,-\sum_{0\le m\le a}\nabla\bI_1^m\:.
\ee
By definition, $|\mu(x)|\,T_0(x)=\sum_{m\ge a+1}C_m\,\overline{x}^{\,m}$, so that (\ref{gizmo}) gives
\be\label{alex1}
\text{Tr}\;\nabla\big(|\mu(x)|\,T_0(x)\big)\;=\;0\:.
\ee
Owing to the fact that $P^m(x,y)=\overline{x}^{\,m}\,\overline{y}^{\,-(m+1)}$, it then easily follows from (\ref{gizmo}) and (\ref{alex0}) that
\bes
\text{Tr}\;\nabla I^m_1(x)\;=\;\dfrac{1}{2\pi}\,\text{Tr}\;\dfrac{x}{|x|}\,\otimes\int_{\partial B_{2|x|}(0)}P^m(x,y)\,\mu(y)f(y)\,dy\:;
\ees
whence the estimate
\bes
|\mu(x)|^{-1}\big|\text{Tr}\;\nabla I^m_1(x)\big|\;\lesssim\;2^{a-m-1}\dfrac{1}{|x|}\int_{\partial B_{2|x|}(0)}|f(y)|\,dy\:,
\ees
and thus
\be\label{alex2}
\big\Vert|\mu|^{-1}\text{Tr}\,\nabla I^m_1\big\Vert_{L^p}\;\lesssim\;2^{a-m-1}\,\Vert f\Vert_{L^p}\:.
\ee
In exactly the same fashion, one finds
\be\label{alex3}
\big\Vert|\mu|^{-1}\text{Tr}\,\nabla I^m_2\big\Vert_{L^p}\;\lesssim\;2^{a-m-1}\,\Vert f\Vert_{L^p}\:.
\ee
%Accordingly, there holds
%\be\label{estSmspe}
%\sum_{m\ge a+1}\big\Vert\nabla\bI^m_2(x)\big\Vert_{L^{\infty}(B_1(0))}\,+\sum_{0\le m\le a}\big\Vert\nabla\bI^m_1(x)\big\Vert_{L^{\infty}(B_1(0))}\:\:\lesssim\:\:|\mu(x)|\:.
%\ee
Finally, there remains to handle the term $|\mu|^{-1}\text{Tr}\,\nabla I_1$. But this term belongs to $L^p$, as we have shown that $|\mu|^{-1}\nabla I_1$ does. Combining this altogether with (\ref{alex1})-(\ref{alex3}) into (\ref{aalex}) yields the announced result.\\[-3ex]

$\hfill\blacksquare$ \\[-3ex]

\begin{Co}\label{CZcoro}
Let $u\in C^2(B_1(0)\setminus\{0\})$ solve
\bes\label{equu}
\Delta u(x)\;=\;\mu(x) f(x)\qquad\text{in\:\:\:} B_1(0)\:,
\ees
where
\bes\label{hypuu}
|f(x)|\;\lesssim\;|x|^{n+r}\qquad\text{and}\qquad|\mu(x)|\;\simeq\;|x|^{a}\:,
\ees
for two non-negative natural numbers $n$ and $a$ ; and $r\in(0,1)$. \\
Then
\be\label{concuu}
\nabla u(x)\;=\;P(\overline{x})\,+\,|\mu(x)|\,T(x)\:,
\ee
where $P$ is a complex-valued polynomial of degree at most $(a+n+1)$, and near the origin $T(x)=\text{O}(|x|^{n+1+r-\epsilon})$ for every $\epsilon>0$.\\[1ex]
If in addition $\mu$ satisfies (\ref{hypw2}), then $\,|x|^{-(n+r)}|\mu|^{-1}\nabla\big(|\mu|T\big)$ belongs to $L^p$ for all finite $p$. Furthermore, there holds the estimate
\be\label{conctrace}
\big|\text{Tr}\;\nabla\big(|\mu(x)|T(x)\big)\big|\;\;\lesssim\;\;|x|^{n+r}|\mu(x)|\:.
\ee
\end{Co}
$\textbf{Proof.}$ The argument goes along the same lines as that of Proposition \ref{CZpondere}. We set
\bes
\om(x)\;:=\;|x|^{n+r}\mu(x)\qquad\text{and}\qquad h(x)\;:=\;|x|^{-(n+r)}f(x)\:.
\ees
From the given hypotheses, we see that $h\in L^\infty$, and $\om$ satisfies (\ref{hypw}) with $(a+n+r)$ in place of $a$. If $\mu$ satisfies (\ref{hypw2}), then so does $\om$, again with $(a+n+r)$ in place of $a$. \\
Using the representation (\ref{green}) gives
\begin{eqnarray*}\label{greenuu}
\nabla u(x)&=&\dfrac{1}{2\pi}\,\int_{\partial B_1(0)}\bigg[\dfrac{x-y}{|x-y|^2}\,\partial_{\vec{\nu}}\,u(y)\,-\,u(y)\,\partial_{\vec{\nu}}\,\dfrac{x-y}{|x-y|^2}\bigg]\,d\sigma(y)\nonumber\\[1ex]
&&\hspace{1cm}-\:\:\dfrac{1}{2\pi}\int_{B_1(0)}\dfrac{x-y}{|x-y|^2}\,\om(y)\,h(y)\,dy\nonumber\\[1.2ex]
&=:&J_0(x)\,+\,J_1(x)\:,\hspace{2.25cm}\forall\:\:x\in B_1(0)\:,
\end{eqnarray*}
where $\vec{\nu}$ is the outer normal unit-vector to the boundary of $B_1(0)$. \\
The integral $J_0$ is estimated as in (\ref{estJ0}) so as to yield
\bes\label{estJ0uu}
J_0(x)\:=\:P_0(\overline{x})\,+\,|\mu(x)|\,T_0(x)\:,
\ees
where $P_0$ is a polynomial of degree at most $(a+n+1)$, and $T_0(x)=\text{O}(|x|^{n+2})$ with $|\mu|^{-1}\nabla\big(|\mu|T_0\big)=\text{O}(|x|^{n+1})$. \\[1ex]
We next estimate the integral $J_1$. We proceed again as in the proof of Proposition \ref{CZpondere}. Namely,
\bes\label{J1uu}
J_1(x)\;=\;I_1(x)\,+\sum_{m=a+n+2}^\infty\!\!I_2^m(x)\,-\sum_{m=0}^{a+n+1}\!\!I_1^m(x)\,+\sum_{m=0}^{a+n+1}\!\!I_1^m(x)+I_2^m(x)\:,
\ees
where we have put
\bes
I_1(x)\,:=\,\dfrac{1}{2\pi}\int_{B_1(0)\cap B_{2|x|}(0)}\dfrac{x-y}{|x-y|^2}\,\om(y)\,h(y)\,dy\:,
\ees
\bes
I^m_1(x)\,:=\,\dfrac{1}{2\pi}\int_{B_1(0)\cap B_{2|x|}(0)}P^m(x,y)\,\om(y)\,h(y)\,dy\:,
\ees
\bes
I^m_2(x)\,:=\,\dfrac{1}{2\pi}\int_{B_1(0)\setminus B_{2|x|}(0)}P^m(x,y)\,\om(y)\,h(y)\,dy\:.
\ees\\
As before, $\,P^m(x,y):=\overline{x}^{\,m}\,\overline{y}^{\,-(m+1)}$.
We first observe that the last sum in the expression for $J_1$ may be written
\begin{eqnarray*}\label{polyuu}
P_1(x)&:=&\sum_{0\le m\le a+n+1}I_1^m(x)+I_2^m(x)\nonumber\\
&=&\!\!\!\sum_{0\le m\le a+n+1}\,\int_{B_1(0)}P^m(x,y)\,\om(y)\,h(y)\,dy\:\:=\sum_{0\le m\le a+n+1}A_m\,\overline{x}^{\,m}\:,
\end{eqnarray*}
where
\bes
A_m\;:=\;\int_{B_1(0)}\overline{y}^{\,-(m+1)}\,\om(y)\,h(y)\,dy\:.
\ees\\
From the boundedness of $h$ and the hypothesis $|\om(y)|\simeq|y|^{a+n+r}$, it follows easily that $|A_m|<\infty$ for $m<a+n+1+r$, and thus since $r>0$, that $P_1$ is a polynomial of degree at most $(a+n+1)$.
Once this has been observed, the remainder of the proof is found {\it mutatis mutandis} that of Proposition \ref{CZpondere}. Namely, we write 
\bes
J_1(x)\;=\;P_1(\overline{x})\,+\,|\om(x)|\,T_1(x)\:,
\ees
with $T_1(x)=\text{O}(|x|^{1-\epsilon})$ for all $\epsilon>0$. Moreover, and $|\om|^{-1}\nabla \big(|\om|T_1\big)\in L^{p}$ for all $p<\infty$ ; and $\,|\om|^{-1}\text{Tr}\,\nabla \big(|\om| T_1\big)\in L^{\infty}$.
 \\[1ex]
Finally, setting $P=P_0+P_1$ and $T=T_0+|x|^{n+r\,}T_1=\text{O}(|x|^{n+r+1-\epsilon})$ gives the desired representation (\ref{concuu}). Clearly, from (\ref{hypw2}) and the above, there holds 
\bes
\big||\mu|^{-1}\nabla\big(|\mu|T\big)\big|\;\lesssim\;\big||\mu|^{-1}\nabla\big(|\mu|T_0\big)\big|\,+\,|x|^{n+r}\big||\om|^{-1}\nabla\big(|\om|T_1\big)\big|
%\big|\nabla T_0(x)\big|\,+\,|x|^{n+r}|\om|^{-1}\big|\nabla\big(|\om|T_1(x)\big)\big|\,+\,|x|^{n+r-1}|T_1(x)|\:,
\ees
%\bes
%\big|\nabla T(x)\big|\;\lesssim\;\big|\nabla T_0(x)\big|\,+\,|x|^{n+r}|\om|^{-1}\big|\nabla\big(|\om|T_1(x)\big)\big|\,+\,|x|^{n+r-1}|T_1(x)|\:,
%\ees
so that indeed $\,|x|^{-(n+r)}|\mu|^{-1}\nabla\big(|\mu|T\big)$ belongs to $L^p$ for all finite $p$. Furthermore, we have
\begin{eqnarray*}
\big|\text{Tr}\,\nabla\big(|\mu|T\big)\big|&\leq&\big|\text{Tr}\,\nabla\big(|\mu|T_0\big)\big|\,+\,\big|\text{Tr}\,\nabla\big(|\om|T_1\big)\big|\\[1ex]
&\lesssim&|x|^{n+1}|\mu|\,+\,|\om|\;\;\lesssim\;\;|x|^{n+r}|\mu|\:,
\end{eqnarray*}
as announced. \\[-3ex]

$\hfill\blacksquare$ \\

We may further iterate the previous result to obtain the next one. 

\begin{Co}\label{CZcoro2}
Let $u\in C^\infty(B_1(0)\setminus\{0\})$ solve
\bes
\Delta u(x)\;=\;\mu(x) f(x)\qquad\text{in\:\:\:} B_1(0)\:,
\ees
where $\,|\mu(x)|\simeq|x|^a$, for some $a\in\mathbb{N}^*$. In addition, we assume that
\bes
|\nabla^jf(x)|\;\lesssim\;|x|^{n-j+r}\qquad\text{and}\qquad|\nabla^j\mu(x)|\;\lesssim\;|x|^{a-j}\:,
\ees
for some $n\in\mathbb{N}$ and $r\in(0,1)$, and for all $j$ satisfying
\bes
0\;\le\;j\;\le\;J\;\le\;\min\,\{a\,,n+1\}\:,\qquad\text{for some}\:\:\:J\in\mathbb{N}^*\:.
\ees
Then there holds for all $j\le J$\,:
\bes
\nabla^{j+1} u(x)\;=\;\nabla^{j}P(\overline{x})\,+\,|\mu(x)|V_j(x)\:,
\ees
where $P$ is a two-component real-valued polynomial of degree at most $(a+n+1)$, and near the origin $V_j(x)=\text{O}(|x|^{n+r-j+1-\epsilon})$ for every $\epsilon>0$.\\[.5ex]
Furthermore\footnote{note that $\,|\mu|V_{j+1}=\nabla\big(|\mu|V_j\big)$},
\bes
|x|^{-(a+n+r-j)}\,\nabla\big(|\mu(x)|V_j(x)\big)\:\in\:\bigcap_{p<\infty}L^p(B_1(0))\:.
\ees
and
\bes
\big|\text{Tr}\;\nabla\big(|\mu(x)|V_j(x)\big)\big|\;\;\lesssim\;\;|x|^{a+n+r-j}\:.
\ees
\end{Co}

%$\log(1-2x\cos\varphi+x^2)\;=\;-2\sum_{k\ge 1}\dfrac{\cos(k\varphi)}{k}\,x^k\qquad,\qquad |x|\le1\quad x\cos\varphi\ne1$

\eject


\begin{thebibliography}{99}
\bibitem[Ad]{Ad} Adams, Robert ``A note on Riesz potentials." Duke Math. J. 42 (1975), no. 4, 765-778.
\bibitem[Al]{Al} Almeida, Lu\'is ``The regularity problem for generalized harmonic maps into homogeneous spaces." Calc. Var. 3 (1995), 193-242.
\bibitem[BR1]{BR1} Bernard, Yann ; Rivi\`ere Tristan ``Local Palais-Smale sequences for the Willmore functional" arXiv:0904.0360v1 (2009). 
\bibitem[BR2]{BR2} Bernard, Yann ; Rivi\`ere Tristan ``Energy quantization for Willmore surfaces and applications" arXiv:1106.3780v1 (2011). 
%\bibitem[Bla]{Bla} Blaschke, Wilhelm ``Vorlesungen \"Uber Differential Geometrie III." Springer (1929).
%\bibitem[Che]{Che} Chen, Bang-Yen ``Some conformal invariants of submanifolds and their applications."  Boll. Un. Mat. Ital. (4) 10 (1974), 380--385.
\bibitem[DL]{DL} Dautray, Robert ; Lions, Jacques-Louis ``Mathematical Analysis and Numerical Methods for Science and Technology ; volume 3: Spectral Theory and Applications." Springer-Verlag (1992).
%\bibitem[Gi]{Gi} Giaquinta, Mariano ``Multiple Integrals in the Calculus of Variations and Nonlinear Elliptic Systems." Annals of Mathematics Studies, 105, PUP (1983).
\bibitem[GW]{GW} Gr\"uter, Michael ; Widman Kjell-Ove ``The Green function for uniformly elliptic equations." Manuscripta Math. 37 (1982), 303-342.
\bibitem[He]{He} H\'elein, Fr\'ed\'eric ``Harmonic Maps, Conservation Laws, and Moving Frames." Cambridge Tracts in Mathematics, 150, CUP (2002).
\bibitem[Hu]{Hu} Huber, Alfred  ``On subharmonic functions and differential geometry in the large."  Comment. Math. Helv. 32 (1957) 13--72.
\bibitem[KS1]{KS1} Kuwert, Ernst; Sch\"atzle Reiner ``Removability of point singularities of Willmore surfaces." Ann. Math. 160 (2004), 315--357.
\bibitem[KS2]{KS2} Kuwert, Ernst; Sch\"atzle Reiner ``Branch points of Willmore surfaces." Duke Math. J. 138 (2007), no. 2, 179--201.
\bibitem[KS3]{KS3} Kuwert, Ernst; Sch\"atzle Reiner ``The Willmore flow with small initial energy." J. Diff. Geom. 57 (2001), 409--441.
\bibitem[MS]{MS}  M\"uller, Stefan;  \v{S}ver\'ak, Vladim\'ir ``On surfaces of finite total curvature." J. Diff. Geom. 42 (1995), no. 2, 229--258.
\bibitem[Ri1]{Ri1} Rivi\`ere, Tristan ``Analysis aspects of the Willmore functional." Invent. Math. 174 (2008), no. 1, 1-45.
\bibitem[Ri2]{Ri2} Rivi\`ere, Tristan  ``Variational principles for immersed surfaces with $L^2$-bounded second fundamental form'' arXiv:1007.2997 (2010).
\bibitem[Ri3]{Ri3} Rivi\`ere, Tristan ''Conformally Invariant 2-dimensional Variational Problems'' Cours joint de l'Institut Henri Poincar\'e - Paris XII Creteil, Novembre 2010.
\bibitem[Ta]{Ta} Tartar, Luc ``An Introduction to Sobolev Spaces and Interpolation Spaces." Lectures notes of the Unione Matematica Italiana, no. 3 (2007).
%\bibitem[Ta3]{Ta3} Tartar, Luc ``Imbedding theorems of Sobolev spaces into Lorentz spaces." 
%Boll. Unione Mat. Ital. Sez. B Artic. Ric. Mat. (8) 1 (1998), no. 3, 479--500. 
%\bibitem[To1]{To1} Toro, Tatiana ``Surfaces with generalized second fundamental form in $L^2$ are Lipschitz manifolds." J. Diff. Geom. 39 (1994), 65-101.
\bibitem[To]{To} Toro, Tatiana ``Geometric conditions and existence of bilipschitz parametrisations." Duke Math. J. 77 (1995), no. 1, 193-227.
%\bibitem[Wen]{Wen} Wente, Henry C. ``An existence theorem for surfaces of constant mean \bibitem[Wil1]{Wil1} Willmore, T. J. ``Note on embedded surfaces." Ann. Stiint. Univ. "Al. I. Cuza" Iasi. Sect. I a Mat. (N.S.) 11B 1965 493--496.
%\bibitem[Wil2]{Wil2} Willmore, T. J. ``Riemannian Geometry." Oxford University Press (1997).
\bibitem[We]{We} Weiner, Joel ``On a problem of Chen, Willmore, {\it et al}." Indiana U. Math. J. 27 (1978) no. 1, 19--35.
\bibitem[Zi]{Zi} Ziemer, William ``Weakly Differentiable Functions." GSM , Springer-Verlag (1989).
\end{thebibliography}
\end{document}